\documentclass[12pt]{amsart}
\usepackage[T1]{fontenc}
\usepackage[utf8]{inputenc}
\usepackage[english]{babel}
\usepackage{lmodern,amsmath,amssymb,amsthm,mathtools}
\usepackage{microtype}
\usepackage{mathrsfs,esint,enumitem,booktabs,longtable,array,aliascnt}
\usepackage[colorlinks=true,linkcolor=blue!55!black,citecolor=green!40!black,urlcolor=blue!55!black]{hyperref}
\usepackage{xcolor}
\usepackage[nameinlink,noabbrev]{cleveref}
\numberwithin{equation}{section}

\newtheorem{introtheorem}{Theorem}

\newtheorem{restatedtheorem}{Theorem}

\crefname{introtheorem}{Theorem}{Theorems}
\Crefname{introtheorem}{Theorem}{Theorems}
\crefname{restatedtheorem}{Theorem}{Theorems}
\Crefname{restatedtheorem}{Theorem}{Theorems}

\newaliascnt{proposition}{theorem}
\newtheorem{proposition}[proposition]{Proposition}
\aliascntresetthe{proposition}
\newaliascnt{lemma}{theorem}
\newtheorem{lemma}[lemma]{Lemma}
\aliascntresetthe{lemma}
\newaliascnt{corollary}{theorem}
\newtheorem{corollary}[corollary]{Corollary}
\aliascntresetthe{corollary}
\theoremstyle{definition}
\newaliascnt{definition}{theorem}
\newtheorem{definition}[definition]{Definition}
\aliascntresetthe{definition}
\theoremstyle{remark}
\newaliascnt{remark}{theorem}

\aliascntresetthe{remark}
\newcommand{\R}{\mathbb{R}}
\newcommand{\Sph}{\mathbb{S}}
\newcommand{\eps}{\varepsilon}
\newcommand{\dd}{\,\mathrm{d}}

\allowdisplaybreaks[2]
\title[Singular Escobar--Yamabe asymptotics]{Delaunay Asymptotics for the Singular Escobar--Yamabe Equation}
\author{Felipe G. Bogo}
\address{Instituto de Matem\'atica e Estat\'istica, Universidade Federal
Fluminense, Niter\'oi, RJ, Brasil, 24210-201.}
\email{felipebogo@id.uff.br}
\hypersetup{pdfauthor={Felipe G. Bogo}}
\date{}
\begin{document}
\raggedbottom
\begin{abstract}
In dimensions $n\geq3$, we prove that near an isolated boundary
singularity a positive solution of the Escobar--Yamabe equation
for a smooth metric approaches a fixed periodic Delaunay solution
or the cylindrical solution if its trace lies between positive
multiples of $|\overline x|^{-(n-2)/2}$, where $\overline x$ is the
tangential Fermi coordinate. The relative error decays at a positive
power rate, with weighted estimates through two derivatives in Fermi
coordinates. The necksize is unique, and the translation parameter
is unique modulo the least positive period in the periodic case.
The Hamiltonian smoothly parametrizes the periodic family modulo
translation. We classify Jacobi fields of polynomial growth at both
ends without assuming symmetry. Translation and energy differentiation
generate these fields along periodic solutions, and two oscillatory
fields generate them along the cylindrical solution. We also exhibit
a first-order obstruction to improving these asymptotics by
boundary-preserving conformal deformations of the flat models.
\end{abstract}
\maketitle

\setcounter{tocdepth}{1}
\tableofcontents

\section{Introduction}\label{sec:introduction}

Let $(M,g)$ be a smooth Riemannian manifold with boundary of
dimension $n\geq3$, and fix $p\in\partial M$. Choose boundary
Fermi coordinates $x=(\overline x,x_n)$ centered at $p$, with
$x_n$ increasing inward. For sufficiently small $\rho>0$, these
coordinates identify a neighborhood of $p$ with the half-ball
$B_\rho^+=B_\rho(0)\cap\R^n_+$, where
$\R^n_+=\{x\in\R^n:x_n\geq0\}$ and $B_\rho(0)$ is the open
Euclidean ball of radius $\rho$. We write
$\partial_0B_\rho^+=B_\rho^+\cap\partial\R^n_+$. The coordinate
conventions and the expansion of $g$ are given in
Section~\ref{subsec:notation-fermi}, following
Almaraz~\cite[Section~2.3]{Almaraz} and Marques~\cite{MarquesFermi}.

For a positive function $U\in C^\infty(B_\rho^+\setminus\{0\})$,
write $u(\overline x)=U(\overline x,0)$ for its boundary trace
on $\partial_0B_\rho^+\setminus\{0\}$. We study the singular
Escobar--Yamabe problem
\begingroup
\predisplaypenalty=0\relax
\begin{equation}\label{eq:intro-geometric}
 \begin{cases}
 -\Delta_g U+\dfrac{n-2}{4(n-1)}R_gU=0
   & \text{in }\operatorname{int}B_\rho^+,\\
 \partial_{\nu_g}U+\dfrac{n-2}{2}h_g U
   =\dfrac{n-2}{2}u^{n/(n-2)}
   & \text{on }\partial_0B_\rho^+\setminus\{0\}.
 \end{cases}
\end{equation}
\endgroup
Here $R_g$ is the scalar curvature, $\nu_g$ is the outward
unit normal, and $h_g$ is the average of the principal
curvatures with respect to $\nu_g$. By conformal
covariance~\cite{Escobar}, the metric $U^{4/(n-2)}g$ is
scalar flat and has boundary mean curvature $1$ on
$\partial_0B_\rho^+\setminus\{0\}$. We call the singularity
at $0$ removable if $u$ extends continuously across $0$,
and nonremovable otherwise.

For the Yamabe equation without boundary, Caffarelli, Gidas
and Spruck~\cite{CGS} established radial symmetry for positive
singular solutions in punctured Euclidean space and their
asymptotic behavior near an isolated singularity. The resulting
periodic solutions in Emden--Fowler coordinates are the
Fowler solutions. Korevaar, Mazzeo, Pacard and Schoen~\cite{KMPS}
obtained refined asymptotics using conformal deformations of
these solutions. For a smooth metric, Marques~\cite{Marques}
proved Fowler asymptotics in dimensions $3\leq n\leq5$, and
Xiong and Zhang~\cite{XiongZhang} treated dimension $6$.
Han, Jiang, Wang, Xiong and Zhang~\cite{HJWXZ} extended the
classification of isolated singularities for smooth metrics
to dimensions $7\leq n\leq29$.

We refer to~\eqref{eq:intro-geometric} with the Euclidean
metric as the flat Escobar--Yamabe problem. A function $U$ on
$\R^n_+\setminus\{0\}$ is \emph{cylindrically symmetric} if
$U(\overline x,x_n)=U(|\overline x|,x_n)$. The extension
formulation of Caffarelli and Silvestre~\cite{CS} relates the
flat problem to the critical equation of order $1/2$ on
$\R^{n-1}$. Caffarelli, Jin, Sire and
Xiong~\cite[Theorems~1.1 and~1.3]{CJXS} proved the critical
two-sided estimates near a nonremovable singularity and
cylindrical symmetry for positive global solutions with a
nonremovable singularity. DelaTorre and Gonz\'alez~\cite{DG}
studied the cylindrical equation and its Hamiltonian.
Periodic solutions of the fractional equation were constructed
by DelaTorre, del Pino, Gonz\'alez and Wei~\cite{DDPGW}.

To describe the flat solutions used below, write $x=r\theta$,
where $r=|x|>0$ and $\theta\in\Sph^{n-1}_+$, and set $t=-\log r$.
The normalized Emden--Fowler transform is
$v(t,\theta)=b_n^{-1}r^{(n-2)/2}U(r\theta)$, with $b_n$ and
$c$ as in~\eqref{eq:prelim-normalizing-constants}. For a
function $v$ on $\R\times\Sph^{n-1}_+$, its trace is
$\tilde v(t,\vartheta)=v(t,(\vartheta,0))$ for
$\vartheta\in\Sph^{n-2}$. If $U$ is cylindrically symmetric,
$v$ depends only on $t$ and $\beta=\arccos\theta_n$. In this
case, we write $\tilde v(t)=v(t,L)$, where $L=\pi/2$.
The normalized flat equation is
\begin{equation}\label{eq:intro-normalized-model}
 \begin{cases}
 v_{tt}+\Delta_{\Sph^{n-1}_+}v-\dfrac{(n-2)^2}{4}v=0
   & \text{in }\R\times\operatorname{int}\Sph^{n-1}_+,\\
 \partial_{\nu_{\Sph^{n-1}_+}}v=c\,\tilde v^{n/(n-2)}
   & \text{on }\R\times\partial\Sph^{n-1}_+,
 \end{cases}
\end{equation}
where $\nu_{\Sph^{n-1}_+}$ is the outward unit conormal
to the equator.

A \emph{Delaunay solution} is a positive cylindrically symmetric
solution of the flat Escobar--Yamabe problem on
$\R^n_+\setminus\{0\}$ with a nonremovable singularity at $0$
whose Emden--Fowler transform is periodic and nonconstant in $t$.
The \emph{necksize} is the minimum of the normalized trace on
$\R$. Theorem~A of~\cite{PartI} gives, for each $\eps\in(0,1)$,
a unique normalized Delaunay solution $v_\eps$ with
$\tilde v_\eps(0)=\min_{t\in\R}\tilde v_\eps(t)=\eps$.
We denote its least positive period by $T_\eps$.

A solution of~\eqref{eq:intro-normalized-model} is
\emph{cylindrical} if it is independent of $t$. We denote the
positive cylindrical solution by $C_n(\theta)$ and set
$v_1(t,\theta)=C_n(\theta)$. The boundary trace of $C_n$ is
identically $1$.

Theorem~B of~\cite{PartI} classifies all positive global solutions
of the flat Escobar--Yamabe problem with a possible isolated
boundary singularity. In
the nonremovable case, the normalized transform is either
$C_n$ or a translate $v_\eps(t+\tau,\theta)$.

For a smooth background metric in dimensions $3$ and $4$,
Liao and Ma~\cite[Theorems~1.1 and~1.2]{LiaoMa} obtained
the critical upper and lower bounds near a nonremovable
boundary singularity. They also showed that $U(x)$ approaches
its average over $\{(\overline y,x_n):|\overline y|=|\overline x|\}$,
with relative error $O(|x|^\alpha)$ for some $\alpha\in(0,1)$.
In their preprint of 5 November 2026,
Jin, Liao and Xiong~\cite[Theorem~1.4]{JinLiaoXiong} proved
that a positive local solution of the flat Escobar--Yamabe
problem with a nonremovable boundary singularity approaches
a unique global singular solution with a relative error
decaying at a positive power rate. Under the critical
two-sided bounds on $u$, we prove convergence to a fixed
flat solution in boundary Fermi coordinates for a general
smooth metric in dimensions $n\geq3$. In Emden--Fowler
coordinates, the flat solution is either a periodic member
of the family or $C_n$. The estimates include the derivatives
of orders $1$ and $2$.

\subsection*{Main results}
\label{subsec:intro-main-results}

For a solution $v$ of~\eqref{eq:intro-normalized-model}, the
Hamiltonian in~\eqref{eq:proof-hamiltonian} is given by
\[
 \begin{aligned}
 |\Sph^{n-2}|H(v)(t)
 &=\frac12\int_{\Sph^{n-1}_+}
    \left(v_t^2-|\nabla_{\Sph^{n-1}_+}v|^2
                 -\frac{(n-2)^2}{4}v^2\right)\dd\theta\\
 &\quad+\frac{(n-2)c}{2(n-1)}
    \int_{\partial\Sph^{n-1}_+}
       \tilde v^{\frac{2(n-1)}{n-2}}\,\dd\sigma .
 \end{aligned}
\]
Both integrals are evaluated at $t$. By
\eqref{eq:intro-normalized-model}, $H(v)(t)$ is independent of
$t$, and we write $H(v)$ for its value. For the cylindrical
solution, $H_c=H(C_n)=-c/(2(n-1))$.

\begin{introtheorem}\label{intro:energy}
The map
\[
 (0,1)\longrightarrow(H_c,0),\qquad
 \eps\longmapsto H(v_\eps)
\]
is a smooth diffeomorphism with strictly negative derivative.
Two normalized noncylindrical Delaunay solutions have the same
Hamiltonian if and only if they differ by a translation in $t$.
The unique normalized singular solution with Hamiltonian $H_c$
is $C_n$. If $\eps(h)$ is the inverse of this map, then
$v_h=v_{\eps(h)}$ has least positive period $T(h)=T_{\eps(h)}$.
Let $v_h^{(1)}(s,\theta)=v_h(T(h)s,\theta)$. For every integer
$k\geq0$ and $0<\gamma<1$, the map $h\mapsto(T(h),v_h^{(1)})$
is smooth from $(H_c,0)$ into
$(0,\infty)\times C^{k,\gamma}_{1,\mathrm{sym,even}}
(\R\times\Sph^{n-1}_+)$.
\end{introtheorem}

\noindent\emph{Idea of the proof.}
The derivative of $H(v_\eps)$ in $\eps$ is the pairing
$\Omega(\partial_t v_\eps,\partial_\eps v_\eps)$
in~\eqref{eq:s4-symplectic-pairing}. Proposition~3.1(iii)
of~\cite{PartI} determines an odd cylindrically symmetric
Jacobi field from its $t$-derivative at $(0,L)$. If the pairing
were zero, the compatibility condition for a periodic
linear problem would hold. A solution of that problem
would produce an odd Jacobi field with polynomial growth
that changes under translation by $T_\eps$. This contradicts
the fact that every such odd field is a multiple of
$\partial_t v_\eps$. The limits of $H(v_\eps)$ at $\eps=0$
and $\eps=1$ give the sign of its derivative and its range.
The proof is given in Section~\ref{sec:energy-coordinate}.

A Jacobi field along $v$ is a smooth solution $\Phi$ of
the linearization of~\eqref{eq:intro-normalized-model}, given
in~\eqref{eq:jacobi-linearized-equation}. The space
$\mathcal K_{\mathrm{pol}}(v)$ consists of these fields on
$\R\times\Sph^{n-1}_+$ satisfying the polynomial growth
condition~\eqref{eq:jacobi-pointwise-growth} for some $q\geq0$.
The bound is required both as $t\to+\infty$ and as $t\to-\infty$.
We write $\mathcal K_{\mathrm{bd}}(v)$ for the space of
Jacobi fields bounded on $\R\times\Sph^{n-1}_+$.

\begin{introtheorem}\label{intro:polynomial-kernel}
Let $H_c<h<0$, and set $\Psi_h=\partial_h v_h$, where
the derivative is taken at fixed $t$. Then
\[
 \mathcal K_{\mathrm{pol}}(v_h)
 =\operatorname{span}\{\partial_t v_h,\Psi_h\}.
\]
For $\Phi=a\partial_t v_h+d\Psi_h$ with $a,d\in\R$, the
derivative of the Hamiltonian is $D H_{v_h}[\Phi](t)=d$
for every $t\in\R$. In particular, this derivative vanishes
if and only if $\Phi$ is a multiple of $\partial_t v_h$.
\end{introtheorem}

\noindent\emph{Idea of the proof.}
Differentiating tangential translations and the maps obtained
by composing a tangential translation with the Kelvin transform
on each side gives two positive coefficients $z_+$ and $z_-$
in the equation for spherical harmonics of degree $1$.
Their geometric mean satisfies~\eqref{eq:s5-positive-supersolution}.
Picone's identity~\eqref{eq:s5-picone} gives the lower bound
\eqref{eq:s5-nonaxial-gap} for degree $1$. Since
$\lambda_\ell\geq\lambda_1$ for $\ell\geq1$, the same bound holds
in every positive degree. Applying it to cutoffs of a
polynomially growing Jacobi field shows that its coefficients
of positive degree vanish.
The cylindrically symmetric fields are then determined by
the parity uniqueness in Proposition~3.1(iii) of~\cite{PartI}.
Differentiating $H(v_h)=h$ gives the derivative in the
$\Psi_h$ direction. The proof is given in
Section~\ref{subsec:s5-polynomial}.

\begin{introtheorem}\label{intro:cylindrical-kernel}
Let $\xi_*>0$ and let $A_{\xi_*}$ be the positive
cylindrically symmetric function determined in
\eqref{eq:s5-stationary-ground-state}, with $A_{\xi_*}(L)=1$.
Then
\[
 \begin{aligned}
 \mathcal K_{\mathrm{pol}}(C_n)
 &=\mathcal K_{\mathrm{bd}}(C_n)\\
 &=\operatorname{span}\bigl\{
   A_{\xi_*}(\beta)\cos(\xi_*t),\,
   A_{\xi_*}(\beta)\sin(\xi_*t)\bigr\}.
 \end{aligned}
\]
\end{introtheorem}

\noindent\emph{Idea of the proof.}
The same estimate in positive spherical harmonic degrees
shows that every $\Phi\in\mathcal K_{\mathrm{pol}}(C_n)$ is
cylindrically symmetric. The equation for $A_{\xi_*}$ gives
the cosine and sine fields. Subtracting multiples of them
from the even and odd parts of $\Phi$ makes the value of the
even part and the $t$-derivative of the odd part vanish at
$(0,L)$, respectively. Proposition~3.1(iii) of~\cite{PartI}
then implies that both remainders vanish. The proof is given in
Section~\ref{subsec:s5-stationary-kernel}. At $C_n$, the
derivative of $H$ vanishes, and
Proposition~\ref{s5:stationary-energy} gives a strictly
positive quadratic term on each nonzero field in the theorem.

These classifications identify the directions in which a
singular end can be deformed to first order. In the Yamabe
problem without boundary, Mazzeo, Pollack and
Uhlenbeck~\cite[Sections~IV and~V]{MPU} add the corresponding
Jacobi fields on each end to the domain of the linearized
operator in weighted spaces. If the global linearized
operator has trivial $L^2$ kernel, the enlarged domain
gives surjectivity and permits the use of the implicit
function theorem to describe the moduli space locally. Here
the two directions along a periodic solution describe changes
in its translation and Hamiltonian, and only the translation
remains when the Hamiltonian is fixed.

We also consider the self-adjoint operator $\mathcal A_v$
on $L^2(\R\times\Sph^{n-1}_+)$, defined by the interior
operator in~\eqref{eq:jacobi-linearized-equation} with its
boundary condition, and the self-adjoint operator
$\mathcal J_v$ on $L^2(\R\times\partial\Sph^{n-1}_+)$
obtained by subtracting multiplication by
$(nc/(n-2))\tilde v^{2/(n-2)}$ from the Dirichlet-to-Neumann
operator for the same interior equation. For each periodic
Delaunay solution and for $C_n$, both operators have trivial
kernel. In both realizations, a kernel element gives a bounded
Jacobi field. The classifications make this field periodic,
so square integrability forces it to vanish. Cutoffs of a
nonzero periodic Jacobi field give a sequence of unit vectors
that converges weakly to zero and whose images under
$\mathcal A_v$ converge strongly to zero. Thus $0$ belongs
to the essential spectrum of $\mathcal A_v$. See
Section~\ref{subsec:s5-l2} for these statements and the
operator domains. For the fractional Yamabe equation on
$\R\times\Sph^{n-1}$, Andrade, DelaTorre, do~\'O, Ratzkin and
Wei~\cite[Theorem~3]{Andrade} proved triviality of the $L^2$
kernel of the linearized operator for orders $s<1$
sufficiently close to $1$.

The flat Hamiltonian parametrization and the Jacobi
classification in a growth class containing fields of
polynomial growth were also obtained by Jin, Liao and
Xiong~\cite[Theorems~1.1 and~1.3 and Proposition~A.2]{JinLiaoXiong}.
Our proofs use the flat classification and the parity
uniqueness of~\cite{PartI}, together with the estimates
for the equations in positive spherical harmonic degrees.

For $\eps\in(0,1]$ and $\tau\in\R$, let $U_{\eps,\tau}$ be
the flat solution in~\eqref{eq:main-comparison}. In the
normalization of~\eqref{eq:intro-geometric}, it is given by
\[
 U_{\eps,\tau}(r\theta)
 = b_n r^{-\frac{n-2}{2}}
        v_\eps(-\log r+\tau,\theta).
\]

\begin{introtheorem}\label{intro:asymptotics}
Let $n\geq3$, and let $U$ be a positive smooth solution of
\eqref{eq:intro-geometric} on $B_\rho^+\setminus\{0\}$ in
boundary Fermi coordinates centered at $p$. Suppose there
are $c_0,C_0>0$ such that, for $0<|\overline x|<\rho$,
\[
 c_0|\overline x|^{-\frac{n-2}{2}}
 \leq u(\overline x)
 \leq C_0|\overline x|^{-\frac{n-2}{2}} .
\]
Then there exist $\delta>0$, $\rho_0\in(0,\rho)$,
$\eps_*\in(0,1]$, $\tau_*\in\R$, and $C>0$ such that,
for $x\in B_{\rho_0}^+\setminus\{0\}$,
\[
 \sum_{j=0}^2 |x|^{\frac{n-2}{2}+j}
       |D^j(U-U_{\eps_*,\tau_*})(x)|
 \leq C|x|^\delta .
\]
Here $D^j$ denotes derivatives in the fixed Fermi chart.
In particular, $U(x)=(1+O(|x|^\delta))U_{\eps_*,\tau_*}(x)$
as $x\to0$. The necksize $\eps_*$ is unique. If $\eps_*<1$,
$\tau_*$ is unique modulo $T_{\eps_*}$. If $\eps_*=1$, the
flat solution is $b_n r^{-(n-2)/2}C_n(\theta)$ for every
$\tau_*$, and its normalized transform is independent of $t$.
\end{introtheorem}

\noindent\emph{Idea of the proof.}
The bounds on the trace give uniform positive upper and
lower bounds and derivative estimates for the rescalings of
$U$ on fixed half-annuli. The equation for its transform $w$
in~\eqref{eq:main-change} differs from
\eqref{eq:intro-normalized-model} by terms of order $e^{-t}$.
The same bound for the derivative of the Hamiltonian gives
convergence of the Pohozaev integral as the radius tends to
zero. For any sequence $t_j\to\infty$, the translated functions
$w(t_j+\cdot,\cdot)$ have a subsequence converging on compact
subsets to a positive global solution of
\eqref{eq:intro-normalized-model}. Theorem~B of~\cite{PartI}
identifies that solution as $C_n$ or a translate of a
Delaunay solution. The Pohozaev invariant gives the same
Hamiltonian for all such subsequential limits, and
Theorem~\ref{intro:energy} determines the same necksize.

If $\eps_*<1$, an orthogonality condition on a fixed interval
of the cylinder selects the translation there. If the error
failed to decrease uniformly on successive intervals, its
normalization would give a nonzero Jacobi field of polynomial
growth with zero Hamiltonian derivative and orthogonal to
$\partial_t v_{\eps_*}$. Theorem~\ref{intro:polynomial-kernel}
excludes such a field. The resulting exponential estimate
gives $|a'(t)-1|=O(e^{-\delta t})$ for some $\delta>0$,
where $a(t)$ is the translation selected on the interval
centered at $t$. Consequently, $a(t)-t\to\tau_*$. When
$\eps_*=1$, the positive quadratic term of the Hamiltonian
on the space in Theorem~\ref{intro:cylindrical-kernel} gives
the corresponding decay estimate for $w-C_n$. Returning to
the Fermi chart gives the estimates through order $2$.
The proof is given in Section~\ref{subsec:asymptotic-convergence}.

\subsection*{Structure of the paper}
\label{subsec:intro-organization}

Section~\ref{sec:preliminaries} fixes the geometric and
cylindrical conventions, establishes the local estimates,
and proves the existence of the Pohozaev invariant.
Section~\ref{sec:energy-coordinate} classifies the
cylindrically symmetric Jacobi fields with polynomial growth
and proves the Hamiltonian parametrization.
Section~\ref{sec:full-jacobi} classifies the Jacobi fields
without a symmetry assumption and studies the two
self-adjoint operators on the cylinder and its boundary.
Section~\ref{sec:fixed-profile} proves the asymptotic theorem
and studies the effect of boundary-preserving conformal
deformations of the flat solutions on the error in Fermi
coordinates. Appendix~\ref{app:cylindrical-computations}
gives the computations used in these proofs.

\subsection*{Relation to concurrent work}

While this manuscript was being completed, Jin, Liao and Xiong
submitted their preprint~\cite{JinLiaoXiong} to arXiv on
November~5, 2026. Their Theorem~1.1 and Proposition~A.2
give the Hamiltonian parametrization of the flat Delaunay family
and its real-analytic relation to the necksize. Their Theorem~1.3
gives the Jacobi field classifications in
Theorems~\ref{intro:polynomial-kernel}
and~\ref{intro:cylindrical-kernel} for a class that also allows
sufficiently slow exponential growth at both ends. Their
Theorem~1.4 proves convergence to a unique global singular
solution in a Euclidean half-ball, with a relative error that
decays at a positive power rate.

In the flat problem, we use the classification and parity
uniqueness of~\cite{PartI}. For the positive spherical harmonic
degrees, we use the geometric mean of the two positive
degree-$1$ coefficients in Picone's identity to obtain the
coercive estimate in Proposition~\ref{s5:nonaxial-gap}.
Section~\ref{subsec:s5-l2} proves that both $L^2$ kernels
are trivial and that zero belongs to the essential spectrum
of the linearized operator on the cylinder.

In the periodic case, both asymptotic arguments use
orthogonality to the translation field and the first variation
of the Hamiltonian. Under the critical two-sided bounds on
the trace, Theorem~\ref{intro:asymptotics} extends convergence
to a fixed flat solution to smooth background metrics, with
weighted estimates through two derivatives in boundary Fermi
coordinates. The limit of the Pohozaev integral determines
the necksize. Proposition~\ref{prop:first-order-geometric-obstruction}
shows that the boundary-preserving conformal deformations
considered here can leave a first-order error even for a
conformally flat background metric.

\subsection*{Statement on the use of AI}

The need for Proposition~3.1 of~\cite{PartI} was identified
before we used AI. In the mathematical work, we used the
model GPT-Sol 6 only to assist with the calculation in
that proposition. The proposition is also used in the
calculations in the present article. We also used the model
to revise the manuscript. All mathematical statements,
proofs, and citations were independently checked by the
author, who takes full responsibility for the contents of
the article. The final manuscript was written and revised
by the author.

\subsection*{Acknowledgments}

This work forms part of my Ph.D.\ thesis at Universidade Federal
Fluminense. I am grateful to my advisor, S\'ergio Almaraz, for his
guidance and advice throughout this work. I also thank
Xuezhang Chen for his warm welcome at Nanjing University and his
generous support during my stay.

This study was financed in part by the Coordena\c{c}\~ao de
Aperfei\c{c}oamento de Pessoal de N\'{\i}vel Superior - Brasil
(CAPES) - Finance Code 001. I received support under process
No.~88887.815246/2023-00 and a scholarship through the Programa
Institucional de Doutorado Sandu\'{\i}che no Exterior (PDSE),
process No.~88881.127484/2025-01.

\section{Preliminaries}\label{sec:preliminaries}

We first fix the boundary Fermi coordinates at $p$ for the Escobar--Yamabe
equation~\eqref{eq:main-geometric}. We then pass from the flat
Escobar--Yamabe equation~\eqref{eq:flat-boundary} to its normalized
cylindrical form~\eqref{eq:s3-normalized-model}, recall the results
on Delaunay solutions from the companion paper, and give estimates
for the linearized equation. In
Section~\ref{subsec:geometric-estimates}, we establish local
estimates for~\eqref{eq:main-geometric} and prove the existence
of the limiting Pohozaev invariant in~\eqref{eq:proof-pohozaev}.

\subsection{Notation}\label{subsec:notation-fermi}

We denote by $\bar g$ the metric induced by $g$ on $\partial M$
and by $\nu_g$ the outward unit normal to $\partial M$. Let
$\nabla^g$, $\Delta_g$, and $R_g$ denote the Levi-Civita connection,
the Laplace--Beltrami operator, and the scalar curvature of $g$,
respectively. For vectors $X,Y$ tangent to $\partial M$, we set
\[
 \mathrm{II}_g(X,Y)=g(\nabla^g_X\nu_g,Y),
 \qquad
 h_g=\frac{1}{n-1}\operatorname{tr}_{\bar g}\mathrm{II}_g .
\]

We use the half-space notation of Almaraz~\cite[Section~2.1]{Almaraz}.
Write $x=(\overline x,x_n)\in\R^{n-1}\times\R$ and let
$\R^n_+=\{x\in\R^n:x_n\geq0\}$. For $r>0$, let $B_r(0)$ be
the open Euclidean ball of radius $r$ centered at $0$ and set
$B_r^+=B_r(0)\cap\R^n_+$. Its flat and spherical faces are denoted
by $\partial_0B_r^+=B_r^+\cap\partial\R^n_+$ and
$\partial^+B_r^+=\partial B_r(0)\cap\R^n_+$, respectively. We
write $\Sph^{n-1}_+=\Sph^{n-1}\cap\R^n_+$, so
$\partial^+B_r^+=r\Sph^{n-1}_+$. The standard Euclidean metric
is denoted by $\delta$, and its Laplacian by
$\Delta=\Delta_\delta=\sum_{i=1}^n\partial_{x_i}^2$.

Choose a $\bar g$-orthonormal basis $E_1,\ldots,E_{n-1}$ of
$T_p\partial M$. For $\overline x$ near $0$, set
$q(\overline x)=\exp_p^{\bar g}(\sum_{a=1}^{n-1}x_aE_a)$. The
boundary Fermi map centered at $p$ is
\begin{equation}\label{eq:prelim-fermi-map}
 \Psi(\overline x,x_n)=
 \exp_{q(\overline x)}^g\bigl(-x_n\nu_g(q(\overline x))\bigr),
 \qquad x_n\geq0 .
\end{equation}
See Almaraz~\cite[Section~2.3]{Almaraz} and
Marques~\cite[Section~3]{MarquesFermi}. For sufficiently small
$\rho>0$, the map $\Psi$ identifies $B_\rho^+$ with a
neighborhood of $p$ in $M$ and sends $0$ to $p$. We also
write $g$ for the pullback $\Psi^*g$ on $B_\rho^+$, with
coefficients
\[
 g_{ij}(x)=g_{\Psi(x)}(d\Psi_x e_i,d\Psi_x e_j).
\]
The Gauss lemma gives
\begin{equation}\label{eq:proof-fermi}
 g_{nn}=1,\qquad g_{an}=0
\end{equation}
for $1 \leq a \leq n-1$. The coordinate frame at $0$ is orthonormal, so
$g_{ij}(0)=\delta_{ij}$. Since $x_n$ increases inward, the outward unit normal on the
flat face is $-\partial_{x_n}$. In Theorem~\ref{intro:asymptotics},
$|x|$ and $|\overline x|$ are Euclidean norms, and $D^j$ denotes
derivatives of order $j$ with respect to the Fermi coordinates $x$.

The tangential derivatives of $g_{ab}$ vanish at $0$ because
$\overline x$ consists of geodesic normal coordinates for $\bar g$.
With the outward normal convention above,
$\partial_{x_n}g_{ab}(0)=-2\mathrm{II}_{ab}(p)$, where
$\mathrm{II}_{ab}(p)=\mathrm{II}_g(E_a,E_b)$. Thus
\begin{equation}\label{eq:prelim-fermi-expansion}
 g_{ab}(x)=\delta_{ab}-2\mathrm{II}_{ab}(p)x_n+O(|x|^2),
\end{equation}
for $1\leq a,b\leq n-1$. After writing $x=rz$ on a fixed
half-annulus, the remainder is $O(r^2)$ in every fixed $C^k$
norm in $z$.

The boundary is \emph{umbilic at $p$} if
$\mathrm{II}_g(p)=h_g(p)\bar g_p$. In that case,
\[
 g_{ab}(x)=(1-2h_g(p)x_n)\delta_{ab}+O(|x|^2).
\]
\begin{definition}\label{def:conformally-flat}
A smooth metric $g$ is \emph{conformally flat} if every point
has a neighborhood on which $g=\varphi^{4/(n-2)}g_0$ for some
positive smooth function $\varphi$ and some flat metric $g_0$.
\end{definition}

The conformal Laplacian and the boundary operator are
\[
 L_g=-\Delta_g+\frac{n-2}{4(n-1)}R_g,
 \qquad B_g=\partial_{\nu_g}+\frac{n-2}{2}h_g.
\]
Let $U$ be a positive smooth function on
$B_\rho^+\setminus\{0\}$, and write $u(\overline x)=U(\overline x,0)$
for its trace on $\partial_0B_\rho^+\setminus\{0\}$. For $\kappa>0$,
the singular Escobar--Yamabe equation in this chart is
\begin{equation}\label{eq:main-geometric}
 \begin{cases}
 L_gU=0 & \text{in } \operatorname{int}(B_\rho^+),\\
 B_gU=\kappa u^{n/(n-2)}
       & \text{on } \partial_0B_\rho^+\setminus\{0\}.
 \end{cases}
\end{equation}
By conformal covariance~\cite{Escobar}, the metric $U^{4/(n-2)}g$
is scalar flat and has constant boundary mean curvature
$2\kappa/(n-2)$ on $\partial_0B_\rho^+\setminus\{0\}$.
Multiplying $U$ by $(2\kappa/(n-2))^{(n-2)/2}$ changes the
boundary coefficient to $(n-2)/2$ and gives
\eqref{eq:intro-geometric}, with boundary mean curvature one.
We call the singularity
of the trace removable if $u$ extends continuously to
$\overline x=0$.

For a nonnegative integer $k$ and $0<\gamma<1$, we use
$C^k(M)$ for the space of functions with continuous derivatives
through order $k$ and $C^{k,\gamma}(M)$ when the derivatives of
order $k$ are also $\gamma$-H\"older continuous. We write
$C^\infty(M)$ for smooth functions and $C_c^\infty(M)$ for smooth
functions of compact support. The subscript $\mathrm{loc}$
denotes regularity on compact subsets. On a closed
half-annulus or cylinder, the indicated norms include the
boundary. We use $L^q(M)$ for the Lebesgue spaces and $W^{k,q}(M)$
for the Sobolev spaces with weak derivatives through order
$k$ in $L^q(M)$. The notation $H^s(M)$ denotes the $L^2(M)$-based
Sobolev space of order $s$. In particular, $H^{-1/2}(M)$ is the
dual of $H^{1/2}(M)$ under the $L^2(M)$ pairing. The spaces are real
unless we take their complexifications for Fourier analysis.
The periodic, uniformly local, and weighted spaces on the
cylinder are specified in Section~\ref{subsec:linearized-spaces}.

\subsection{The flat equation and the cylindrical normalization}
\label{subsec:flat-model}

We write the flat Escobar--Yamabe equation in Emden--Fowler
coordinates and describe the action of dilations and Kelvin
inversion. We then normalize the cylindrical solution to have
trace one and use the same normalization to define the Hamiltonian.

We use the flat Escobar--Yamabe equation in the normalization
of the companion paper~\cite{PartI},
\begin{equation}\label{eq:flat-boundary}
 \begin{cases}
 \Delta U=0 & \text{for }x_n>0,\\
 \partial_{x_n}U+\dfrac{n-2}{2}u^{n/(n-2)}=0
     & \text{on }\partial\R^n_+\setminus\{0\},
 \end{cases}
\end{equation}
where $u(\overline x)=U(\overline x,0)$. A classical global
solution of~\eqref{eq:flat-boundary} is a positive function
$U\in C^\infty(\R^n_+\setminus\{0\})$ satisfying both equations.
This smoothness includes $\partial\R^n_+\setminus\{0\}$.
We impose no growth bound on $U$ as $|x|\to\infty$ and assume
cylindrical symmetry only where stated.
The metric $U^{4/(n-2)}\delta$ is scalar flat and has mean
curvature one on $\partial\R^n_+\setminus\{0\}$.

Write $x=r\theta$, with $r=|x|>0$ and $\theta\in\Sph^{n-1}_+$,
and set $t=-\log r$. We equip $\R\times\Sph^{n-1}_+$ with
the product metric and volume element $\dd\mu=\dd t\,\dd\theta$,
where $\dd\theta$ is the round measure on $\Sph^{n-1}_+$.
We write $\dd\sigma$ for the round measure on
$\partial\Sph^{n-1}_+=\Sph^{n-2}$. For a smooth function $f$
on $\R\times\Sph^{n-1}_+$, we denote its trace on
$\R\times\partial\Sph^{n-1}_+$ by $\tilde f$. Thus
$\tilde f(t,\vartheta)=f(t,(\vartheta,0))$ for
$\vartheta\in\Sph^{n-2}$.

Let $\beta=\arccos\theta_n$ and $L=\pi/2$. We write
\[
 \theta=(\sin\beta\,\vartheta,\cos\beta),
 \qquad 0\leq\beta\leq L,\qquad \vartheta\in\Sph^{n-2}.
\]
The outward unit conormal to $\partial\Sph^{n-1}_+$ is
$\nu_{\Sph^{n-1}_+}$, with
$\partial_{\nu_{\Sph^{n-1}_+}}=\partial_\beta$ at $\beta=L$.
We define the inward unit conormal $\eta_{\Sph^{n-1}_+}$
as $-\nu_{\Sph^{n-1}_+}$.
For $0<\beta<L$, the round Laplacian is
\[
 \Delta_{\Sph^{n-1}_+}=\partial_{\beta\beta}
      +(n-2)\cot\beta\,\partial_\beta
      +\sin^{-2}\beta\,\Delta_{\Sph^{n-2}}.
\]
The round measure is
$\dd\theta=\omega(\beta)\,\dd\beta\,\dd\sigma(\vartheta)$,
where $\omega(\beta)=\sin^{n-2}\beta$.

For a solution $U$ of~\eqref{eq:flat-boundary}, set
$V(t,\theta)=r^{(n-2)/2}U(r\theta)$ on $\R\times\Sph^{n-1}_+$.
The Euclidean polar formula gives the interior equation for $V$.
On $\R\times\partial\Sph^{n-1}_+$, the outward conormal
derivative of $V$ is $-r^{n/2}\partial_{x_n}U(r\vartheta,0)$,
and $\tilde V(t,\vartheta)=r^{(n-2)/2}u(r\vartheta)$. Hence
\eqref{eq:flat-boundary} becomes
\begin{equation}\label{eq:flat-cylinder-unnormalized}
 \begin{cases}
 V_{tt}+\Delta_{\Sph^{n-1}_+}V-\dfrac{(n-2)^2}{4}V=0
       & \text{in }\R\times\operatorname{int}\Sph^{n-1}_+,\\
 \partial_{\nu_{\Sph^{n-1}_+}}V=\dfrac{n-2}{2}\tilde V^{n/(n-2)}
       & \text{on }\R\times\partial\Sph^{n-1}_+.
 \end{cases}
\end{equation}
The dilation $U_\lambda(x)=\lambda^{(n-2)/2}U(\lambda x)$,
for $\lambda>0$, corresponds to
\begin{equation}\label{eq:flat-dilation-translation}
 V_\lambda(t,\theta)=V(t-\log\lambda,\theta).
\end{equation}
For $y=(\overline y,0)$ and $\rho>0$, the Kelvin transform is
\begin{equation}\label{eq:flat-kelvin}
 U_{y,\rho}(x)=\left(\frac{\rho}{|x-y|}\right)^{n-2}
 U\left(y+\frac{\rho^2(x-y)}{|x-y|^2}\right).
\end{equation}
It preserves~\eqref{eq:flat-boundary} on the inverted domain.
When $y=0$, its action is $t\mapsto-t-2\log\rho$.

Following Caffarelli, Jin, Sire and Xiong~\cite[Theorem~1.3]{CJXS},
we call $U$ \emph{cylindrically symmetric} if
$U(\overline x,x_n)=U(|\overline x|,x_n)$.
In this case the Emden--Fowler transform is independent of
$\vartheta$ and is given by
\[
 V(t,\beta)=e^{-\frac{n-2}{2}t}
             U(e^{-t}\sin\beta,e^{-t}\cos\beta).
\]
Its trace is $\tilde V(t)=V(t,L)$. The angular part of the
equation reduces to $\omega^{-1}(\omega V_\beta)_\beta$
for $0<\beta<L$. Smoothness in spherical charts and evenness
along meridians give $V_\beta(t,0)=0$.

To normalize~\eqref{eq:flat-cylinder-unnormalized}, let $C_n$
be the solution of
\[
 \begin{cases}
 \left(-\Delta_{\Sph^{n-1}_+}+\dfrac{(n-2)^2}{4}\right)C_n=0
     & \text{in }\operatorname{int}\Sph^{n-1}_+,\\
 C_n=1 & \text{on }\partial\Sph^{n-1}_+.
 \end{cases}
\]
The Dirichlet problem has a unique smooth solution. The maximum
principle gives $C_n>0$, and uniqueness gives cylindrical
symmetry. Thus $C_n=C_n(\beta)$, with $C_n'(0)=0$ and $C_n(L)=1$.
Following~\cite[Section~3.1]{PartI}, we set
\begin{equation}\label{eq:prelim-normalizing-constants}
 c=C_n'(L)=2\left(\frac{\Gamma(n/4)}{\Gamma((n-2)/4)}\right)^2,
 \qquad b_n=\left(\frac{2c}{n-2}\right)^{\frac{n-2}{2}}.
\end{equation}
Then $v=b_n^{-1}V$ satisfies
\begin{equation}\label{eq:s3-normalized-model}
 \begin{cases}
 v_{tt}+\Delta_{\Sph^{n-1}_+}v-\dfrac{(n-2)^2}{4}v=0
       & \text{in }\R\times\operatorname{int}\Sph^{n-1}_+,\\
 \partial_{\nu_{\Sph^{n-1}_+}}v=c\tilde v^{n/(n-2)}
       & \text{on }\R\times\partial\Sph^{n-1}_+.
 \end{cases}
\end{equation}
A solution is \emph{cylindrical} if it is independent of $t$ and
\emph{noncylindrical} otherwise. By~\cite[Theorem~B]{PartI},
$C_n$ is the unique positive smooth cylindrical solution
of~\eqref{eq:s3-normalized-model}. The corresponding solution of
\eqref{eq:flat-boundary} is $b_nr^{-(n-2)/2}C_n(\theta)$.

For a cylindrically symmetric $v$, we write $\tilde v(t)=v(t,L)$.
A two-sided bound on the trace means that
$m\leq\tilde v(t)\leq M$ for all $t\in\R$, for some $M,m>0$.
If $v$ solves~\eqref{eq:s3-normalized-model}, then
$\widehat v=c^{(n-2)/2}v$ satisfies the same interior equation
with boundary coefficient one.

Let $I\subset\R$ be an open interval and let
$v\in C^\infty(I\times\Sph^{n-1}_+)$ satisfy $\tilde v>0$ on
$I\times\partial\Sph^{n-1}_+$. For $a>0$ and $t\in I$, define
the Hamiltonian per unit area of $\Sph^{n-2}$ by
\begin{equation}\label{eq:proof-hamiltonian}
 \begin{aligned}
 |\Sph^{n-2}|H_a(v)(t)
 &=\frac12\int_{\Sph^{n-1}_+}
       \left(v_t^2-|\nabla_{\Sph^{n-1}_+}v|^2
                         -\frac{(n-2)^2}{4}v^2\right)\dd\theta\\
 &\quad+\frac{(n-2)a}{2(n-1)}\int_{\partial\Sph^{n-1}_+}
                  \tilde v^{\frac{2(n-1)}{n-2}}\dd\sigma .
 \end{aligned}
\end{equation}
The integrands are evaluated at $t$. Suppose $v$ solves
\eqref{eq:s3-normalized-model} with boundary coefficient $a$
in place of $c$. Differentiating in $t$ and integrating by
parts on $\Sph^{n-1}_+$ give
\[
 \begin{aligned}
 |\Sph^{n-2}|\frac{\dd}{\dd t}H_a(v)(t)
 &=\int_{\Sph^{n-1}_+}v_t
      \left(v_{tt}+\Delta_{\Sph^{n-1}_+}v-\frac{(n-2)^2}{4}v\right)
      \dd\theta\\
 &\quad+\int_{\partial\Sph^{n-1}_+}\partial_t\tilde v
      \left(a\tilde v^{n/(n-2)}-\partial_{\nu_{\Sph^{n-1}_+}}v\right)
      \dd\sigma .
 \end{aligned}
\]
Both integrals vanish by the equation for $v$. Thus $H_a(v)(t)$
is independent of $t$, and we write $H_a(v)$ for its value.
When $a=c$, we abbreviate the Hamiltonian by $H(v)(t)$ and
its constant value along a solution by $H(v)$. For a
cylindrically symmetric solution of~\eqref{eq:s3-normalized-model},
integration in $\vartheta$ gives
\begin{equation}\label{eq:s3-hamiltonian}
 \begin{aligned}
 H(v)
 &=\frac12\int_0^L\omega(\beta)
       \left(v_t^2-v_\beta^2-\frac{(n-2)^2}{4}v^2\right)\dd\beta\\
 &\quad+\frac{(n-2)c}{2(n-1)}
                  \tilde v(t)^{\frac{2(n-1)}{n-2}} .
 \end{aligned}
\end{equation}
Multiplying the angular equation for $C_n$ by $C_n$ and integrating
by parts gives
\[
 \int_0^L\omega(\beta)
       \left((C_n')^2+\frac{(n-2)^2}{4}C_n^2\right)\dd\beta
 =[\omega C_n C_n']_0^L=c .
\]
It follows that $H(C_n)=-c/(2(n-1))$. We denote this value by $H_c$.

The identity $((n-2)/2)b_n^{2/(n-2)}=c$ gives
\[
 H_{(n-2)/2}(b_nv)(t)=b_n^2H(v).
\]
Section~\ref{subsec:s4-physical-energy} uses this relation to
classify the dilation orbits of singular solutions of
\eqref{eq:flat-boundary}. The relation between the Hamiltonian
and the Pohozaev integral for the geometric equation is given in
Section~\ref{subsec:geometric-estimates}.

\subsection{Linearized spaces and local estimates}
\label{subsec:linearized-spaces}
\label{subsec:estimates-compactness}

We give local estimates for the normalized equation and its
linearization. We then compare the growth conditions used to
classify Jacobi fields and prove the periodic Robin solvability
result used in Section~\ref{subsec:s4-axial-jacobi}.

Let $k\geq0$ be an integer and let $0<\gamma<1$.
The space $C^{k,\gamma}_{\mathrm{loc}}(\R\times\Sph^{n-1}_+)$
consists of functions with $C^{k,\gamma}$ regularity on
$I\times\Sph^{n-1}_+$ for every compact interval $I\subset\R$.
Near the equator, write $s=L-\beta$. Then
\[
 \theta=(\cos s\,\vartheta,\sin s),
 \qquad \vartheta\in\Sph^{n-2},\quad s\geq0.
\]
Thus $s=0$ on the equator, and $s$ measures distance from the
equator along the inward meridian. A local chart on the cylinder
uses $t$, local coordinates for $\vartheta$ on $\Sph^{n-2}$,
and $s$. In these coordinates, regularity up to the equator means that
all derivatives through order $k$ extend continuously to $s=0$
and that the derivatives of order $k$ are H\"older continuous
with exponent $\gamma$ up to and including $s=0$.
We write $C^\infty_{\mathrm{loc}}(\R\times\Sph^{n-1}_+)$ for
the space of functions that are smooth in these charts on every
product $I\times\Sph^{n-1}_+$ with $I\subset\R$ compact.

The space $C^{k,\gamma}_{\mathrm{uloc}}(\R\times\Sph^{n-1}_+)$
consists of functions in
$C^{k,\gamma}_{\mathrm{loc}}(\R\times\Sph^{n-1}_+)$ with finite norm
\[
 \|z\|_{C^{k,\gamma}_{\mathrm{uloc}}(\R\times\Sph^{n-1}_+)}
 =\sup_{\tau\in\R}
     \|z\|_{C^{k,\gamma}([\tau-1,\tau+1]\times\Sph^{n-1}_+)} .
\]
The boundary space
$C^{k,\gamma}_{\mathrm{uloc}}(\R\times\partial\Sph^{n-1}_+)$
is defined using $[\tau-1,\tau+1]\times\partial\Sph^{n-1}_+$.

Let $v$ be a smooth solution of~\eqref{eq:s3-normalized-model}.
If $0<c_0\leq v\leq C_0$ on $\R\times\Sph^{n-1}_+$, then
\[
 \|v\|_{C^{k,\gamma}([\tau-1,\tau+1]\times\Sph^{n-1}_+)}
 \leq C(n,c_0,C_0,k,\gamma)
 \qquad\text{for all }\tau\in\R .
\]
The estimates follow from the local Neumann and Schauder
estimates on nested intervals inside $[\tau-2,\tau+2]$,
as in~\cite[Proposition~2.1]{PartI}. Translation in $t$ makes
the constants independent of $\tau$. A diagonal argument
gives a subsequence converging in
$C^\infty_{\mathrm{loc}}(\R\times\Sph^{n-1}_+)$ for any sequence
of solutions with the same bounds.

Under the bound $0\leq v\leq C_0$, bounded Neumann data give
a uniform local H\"older estimate for $v$. Since
$s\mapsto cs^{n/(n-2)}$ is Lipschitz on $[0,C_0]$, the Neumann
data are also uniformly H\"older continuous. The Neumann
Schauder estimates then give a uniform bound in
$C^{1,\delta_0}([\tau-1,\tau+1]\times\Sph^{n-1}_+)$ for some
$0<\delta_0<1$. We use this bound in
Section~\ref{subsec:s4-global-energy} when the necksize tends
to zero.

For a positive smooth solution $v$ of~\eqref{eq:s3-normalized-model},
a smooth function $\Phi$ on $\R\times\Sph^{n-1}_+$ is a
\emph{Jacobi field} along $v$ if
\begin{equation}\label{eq:jacobi-linearized-equation}
 \begin{aligned}
 \mathcal A\Phi
 &=-\Phi_{tt}-\Delta_{\Sph^{n-1}_+}\Phi
       +\frac{(n-2)^2}{4}\Phi=0
       &&\text{in }\R\times\operatorname{int}\Sph^{n-1}_+,\\
 \mathcal B_v\Phi
 &=\partial_{\nu_{\Sph^{n-1}_+}}\Phi
       -\frac{nc}{n-2}\tilde v^{2/(n-2)}\tilde\Phi=0
       &&\text{on }\R\times\partial\Sph^{n-1}_+.
 \end{aligned}
\end{equation}
Smoothness at the pole is understood in spherical charts.

For $q\geq0$, set $\langle t\rangle=(1+t^2)^{1/2}$. The space
$C^{k,\gamma}_{\mathrm{pol},q}(\R\times\Sph^{n-1}_+)$ consists
of functions in $C^{k,\gamma}_{\mathrm{loc}}(\R\times\Sph^{n-1}_+)$
with finite norm
\begin{equation}\label{eq:jacobi-polynomial-space}
 \|z\|_{C^{k,\gamma}_{\mathrm{pol},q}(\R\times\Sph^{n-1}_+)}
 =\sup_{\tau\in\R}\langle\tau\rangle^{-q}
       \|z\|_{C^{k,\gamma}([\tau-1,\tau+1]\times\Sph^{n-1}_+)} .
\end{equation}
We define $C^{k,\gamma}_{\mathrm{pol},q}(\R\times\partial\Sph^{n-1}_+)$
using $[\tau-1,\tau+1]\times\partial\Sph^{n-1}_+$ in the norm.
For $q=0$ the corresponding spaces are the uniformly local ones.

A Jacobi field has polynomial growth of order at most $q$ if
there is a constant $C>0$ such that
\begin{equation}\label{eq:jacobi-pointwise-growth}
 |\Phi(t,\theta)|\leq C\langle t\rangle^q
 \qquad\text{for all }(t,\theta)\in\R\times\Sph^{n-1}_+ .
\end{equation}
The following lemma allows us to use the same growth class in
Section~\ref{subsec:s4-axial-jacobi} and
Section~\ref{sec:full-jacobi}.

\begin{lemma}\label{lem:jacobi-growth-equivalence}
Let $v$ be a positive smooth periodic solution
of~\eqref{eq:s3-normalized-model}, and let $\Phi$ be a Jacobi
field along $v$. For each $q\geq0$, the following are equivalent.
\begin{enumerate}[label=\textup{(\roman*)}]
 \item The pointwise bound~\eqref{eq:jacobi-pointwise-growth} holds.
 \item There is a constant $C>0$ such that
 \[
  \|\Phi\|_{L^2([\tau-1,\tau+1]\times\Sph^{n-1}_+)}
       \leq C\langle\tau\rangle^q\qquad\text{for all }\tau\in\R .
 \]
 \item $\Phi\in C^{k,\gamma}_{\mathrm{pol},q}(\R\times\Sph^{n-1}_+)$
 for every integer $k\geq0$ and $0<\gamma<1$.
\end{enumerate}
\end{lemma}

\begin{proof}
On $[\tau-1,\tau+1]$, we have
$\langle t\rangle\leq 2\langle\tau\rangle$. The measure of
$[\tau-1,\tau+1]\times\Sph^{n-1}_+$ is independent of $\tau$,
so \textup{(i)} implies \textup{(ii)}.

For every $\delta>0$ there is $C_\delta>0$ such that, for any
compact interval $I\subset\R$ and smooth $z$ on
$I\times\Sph^{n-1}_+$,
\begin{equation}\label{eq:jacobi-trace-inequality}
 \begin{aligned}
 \int_I\int_{\partial\Sph^{n-1}_+}|\tilde z|^2\,\dd\sigma\,\dd t
 &\leq \delta\int_I\int_{\Sph^{n-1}_+}
                   |\nabla_{\Sph^{n-1}_+}z|^2\,\dd\theta\,\dd t\\
 &\quad+C_\delta\int_I\int_{\Sph^{n-1}_+}
                   |z|^2\,\dd\theta\,\dd t .
 \end{aligned}
\end{equation}
To see this, write $z(t,L,\vartheta)$ as $z(t,\beta,\vartheta)$
plus the integral of $\partial_\beta z$ from $\beta$ to $L$.
Apply Cauchy--Schwarz and average in $\beta$ over a collar of the
equator. The coefficient of the derivative term is proportional
to the width of the collar, which can be chosen to give $\delta$.

Fix $\tau\in\R$ and let $I_\tau=[\tau-2,\tau+2]$. Choose a smooth
cutoff $\chi$ supported in the interior of $I_\tau$, equal to
one on $[\tau-3/2,\tau+3/2]$, with derivative bounded
independently of $\tau$. Multiplying the interior equation
in~\eqref{eq:jacobi-linearized-equation} by $\chi^2\Phi$,
integrating by parts and using the boundary equation give
\[
 \begin{aligned}
 &\int_{I_\tau\times\Sph^{n-1}_+}
       \left(|\nabla(\chi\Phi)|^2+
              \frac{(n-2)^2}{4}|\chi\Phi|^2\right)\dd\mu\\
 &\qquad=\int_{I_\tau\times\Sph^{n-1}_+}
               |\chi'|^2|\Phi|^2\,\dd\mu
    +\frac{nc}{n-2}\int_{I_\tau}\int_{\partial\Sph^{n-1}_+}
         \tilde v^{2/(n-2)}|\chi\tilde\Phi|^2\,\dd\sigma\,\dd t ,
 \end{aligned}
\]
where $\nabla$ is the product gradient. Since $\tilde v$ is bounded, we can choose $\delta$ sufficiently
small in~\eqref{eq:jacobi-trace-inequality} to move the angular
derivative term to the left-hand side. The remaining terms are bounded by a constant
times $\|\Phi\|_{L^2(I_\tau\times\Sph^{n-1}_+)}^2$. This bounds
$\|\Phi\|_{H^1([\tau-3/2,\tau+3/2]\times\Sph^{n-1}_+)}$ by
$C\|\Phi\|_{L^2(I_\tau\times\Sph^{n-1}_+)}$.

By periodicity, the coefficient in the boundary condition and
all its derivatives are uniformly bounded. Local elliptic
estimates up to the equator on smaller intervals, followed by
Sobolev embedding, then give
\[
 \|\Phi\|_{C^{k,\gamma}([\tau-1,\tau+1]\times\Sph^{n-1}_+)}
 \leq C_{v,k,\gamma}\|\Phi\|_{L^2(I_\tau\times\Sph^{n-1}_+)} .
\]
The two intervals of length two centered at $\tau-1$ and
$\tau+1$ cover $I_\tau$. Since
$\langle\tau\pm1\rangle\leq2\langle\tau\rangle$, condition
\textup{(ii)} bounds the last norm by $C'\langle\tau\rangle^q$.
This proves \textup{(iii)}. Taking $k=0$ and $\tau=t$ in
\eqref{eq:jacobi-polynomial-space} gives \textup{(i)}.
\end{proof}

For $T>0$, let $C^{k,\gamma}_T(\R\times\Sph^{n-1}_+)$ be the
subspace of $T$-periodic functions in
$C^{k,\gamma}_{\mathrm{uloc}}(\R\times\Sph^{n-1}_+)$, with the
norm on $[0,T]\times\Sph^{n-1}_+$. We define
$C^{k,\gamma}_T(\R\times\partial\Sph^{n-1}_+)$ using
$[0,T]\times\partial\Sph^{n-1}_+$.
The subscripts $\mathrm{sym}$, $\mathrm{even}$ and $\mathrm{odd}$
indicate cylindrical symmetry, evenness in $t$ and oddness in
$t$, respectively.

If $v=v(t,\beta)$ and $z=z(t,\beta)$, then $\mathcal A z$
depends only on $t$ and $\beta$, and $\mathcal B_v z$ depends
only on $t$. If $v(-t,\theta)=v(t,\theta)$ and
$z(-t,\theta)=\pm z(t,\theta)$, then $\mathcal A z$ and
$\mathcal B_v z$ have the same parity in $t$ as $z$. Thus the linearized operators
preserve cylindrical symmetry when $v$ is cylindrically
symmetric, and preserve both parities when $v$ is even in $t$.

For the periodic weak formulation, write
$H^1_T(\R\times\Sph^{n-1}_+)$ for the space of $T$-periodic
functions in $H^1_{\mathrm{loc}}(\R\times\Sph^{n-1}_+)$, with
the $H^1([0,T]\times\Sph^{n-1}_+)$ norm. We also use
$H^{1/2}(\R\times\partial\Sph^{n-1}_+)$ in
Section~\ref{subsec:s5-l2}. In a boundary chart, let $h$ be a
term localized by a smooth partition of unity, expressed in
the chart and extended by zero to $\R^{n-1}$. Its fractional
seminorm is defined by
\[
 [h]_{H^{1/2}(\R^{n-1})}^2
 =\int_{\R^{n-1}}\int_{\R^{n-1}}
      \frac{|h(x)-h(y)|^2}{|x-y|^n}\,\dd x\,\dd y .
\]
The squared norm on $H^{1/2}(\R\times\partial\Sph^{n-1}_+)$
is the sum of $\|h\|_{L^2(\R^{n-1})}^2+[h]_{H^{1/2}(\R^{n-1})}^2$
over the members of a smooth partition of unity subordinate to
a fixed product atlas of uniformly bounded overlap. The trace theorem
gives a bounded map
\begin{align*}
 H^1(\R\times\Sph^{n-1}_+)
 &\longrightarrow H^{1/2}(\R\times\partial\Sph^{n-1}_+), \\
 z&\longmapsto\tilde z,
\end{align*}
and a bounded right inverse. See~\cite[Theorem~3.37]{McLean}.
We write $H^{-1/2}(\R\times\partial\Sph^{n-1}_+)$ for the dual
of $H^{1/2}(\R\times\partial\Sph^{n-1}_+)$ with respect to the
$L^2(\R\times\partial\Sph^{n-1}_+)$ pairing.

\begin{proposition}\label{prop:periodic-robin}
Let $v$ be a positive smooth $T$-periodic solution of
\eqref{eq:s3-normalized-model} and let $0<\gamma<1$. The map
\[
 \begin{aligned}
 (\mathcal A,\mathcal B_v):
 C_T^{2,\gamma}(\R\times\Sph^{n-1}_+)
 \longrightarrow{}&
 C_T^{0,\gamma}(\R\times\Sph^{n-1}_+)
  {}\times C_T^{1,\gamma}(\R\times\partial\Sph^{n-1}_+)
 \end{aligned}
\]
is Fredholm of index zero. Let $f\in C_T^{0,\gamma}(\R\times\Sph^{n-1}_+)$ and
$g\in C_T^{1,\gamma}(\R\times\partial\Sph^{n-1}_+)$. The equations
$\mathcal A y=f$ and $\mathcal B_v y=g$ have a solution
$y\in C_T^{2,\gamma}(\R\times\Sph^{n-1}_+)$ if and only if
\begin{equation}\label{eq:periodic-robin-compatibility}
 \int_0^T\int_{\Sph^{n-1}_+}f\Phi\,\dd\theta\,\dd t
 +\int_0^T\int_{\partial\Sph^{n-1}_+}g\tilde\Phi\,\dd\sigma\,\dd t
 =0
\end{equation}
for every $T$-periodic Jacobi field $\Phi$ along $v$. When this
condition holds, there is a unique solution orthogonal to these fields in
$L^2([0,T]\times\Sph^{n-1}_+)$.

The same assertions hold in the cylindrically symmetric
subspaces if $v$ is cylindrically symmetric, and in either
parity subspace if $v$ is even in $t$. In each case the fields
in~\eqref{eq:periodic-robin-compatibility} belong to the chosen
subspace.
\end{proposition}

\begin{proof}
For $y,z\in H^1_T(\R\times\Sph^{n-1}_+)$, define
\begin{equation}\label{eq:jacobi-periodic-form}
 \begin{aligned}
 Q_{v,T}(y,z)
 ={}&\int_0^T\int_{\Sph^{n-1}_+}
     \left(y_tz_t+
       \langle\nabla_{\Sph^{n-1}_+}y,\nabla_{\Sph^{n-1}_+}z\rangle
       +\frac{(n-2)^2}{4}yz\right)\dd\theta\,\dd t\\
 &-\frac{nc}{n-2}\int_0^T\int_{\partial\Sph^{n-1}_+}
       \tilde v^{2/(n-2)}\tilde y\,\tilde z\,\dd\sigma\,\dd t .
 \end{aligned}
\end{equation}
Given $y\in H^1_T(\R\times\Sph^{n-1}_+)$, choose smooth
$T$-periodic functions $y_j$ converging to $y$ in
$H^1_T(\R\times\Sph^{n-1}_+)$. Applying
\eqref{eq:jacobi-trace-inequality} on $[0,T]$ to $y_j-y_\ell$
gives
\[
 \|\tilde y_j-\tilde y_\ell\|_{L^2([0,T]\times\partial\Sph^{n-1}_+)}^2
 \leq C\|y_j-y_\ell\|_{H^1([0,T]\times\Sph^{n-1}_+)}^2 .
\]
The right-hand side tends to zero as $j,\ell\to\infty$.
Thus the traces converge in $L^2([0,T]\times\partial\Sph^{n-1}_+)$
to $\tilde y$. Passing to the limit in
\eqref{eq:jacobi-trace-inequality} for $y_j$ gives the same
inequality for $y$. It follows that $Q_{v,T}$ is bounded
below and closed.
Indeed, for $\lambda>0$ sufficiently large, the quadratic form
\[
 Q_{v,T}(y,y)+\lambda\|y\|_{L^2([0,T]\times\Sph^{n-1}_+)}^2
\]
defines a norm squared equivalent to
$\|y\|_{H^1_T(\R\times\Sph^{n-1}_+)}^2$. The compact embedding
into $L^2([0,T]\times\Sph^{n-1}_+)$ gives compact resolvent for
the associated self-adjoint Robin operator.

The Lax--Milgram theorem gives a unique periodic weak solution
for the shifted pair $(\mathcal A+\lambda,\mathcal B_v)$ with
data $(f,g)$. The right-hand side is the continuous functional
\[
 z\longmapsto
 \int_0^T\int_{\Sph^{n-1}_+}fz\,\dd\theta\,\dd t
 +\int_0^T\int_{\partial\Sph^{n-1}_+}g\tilde z\,\dd\sigma\,\dd t,
\]
with $z\in H^1_T(\R\times\Sph^{n-1}_+)$. Robin Schauder regularity gives
$y\in C_T^{2,\gamma}(\R\times\Sph^{n-1}_+)$. Thus the shifted
map is an isomorphism. The inclusion
$C_T^{2,\gamma}(\R\times\Sph^{n-1}_+)
 \hookrightarrow C_T^{0,\gamma}(\R\times\Sph^{n-1}_+)$
is compact, so removing the shift leaves a Fredholm map of
index zero.

Integrating the interior terms in $Q_{v,T}(y,z)$ by parts
in $t$ and in the angular variables gives
\[
 \begin{aligned}
 Q_{v,T}(y,z)
 &=\int_0^T\int_{\Sph^{n-1}_+}
       (\mathcal A y)z\,\dd\theta\,\dd t\\
 &\quad+\int_0^T\int_{\partial\Sph^{n-1}_+}
       (\mathcal B_v y)\tilde z\,\dd\sigma\,\dd t
 +\left[\int_{\Sph^{n-1}_+}y_tz\,\dd\theta\right]_{t=0}^{t=T}.
 \end{aligned}
\]
The last term vanishes because $y_t$ and $z$ are $T$-periodic.
If $\mathcal A y=f$ and $\mathcal B_v y=g$, we therefore obtain
\[
 Q_{v,T}(y,z)
 =\int_0^T\int_{\Sph^{n-1}_+}fz\,\dd\theta\,\dd t
 +\int_0^T\int_{\partial\Sph^{n-1}_+}g\tilde z\,\dd\sigma\,\dd t .
\]
For a $T$-periodic Jacobi field $\Phi$, the same integration by
parts gives $Q_{v,T}(\Phi,y)=0$, since
$\mathcal A\Phi=0$ and $\mathcal B_v\Phi=0$. By symmetry,
$Q_{v,T}(y,\Phi)=0$. Taking $z=\Phi$ in the preceding identity
gives~\eqref{eq:periodic-robin-compatibility}.

Let $\Phi_1,\ldots,\Phi_d$ be an orthonormal basis of the
periodic kernel in $L^2([0,T]\times\Sph^{n-1}_+)$. For each
basis field $\Phi_j$, condition
\eqref{eq:periodic-robin-compatibility} takes the form
\[
 \int_0^T\int_{\Sph^{n-1}_+} f\Phi_j\,\dd\theta\,\dd t
 +\int_0^T\int_{\partial\Sph^{n-1}_+}
       g\tilde\Phi_j\,\dd\sigma\,\dd t=0,
 \qquad j=1,\ldots,d.
\]
By linearity in $\Phi$, these $d$ equations suffice to impose
\eqref{eq:periodic-robin-compatibility} on the entire periodic
kernel. To prove that the equations are independent, suppose
that a linear combination of their left-hand sides, with
coefficients $a_1,\ldots,a_d$, vanishes for every pair $(f,g)$.
Set $\Psi=\sum_{j=1}^d a_j\Phi_j$. Choosing $f=\Psi$
and $g=0$ in that linear combination gives
\[
 0=\int_0^T\int_{\Sph^{n-1}_+}|\Psi|^2\,\dd\theta\,\dd t
   =\sum_{j=1}^d |a_j|^2.
\]
These $d$ independent equations define a subspace of
codimension $d$ consisting of the pairs $(f,g)$ that satisfy
\eqref{eq:periodic-robin-compatibility}. The preceding
integration by parts shows that the range lies in this subspace.
Since the operator has index zero and its kernel has dimension
$d$, the range has the same codimension. Hence the inclusion
is an equality.
Given a solution $y$, set
\[
 y_\perp
 =y-\sum_{j=1}^d
   \left(\int_0^T\int_{\Sph^{n-1}_+}y\Phi_j\,\dd\theta\,\dd t\right)
   \Phi_j .
\]
Each term subtracted from $y$ is a Jacobi field, so $y_\perp$
satisfies the same equations and is orthogonal to the kernel.
The difference of two orthogonal solutions is a Jacobi field
orthogonal to every field in the kernel, and hence vanishes.

If $v$, $f$ and $g$ are cylindrically symmetric, rotating the
angular variable of a solution of the shifted problem gives
a solution of the same problem. Uniqueness of the shifted
solution therefore gives cylindrical symmetry. If $v$ is even
and $f$ and $g$ are both even or both odd in $t$, then
$y(-t,\theta)$ or $-y(-t,\theta)$, respectively, solves the same
shifted problem. Uniqueness again gives the required parity.
Thus the shifted map is an isomorphism on each of the stated
subspaces. Removing the shift is still a compact perturbation,
and the preceding argument applies with the kernel taken in
the chosen subspace.
\end{proof}

\subsection{The normalized flat family}
\label{subsec:companion-family}

We recall the classification and estimates for normalized flat
solutions from~\cite{PartI} used in
Sections~\ref{sec:energy-coordinate} and~\ref{sec:fixed-profile}.
We use the function spaces and growth conditions introduced in
Section~\ref{subsec:linearized-spaces}.

For $n\geq3$, let $\mathcal B_n$ be the set of positive,
cylindrically symmetric functions
$v\in C^\infty(\R\times\Sph^{n-1}_+)$ that
solve~\eqref{eq:s3-normalized-model} and satisfy
\[
 m\leq\tilde v(t)\leq M\qquad\text{for all }t\in\R
\]
for some $m,M>0$.

By~\cite[Theorem~A]{PartI}, for each $\eps\in(0,1)$ there is a
unique $v_\eps\in\mathcal B_n$ such that
\[
 \tilde v_\eps(0)=\min_{t\in\R}\tilde v_\eps(t)=\eps .
\]
The minimum $\eps$ is the \emph{necksize}. The solution $v_\eps$
is even in $t$ and has least positive period $T_\eps$. On
$[0,T_\eps)$, its trace attains its minimum only at $0$ and
its maximum only at $T_\eps/2$. Moreover,
$\tilde v_\eps''(0)>0$ and $\tilde v_\eps''(T_\eps/2)<0$.
Every noncylindrical $v\in\mathcal B_n$ has the form
$v(t,\theta)=v_\eps(t+\tau,\theta)$ for a unique $\eps\in(0,1)$
and some $\tau\in\R$. Two values of $\tau$ give the same solution
if and only if their difference is an integer multiple of
$T_\eps$. We set $v_1=C_n$.

Set
\begin{equation}\label{eq:s4-fixed-period-family}
 v_\eps^{(1)}(s,\theta)=v_\eps(T_\eps s,\theta),
 \qquad (s,\theta)\in\R\times\Sph^{n-1}_+.
\end{equation}
For each integer $k\geq0$ and $0<\gamma<1$, the map
\[
 \begin{aligned}
 (0,1)&\longrightarrow
 (0,\infty)\times
 C^{k,\gamma}_{1,\mathrm{sym,even}}(\R\times\Sph^{n-1}_+),\\
 \eps&\longmapsto (T_\eps,v_\eps^{(1)})
 \end{aligned}
\]
is smooth by~\cite[Theorem~A]{PartI}. The estimates in
\cite[Eqs.~(3.8) and (3.27)]{PartI} give
$\eps\leq\tilde v_\eps(t)\leq M_n$ for all $t\in\R$ and
$\eps C_n(\theta)\leq v_\eps(t,\theta)\leq M_nC_n(\theta)$
on $\R\times\Sph^{n-1}_+$, where $M_n=((n-2)/(2c))^{(n-2)/2}$.
By~\cite[Proposition~3.7(i)]{PartI}, $v_\eps\to C_n$ as
$\eps\to1$ in $C^{k,\gamma}_{\mathrm{uloc}}(\R\times\Sph^{n-1}_+)$
for every such $k$ and $\gamma$.

Let $v\in\mathcal B_n$ be even in $t$. Let $\Phi$ be a
cylindrically symmetric Jacobi field along $v$ that is odd in
$t$ and satisfies~\eqref{eq:jacobi-pointwise-growth} for some
$q\geq0$. If $v$ is periodic, we may verify this growth condition
using any of the equivalent bounds in
Lemma~\ref{lem:jacobi-growth-equivalence}.
By~\cite[Proposition~3.1(iii)]{PartI}, if
$\partial_t\Phi(0,L)=0$, then $\Phi\equiv0$.

We use this uniqueness statement to characterize the odd fields
along $v_\eps$ in Section~\ref{subsec:s4-axial-jacobi}.

Theorem~B of~\cite{PartI} classifies the positive global solutions
of~\eqref{eq:flat-boundary}. If the singularity at $0$ is
removable, the solution is a standard bubble of Li and
Zhu~\cite[Theorem~1.2]{LiZhu}.

\subsection{Geometric estimates and the limiting Pohozaev invariant}
\label{subsec:geometric-estimates}

We use Harnack's inequality to extend the two-sided trace bounds
to closed half-annuli and obtain uniform estimates for the
rescaled solutions. In Emden--Fowler coordinates, these estimates
give the decay of the error terms in the geometric equation.
The resulting Hamiltonian limit determines the Pohozaev invariant
used to identify the necksize in Section~\ref{sec:fixed-profile}.

For a solution $U$ of~\eqref{eq:main-geometric}, set
$A_\kappa=(\kappa/c)^{(n-2)/2}$ and define
\begin{equation}\label{eq:main-change}
 w(t,\theta)=A_\kappa e^{-(n-2)t/2}U(e^{-t}\theta),
 \qquad (t,\theta)\in(-\log\rho,\infty)\times\Sph^{n-1}_+.
\end{equation}
When $\kappa=(n-2)/2$, we have $A_\kappa^{-1}=b_n$
by~\eqref{eq:prelim-normalizing-constants}.

\begin{proposition}\label{prop:geometric-inputs}
Let $U$ be a positive solution of~\eqref{eq:main-geometric}
in $C^\infty(B_\rho^+\setminus\{0\})$, with smoothness up to
$\partial_0B_\rho^+\setminus\{0\}$. Suppose that its trace
$u(\overline x)=U(\overline x,0)$ satisfies
\begin{equation}\label{eq:prelim-trace-bounds}
 c_0|\overline x|^{-\frac{n-2}{2}}
 \leq u(\overline x)
 \leq C_0|\overline x|^{-\frac{n-2}{2}}
 \qquad\text{for }0<|\overline x|<\rho
\end{equation}
for some $C_0, c_0>0$. Then $w$ is bounded above and below by
positive constants on $[t_0,\infty)\times\Sph^{n-1}_+$ for some
$t_0>-\log\rho$. For every $\ell>0$, integer $k\geq0$ and
$0<\gamma<1$, the norm of $w$ in
$C^{k,\gamma}([t-\ell,t+\ell]\times\Sph^{n-1}_+)$ is bounded
independently of sufficiently large $t$.

The function $w$ satisfies~\eqref{eq:proof-defect-system}, with
the error terms in~\eqref{eq:proof-defects} satisfying
\eqref{eq:proof-defect-decay}. The Hamiltonian has the limit
in~\eqref{eq:proof-energy-limit}, and the Pohozaev integral
satisfies~\eqref{eq:proof-pohozaev}.
\end{proposition}

\begin{proof}
\emph{Estimates on half-annuli.}
For $r>0$, define $g_r$ by $(g_r)_{ij}(z)=g_{ij}(rz)$ and set
$v_r(z)=A_\kappa r^{(n-2)/2}U(rz)$ for
$z\in\R^n_+$ with $0<|z|<\rho/r$. On this domain the rescaled
equation is
\begin{equation}\label{eq:proof-rescaled}
 \begin{cases}
 -\Delta_{g_r}v_r+r^2\dfrac{n-2}{4(n-1)}R_g(rz)v_r=0
      & \text{for } z_n>0,\\[3pt]
 \partial_{\nu_{g_r}}v_r+\dfrac{n-2}{2}r h_g(r\overline z,0)v_r
      =c\,v_r(\overline z,0)^{n/(n-2)}
      & \text{for } z_n=0 .
 \end{cases}
\end{equation}
Here we used $\kappa A_\kappa^{-2/(n-2)}=c$.
For fixed $0<a<b$, the trace bounds give uniform positive
upper and lower bounds for $v_r(\overline z,0)$ when
$a/2\leq|\overline z|\leq2b$ and $r$ is sufficiently small. Thus the coefficient of $v_r$ in
the boundary condition written in the form
\[
 \partial_{\nu_{g_r}}v_r+
 \left(\frac{n-2}{2}r h_g(r\overline z,0)
       -c\,v_r(\overline z,0)^{2/(n-2)}\right)v_r=0
\]
is uniformly bounded on this set.

Let $B$ be a ball centered on $\partial\R^n_+$ whose closure
lies in $\{z\in\R^n:a/2<|z|<2b\}$. For $z=(\overline z,z_n)\in B$,
define
\[
 \widehat v_r(\overline z,z_n)=v_r(\overline z,|z_n|),
 \qquad
 (\widehat g_r)_{ij}(\overline z,z_n)
   =(g_r)_{ij}(\overline z,|z_n|),\quad 1\leq i,j\leq n.
\]
By the Fermi identities~\eqref{eq:proof-fermi}, $(g_r)_{nn}=1$
and $(g_r)_{an}=0$ for $1\leq a<n$. Consequently, the surface
measure induced by $g_r$ on $B\cap\{z_n=0\}$ is
\[
 \dd S_{g_r}=\sqrt{\det g_r(\overline z,0)}\,\dd\overline z.
\]
To express the boundary term after reflection, define on $B$
\[
 F_r(z)=
 \sqrt{\det g_r(\overline z,0)}
 \left(\frac{n-2}{2}r h_g(r\overline z,0)
       -c\,v_r(\overline z,0)^{2/(n-2)}\right)
       \operatorname{sgn}(z_n)e_n .
\]
The factor multiplying $\operatorname{sgn}(z_n)e_n$ in $F_r$
is independent of $z_n$. Integrating in $z_n$ on either side
of $z_n=0$ therefore writes the boundary term in the reflected
weak equation as $-\int_B F_r\cdot\nabla(\widehat v_r\psi)\,\dd z$.
Integrating by parts in the two half-balls gives, for every
$\psi\in C_c^\infty(B)$,
\[
 \begin{aligned}
 0={}&\int_B\left(
     \sqrt{\det\widehat g_r}\,\widehat g_r^{-1}\nabla\widehat v_r
     -F_r\widehat v_r\right)\cdot\nabla\psi\,\dd z\\
 &+\int_B r^2\frac{n-2}{4(n-1)} R_g(r\overline z,r|z_n|)
     \sqrt{\det\widehat g_r}\,\widehat v_r\psi\,\dd z
   -\int_B (F_r\cdot\nabla\widehat v_r)\psi\,\dd z .
 \end{aligned}
\]
The gradients in this identity are Euclidean. Since $g_{ij}(0)=\delta_{ij}$, the coefficients of
$\widehat g_r$ converge uniformly on $B$ to those of the
Euclidean metric as $r\to0$. Thus, for all sufficiently small $r$,
\[
 \frac12|\xi|^2
 \leq \sqrt{\det\widehat g_r(z)}\,
       \bigl(\widehat g_r(z)^{-1}\xi\bigr)\cdot\xi
 \leq 2|\xi|^2
\]
for all $z\in B\text{ and }\xi\in\R^n$. The trace bounds~\eqref{eq:prelim-trace-bounds} give a uniform
bound for $F_r$ on $B$. The coefficient of $\widehat v_r\psi$
in the scalar curvature term is also uniformly bounded by
the smoothness of $g$ in the fixed Fermi chart. The two
pieces of $\widehat v_r$ are smooth up to $z_n=0$ and agree
there, so $\widehat v_r\in H^1(B)$.

We may therefore apply the Harnack inequality of Gilbarg and
Trudinger~\cite[Theorem~8.20]{GilbargTrudinger} to the positive
weak solution $\widehat v_r$. For each ball $B_R(z_0)$ such
that the closure of $B_{4R}(z_0)$ is contained in $B$, we obtain
\[
 \sup_{B_R(z_0)}\widehat v_r
 \leq C\inf_{B_R(z_0)}\widehat v_r,
\]
where $C$ depends on $n$, $R$, and the uniform coefficient
bounds, and is independent of $r$. In balls contained in
$\{z_n>0\}$, the same theorem applies to $v_r$ directly
through the interior equation in~\eqref{eq:proof-rescaled}.

Choose a finite cover of $\{z\in\R^n_+:a\leq|z|\leq b\}$
by interior balls and half-balls centered on $\{z_n=0\}$.
Choose the cover so that the larger concentric balls required
for Harnack's inequality lie in $\{z\in\R^n:a/2<|z|<2b\}$,
and so that the larger balls used for the interior equation
also lie in $\{z_n>0\}$. The closed half-annulus is connected,
so each member of the cover can be joined to a member
containing $(a+b)e_1/2$ by a finite sequence of members
with consecutive intersections.

At the boundary point $(a+b)e_1/2$,
\eqref{eq:prelim-trace-bounds} gives
\[
 A_\kappa c_0\left(\frac{a+b}{2}\right)^{-(n-2)/2}
 \leq v_r\left(\frac{a+b}{2}e_1\right)
 \leq A_\kappa C_0\left(\frac{a+b}{2}\right)^{-(n-2)/2}.
\]
On each member of the cover, Harnack's inequality bounds
the value of $v_r$ at any point above and below by fixed
multiples of its value at any other point. At an intersection
of two consecutive members, the bounds obtained in the
first therefore give bounds at a point of the second.
Repeating this argument along the sequence gives, for all
sufficiently small $r$,
\begin{equation}\label{eq:proof-full-annulus}
 0<c_{a,b}\leq v_r(z)\leq C_{a,b}
\end{equation}
for $z\in\R^n_+$ with $a\leq|z|\leq b$. The cover is independent
of $r$, so the constants are uniform in $r$.

Apply the interior H\"older estimate
\cite[Theorem~8.24]{GilbargTrudinger} to $v_r$ on interior
balls and to $\widehat v_r$ on balls centered on $\{z_n=0\}$.
The coefficient bounds are uniform in $r$, and
\eqref{eq:proof-full-annulus} on a larger closed half-annulus
bounds the solutions on these balls. Restricting the estimates
for $\widehat v_r$ to $z_n\geq0$ gives a uniform bound for
$v_r$ in $C^{0,\alpha}(\{z\in\R^n_+:a\leq|z|\leq b\})$
for some $0<\alpha<1$.

We next use the boundary condition in~\eqref{eq:proof-rescaled}
to gain one derivative. The bound just obtained makes
$\partial_{\nu_{g_r}}v_r$ uniformly $\alpha$-H\"older continuous
on $\{(\overline z,0):a\leq|\overline z|\leq b\}$.
On a half-ball $B\cap\R^n_+$ with $B$ centered on $\{z_n=0\}$,
the vector field used in the boundary Schauder estimate is
\begin{equation}\label{eq:proof-schauder-field}
 -v_r(\overline z,0)\lim_{s\to0^+}F_r(\overline z,s)
 =\sqrt{\det g_r(\overline z,0)}\,
   \partial_{\nu_{g_r}}v_r(\overline z,0)e_n,
\end{equation}
where the equality follows from~\eqref{eq:proof-rescaled}.
This field points in the $e_n$ direction and is independent
of $z_n$, so its divergence vanishes. On $z_n=0$, its
component in the outward direction $-e_n$ cancels the
normal component of
$\sqrt{\det g_r}\,g_r^{-1}\nabla_\delta v_r$.

The vector field in~\eqref{eq:proof-schauder-field} is uniformly
$\alpha$-H\"older continuous on the half-ball, the coefficients
of $\sqrt{\det g_r}\,g_r^{-1}$ are uniformly smooth. The scalar
term in the divergence equation is uniformly $\alpha$-H\"older
continuous by the preceding estimate for $v_r$ and the smoothness
of the coefficients. The boundary Schauder estimate
\cite[Theorem~1.1 with weight exponent zero]{TerraciniTortoneVita}
therefore gives a uniform $C^{1,\alpha}$ estimate on a
smaller half-ball. Together with the interior Schauder estimates,
this gives a uniform bound for $v_r$ in
$C^{1,\alpha}(\{z\in\R^n_+:a\leq|z|\leq b\})$.

The first derivative estimates give uniform local Lipschitz
bounds. Hence, for any $0<\gamma<1$, the vector field
in~\eqref{eq:proof-schauder-field} is uniformly
$\gamma$-H\"older continuous on each of the fixed half-balls.
Applying the same Schauder estimate with exponent $\gamma$
gives a uniform bound for $v_r$ in
$C^{1,\gamma}(\{z\in\R^n_+:a\leq|z|\leq b\})$.
By~\eqref{eq:proof-full-annulus}, the range of $v_r$ on each
fixed closed half-annulus lies in a compact subset of
$(0,\infty)$. Thus, if the derivatives of $v_r$ through order
$j$ are uniformly bounded and those of order $j$ are uniformly
$\gamma$-H\"older continuous on a larger half-annulus, the
same bounds hold for the vector field
in~\eqref{eq:proof-schauder-field} and for the scalar term on
the half-balls where they are used. Applying the Schauder
estimates on smaller balls and half-balls gives the bounds
through order $j+1$, with the same exponent $\gamma$.
Iteration gives, for every integer $k\geq0$ and $0<\gamma<1$,
\begin{equation}\label{eq:proof-scaled-regularity}
 \|v_r\|_{C^{k,\gamma}(\{z\in\R^n_+:a\leq|z|\leq b\})}
 \leq C_{a,b,k,\gamma}
\end{equation}
for all sufficiently small $r$.

\medskip
\noindent\emph{Equation in Emden--Fowler coordinates.}
Fix $\ell>0$ and take $t_0>-\log\rho+\ell$. Write $r_0=e^{-t_0}$
and $t=t_0+s$, with $|s|\leq\ell$. By the definitions of $w$
and $v_{r_0}$,
\[
 w(t_0+s,\theta)=e^{-(n-2)s/2}v_{r_0}(e^{-s}\theta).
\]
As $(s,\theta)$ varies in $[-\ell,\ell]\times\Sph^{n-1}_+$,
the point $e^{-s}\theta$ lies in the fixed closed half-annulus
$\{z\in\R^n_+:e^{-\ell}\leq|z|\leq e^\ell\}$. Since
$e^{-(n-2)s/2}$ is bounded above and below by positive constants
for $|s|\leq\ell$, \eqref{eq:proof-full-annulus} gives uniform
positive lower and upper bounds for $w$ on
$[t_0-\ell,t_0+\ell]\times\Sph^{n-1}_+$. The same change of
variables and~\eqref{eq:proof-scaled-regularity} give the uniform
bounds in $C^{k,\gamma}([t_0-\ell,t_0+\ell]\times\Sph^{n-1}_+)$.

For $r=e^{-t}$, the Euclidean polar formula gives
\[
 \left(w_{tt}+\Delta_{\Sph^{n-1}_+}w
             -\frac{(n-2)^2}{4}w\right)(t,\theta)
 =A_\kappa r^{(n+2)/2}\Delta_\delta U(r\theta).
\]
At $\theta=(\vartheta,0)$, with $\vartheta\in\Sph^{n-2}$,
the outward conormal derivative is
\[
 \partial_{\nu_{\Sph^{n-1}_+}}w(t,(\vartheta,0))
 =-A_\kappa r^{n/2}\partial_{x_n}U(r\vartheta,0).
\]
In the Fermi chart, $\nu_g=-\partial_{x_n}$ on
$\partial_0B_\rho^+$. Substituting the two equations
in~\eqref{eq:main-geometric} and using
$\kappa A_\kappa^{-2/(n-2)}=c$, we obtain
\begin{equation}\label{eq:proof-defect-system}
 \begin{cases}
 w_{tt}+\Delta_{\Sph^{n-1}_+}w-\dfrac{(n-2)^2}{4}w=f_g,\\[4pt]
 \partial_{\nu_{\Sph^{n-1}_+}}w
 =c\,\tilde w^{n/(n-2)}+f_{\partial,g}.
 \end{cases}
\end{equation}
The first equation holds in
$(-\log\rho,\infty)\times\operatorname{int}\Sph^{n-1}_+$
and the second on
$(-\log\rho,\infty)\times\partial\Sph^{n-1}_+$.
The interior and boundary error terms are, respectively,
\begin{equation}\label{eq:proof-defects}
 \begin{aligned}
 f_g(t,\theta)
 &=A_\kappa r^{(n+2)/2}
   \left[(\Delta_\delta-\Delta_g)U
        +\frac{n-2}{4(n-1)}R_gU\right](r\theta),\\
 f_{\partial,g}(t,\vartheta)
 &=-\frac{n-2}{2}r h_g(r\vartheta,0)\tilde w(t,\vartheta).
 \end{aligned}
\end{equation}

To estimate $f_g$ on $[t_0-\ell,t_0+\ell]\times\Sph^{n-1}_+$,
put $z=e^{-s}\theta$, so that $r\theta=r_0z$. In these
coordinates, the first identity in~\eqref{eq:proof-defects}
becomes
\[
 \begin{aligned}
 f_g(t_0+s,\theta)=e^{-(n+2)s/2}\bigg(
 &[(\Delta_\delta-\Delta_{g_{r_0}})v_{r_0}](z)\\
 &+\frac{n-2}{4(n-1)}r_0^2R_g(r_0z)v_{r_0}(z)\bigg).
 \end{aligned}
\]
The coefficients of $\Delta_\delta-\Delta_{g_{r_0}}$ are
$O(r_0)$ on the fixed half-annulus by~\eqref{eq:prelim-fermi-expansion}
and the coordinate formula for the Laplace--Beltrami operator.
See Section~\ref{subsec:app-geometric-defects} for the explicit
first-order expansion. The scalar coefficient $r_0^2R_g(r_0z)$ is $O(r_0^2)$.
These bounds also hold for derivatives of every fixed order
in $z$. On the boundary, $r=r_0e^{-s}$ with $|s|\leq\ell$, so
$r h_g(r\vartheta,0)$ and its derivatives in $(s,\vartheta)$
are $O(r_0)$. Combining these bounds
with~\eqref{eq:proof-scaled-regularity} gives
\begin{equation}\label{eq:proof-defect-decay}
 \begin{aligned}
 &\|f_g\|_{C^{k,\gamma}([t_0-\ell,t_0+\ell]\times\Sph^{n-1}_+)}\\
 &\qquad
 +\|f_{\partial,g}\|_{
 C^{k,\gamma}([t_0-\ell,t_0+\ell]\times\partial\Sph^{n-1}_+)}
 \leq C_{\ell,k,\gamma}e^{-t_0}
 \end{aligned}
\end{equation}
for all sufficiently large $t_0$.

\medskip
\noindent\emph{The limiting Pohozaev invariant.}
Differentiate the Hamiltonian~\eqref{eq:proof-hamiltonian}
with $a=c$ and integrate by parts on $\Sph^{n-1}_+$. The
angular gradient term contributes the integral of
$w_t\Delta_{\Sph^{n-1}_+}w$ in the hemisphere and the integral
of $-\partial_t\tilde w\,\partial_{\nu_{\Sph^{n-1}_+}}w$ on
the equator. Combining these terms with the other terms in
the derivative gives
\[
 \begin{aligned}
 |\Sph^{n-2}|\frac{\dd}{\dd t}H(w)(t)
 &=\int_{\Sph^{n-1}_+}
    w_t\left(w_{tt}+\Delta_{\Sph^{n-1}_+}w
                  -\frac{(n-2)^2}{4}w\right)\dd\theta\\
 &\quad+\int_{\partial\Sph^{n-1}_+}
    \partial_t\tilde w\left(
       c\,\tilde w^{n/(n-2)}
       -\partial_{\nu_{\Sph^{n-1}_+}}w\right)\dd\sigma\\
 &=\int_{\Sph^{n-1}_+}w_t f_g\,\dd\theta
   -\int_{\partial\Sph^{n-1}_+}
       \partial_t\tilde w\,f_{\partial,g}\,\dd\sigma,
 \end{aligned}
\]
where the last equality follows from~\eqref{eq:proof-defect-system}.
The uniform bounds for $w_t$ and $\partial_t\tilde w$, together
with~\eqref{eq:proof-defect-decay}, bound the absolute value of
this derivative by $Ce^{-t}$ for all sufficiently large $t$.
Since this bound is integrable, $|\Sph^{n-2}|H(w)(t)$ has a
finite limit as $t\to\infty$. We denote the limit by $H_0$.
Integrating the derivative bound from $t$ to infinity gives
\begin{equation}\label{eq:proof-energy-limit}
 |\Sph^{n-2}|H(w)(t)=H_0+O(e^{-t}).
\end{equation}

For $0<r<\rho$, define the Pohozaev integral in the fixed
Fermi chart by
\[
 \begin{aligned}
 P_\kappa(r,U)
 &=\int_{\partial^+B_r^+}
   \left(\frac{n-2}{2}U\partial_rU
          -\frac r2|\nabla_\delta U|^2
          +r(\partial_rU)^2\right)\dd S\\
 &\quad+\frac{(n-2)r\kappa}{2(n-1)}
       \int_{r\Sph^{n-2}}
          u(\overline x)^{\frac{2(n-1)}{n-2}}\dd\overline S ,
 \end{aligned}
\]
where $\dd S$ and $\dd\overline S$ are the Euclidean surface
measures on $\partial^+B_r^+\subset\R^n$ and
$r\Sph^{n-2}\subset\R^{n-1}$, respectively.

Set $t=-\log r$. Differentiating~\eqref{eq:main-change} at
fixed $\theta$ and decomposing the Euclidean gradient into
its radial and angular parts give
\[
 \begin{aligned}
 r\partial_rU(r\theta)
 &=-A_\kappa^{-1}r^{-(n-2)/2}
        \left(w_t+\frac{n-2}{2}w\right)(t,\theta),\\
 |\nabla_\delta U(r\theta)|^2
 &=(\partial_rU(r\theta))^2
   +A_\kappa^{-2}r^{-n}
        |\nabla_{\Sph^{n-1}_+}w(t,\theta)|^2 .
 \end{aligned}
\]
The terms containing $ww_t$ cancel upon substitution into
the hemispherical integral. Since $\dd S=r^{n-1}\dd\theta$,
this integral gives $A_\kappa^{-2}$ times the interior
term in~\eqref{eq:proof-hamiltonian}.

On the equator,
$u(r\vartheta)=A_\kappa^{-1}r^{-(n-2)/2}\tilde w(t,\vartheta)$
and $\dd\overline S=r^{n-2}\dd\sigma$. The power of $u$
in $P_\kappa$ contributes $r^{-(n-1)}$, which cancels the
factor $r^{n-1}$ from $r\,\dd\overline S$. The identity
$\kappa A_\kappa^{-2/(n-2)}=c$ then gives $A_\kappa^{-2}$
times the boundary term in~\eqref{eq:proof-hamiltonian} with
$a=c$. Adding the two integrals and using~\eqref{eq:proof-energy-limit}
gives
\begin{equation}\label{eq:proof-pohozaev}
 P_\kappa(r,U)
 =A_\kappa^{-2}|\Sph^{n-2}|H(w)(-\log r)
 =A_\kappa^{-2}H_0+O(r).
\end{equation}
The Pohozaev invariant of $U$ at the origin is the limit of
$P_\kappa(r,U)$ as $r\to0$, namely $A_\kappa^{-2}H_0$.
\end{proof}

\section{Cylindrically symmetric Jacobi fields and the Hamiltonian}
\label{sec:energy-coordinate}

Let $n\geq3$. We use the family $v_\eps$ from
Section~\ref{subsec:companion-family} to classify the cylindrically
symmetric Jacobi fields with polynomial growth. We then show that
the derivative of the Hamiltonian with respect to the necksize
never vanishes. The endpoint limits determine the range of the
Hamiltonian and give the classification by energy in the half-space.

\subsection{Jacobi fields and the first variation of the Hamiltonian}
\label{subsec:s4-axial-jacobi}
\label{subsec:s4-pairing}

Differentiation of $v_\eps$ with respect to $t$ and $\eps$
gives the translation field and the necksize derivative.
The uniqueness statements for odd and even fields
in~\cite[Proposition~3.1(iii)]{PartI} show that these two
fields span the space of cylindrically symmetric Jacobi fields
with polynomial growth. We then use a periodic problem to show
that their pairing is nonzero and compute the first variation
of the Hamiltonian.

Fix $\eps\in(0,1)$. With the operators $\mathcal A$ and
$\mathcal B_v$ defined in~\eqref{eq:jacobi-linearized-equation},
the equations for a Jacobi field along $v_\eps$ are
\begin{equation}\label{eq:s4-axial-jacobi-equation}
 \begin{cases}
 \mathcal A\Phi=0
   & \text{in }\R\times\operatorname{int}\Sph^{n-1}_+,\\
 \mathcal B_{v_\eps}\Phi=0
   & \text{on }\R\times\partial\Sph^{n-1}_+.
 \end{cases}
\end{equation}
For $v\in\mathcal B_n$, write
$\mathcal K_{\mathrm{pol}}^{\mathrm{sym}}(v)$ for the set of
$\Phi\in C^\infty_{\mathrm{sym}}(\R\times\Sph^{n-1}_+)$
satisfying $\mathcal A\Phi=0$, $\mathcal B_v\Phi=0$ and
\eqref{eq:jacobi-pointwise-growth} for some $q\geq0$.
By Lemma~\ref{lem:jacobi-growth-equivalence}, the growth condition
is equivalent to membership in
$C^{k,\gamma}_{\mathrm{pol},q}(\R\times\Sph^{n-1}_+)$ for every
integer $k\geq0$ and $0<\gamma<1$.

We use the smooth period-one family
$\eps\mapsto(T_\eps,v_\eps^{(1)})$ from
Section~\ref{subsec:companion-family}. Set $F_\eps=\partial_t v_\eps$ and
$G_\eps=\partial_\eps v_\eps$, with the derivative in $\eps$
taken at fixed $t$. Differentiating
$v_\eps(t,\theta)=v_\eps^{(1)}(t/T_\eps,\theta)$, as
given by~\eqref{eq:s4-fixed-period-family}, yields
\begin{equation}\label{eq:s4-necksize-derivative}
 G_\eps(t,\theta)
 =(\partial_\eps v_\eps^{(1)})(t/T_\eps,\theta)
   -\frac{T_\eps'}{T_\eps}\,tF_\eps(t,\theta).
\end{equation}
Here $\partial_\eps v_\eps^{(1)}$ is taken at fixed $s$.

For $y_1,y_2\in C^\infty_{\mathrm{sym}}(\R\times\Sph^{n-1}_+)$,
define the skew-symmetric pairing
\begin{equation}\label{eq:s4-symplectic-pairing}
 \Omega(y_1,y_2)(t)
 =\int_0^L\omega(\beta)
   (y_1\partial_t y_2-(\partial_t y_1)y_2)(t,\beta)\dd\beta .
\end{equation}
If $y_1$ and $y_2$ satisfy~\eqref{eq:s4-axial-jacobi-equation},
differentiation in $t$ and integration by parts in $\beta$ give
\[
 \frac{\dd}{\dd t}\Omega(y_1,y_2)(t)
 =-\bigl[\omega(y_1\partial_\beta y_2
                 -(\partial_\beta y_1)y_2)\bigr]_0^L.
\]
At $\beta=L$, the equations $\mathcal B_{v_\eps}y_i=0$ give
$\partial_\beta y_i(t,L)=(nc/(n-2))
\tilde v_\eps(t)^{2/(n-2)}y_i(t,L)$ for $i=1,2$. The two
products at this endpoint are therefore equal. At $\beta=0$,
the functions and their first derivatives are bounded and
$\omega(0)=0$. Thus the derivative vanishes, and we write
$\Omega(y_1,y_2)$ for the constant value.

\begin{proposition}\label{s4:axial-polynomial-kernel}
For every $\eps\in(0,1)$,
$F_\eps,G_\eps\in\mathcal K_{\mathrm{pol}}^{\mathrm{sym}}(v_\eps)$.
They satisfy
\begin{equation}\label{eq:s4-axial-parity-evaluations}
 \begin{gathered}
 F_\eps(-t,\theta)=-F_\eps(t,\theta),\qquad
 G_\eps(-t,\theta)=G_\eps(t,\theta),\\
 \partial_tF_\eps(0,L)=\tilde v_\eps''(0)>0,\qquad
 G_\eps(0,L)=1.
 \end{gathered}
\end{equation}
The field $F_\eps$ is $T_\eps$-periodic, and $G_\eps$ satisfies
\eqref{eq:jacobi-pointwise-growth} with $q=1$ and
\begin{equation}\label{eq:s4-necksize-monodromy}
 G_\eps(t+T_\eps,\theta)
 =G_\eps(t,\theta)-T_\eps'F_\eps(t,\theta).
\end{equation}
Moreover,
\begin{equation}\label{eq:s4-axial-polynomial-kernel}
 \mathcal K_{\mathrm{pol}}^{\mathrm{sym}}(v_\eps)
   =\operatorname{span}\{F_\eps,G_\eps\}.
\end{equation}
The odd and even subspaces are spanned by $F_\eps$ and
$G_\eps$, respectively, and $\Omega(F_\eps,G_\eps)\ne0$.
\end{proposition}

\begin{proof}
Differentiating~\eqref{eq:s3-normalized-model} in $t$ and in
$\eps$ gives~\eqref{eq:s4-axial-jacobi-equation} for $F_\eps$
and $G_\eps$. The differentiation in $\eps$ at fixed $t$ is
justified on each compact interval by
\eqref{eq:s4-necksize-derivative}. Since $v_\eps$ is even in
$t$, its $t$-derivative is odd and its $\eps$-derivative is even.
The evaluations in~\eqref{eq:s4-axial-parity-evaluations}
follow from $F_\eps=\partial_t v_\eps$,
$\tilde v_\eps''(0)>0$ and $\tilde v_\eps(0)=\eps$.

Both $F_\eps$ and
$(\partial_\eps v_\eps^{(1)})(t/T_\eps,\theta)$ are
$T_\eps$-periodic in $t$. Formula~\eqref{eq:s4-necksize-derivative}
therefore gives the pointwise bound
\eqref{eq:jacobi-pointwise-growth} with $q=1$ for $G_\eps$.
This proves that both fields belong to
$\mathcal K_{\mathrm{pol}}^{\mathrm{sym}}(v_\eps)$.
Differentiating $v_\eps(t+T_\eps,\theta)=v_\eps(t,\theta)$
in $\eps$ and using the periodicity of $F_\eps$ gives
\eqref{eq:s4-necksize-monodromy}.

Let $\Phi\in\mathcal K_{\mathrm{pol}}^{\mathrm{sym}}(v_\eps)$.
If $\Phi$ is odd in $t$, the field
\[
 \Phi-\frac{\partial_t\Phi(0,L)}{\tilde v_\eps''(0)}F_\eps
\]
belongs to $\mathcal K_{\mathrm{pol}}^{\mathrm{sym}}(v_\eps)$,
is odd in $t$ and has vanishing $t$-derivative at $(0,L)$.
If $\Phi$ is even in $t$, the field $\Phi-\Phi(0,L)G_\eps$
belongs to the same space, is even in $t$ and vanishes at
$(0,L)$. In both cases the difference is zero by
\cite[Proposition~3.1(iii)]{PartI}. This result applies for
$n\geq3$, with $n=3,4$ treated in~\cite[Theorem~C.1]{PartI}.
Since $v_\eps$ is even, the odd and even parts of every
$\Phi\in\mathcal K_{\mathrm{pol}}^{\mathrm{sym}}(v_\eps)$
belong to that space. This proves
\eqref{eq:s4-axial-polynomial-kernel}.

Suppose that $\Omega(F_\eps,G_\eps)=0$. Set
$d=-T_\eps'/T_\eps$ and define
\[
 E_\eps=G_\eps-dtF_\eps .
\]
By~\eqref{eq:s4-necksize-monodromy},
$E_\eps(t+T_\eps,\theta)=E_\eps(t,\theta)$.
Since $G_\eps$ and $tF_\eps$ are even in $t$, we have
\[
 E_\eps\in C^{k,\gamma}_{T_\eps,\mathrm{sym,even}}
   (\R\times\Sph^{n-1}_+)
\]
for every integer $k\geq0$ and $0<\gamma<1$. The boundary operator $\mathcal B_{v_\eps}$ commutes with
multiplication by $t$, while
$\mathcal A(tF_\eps)=t\mathcal AF_\eps-2\partial_tF_\eps$.
Consequently,
\[
 \mathcal AE_\eps=2d\,\partial_tF_\eps,
 \qquad \mathcal B_{v_\eps}E_\eps=0.
\]

Substituting $G_\eps=E_\eps+dtF_\eps$ into
\eqref{eq:s4-symplectic-pairing} cancels the two terms
containing $dtF_\eps\partial_tF_\eps$. Integrate the resulting
identity over $[0,T_\eps]$. Integration by parts in $t$ turns
the integral of $-(\partial_tF_\eps)E_\eps$ into that of
$F_\eps\partial_tE_\eps$, because $F_\eps E_\eps$ takes the
same values at $0$ and $T_\eps$. Thus
\begin{equation}\label{eq:s4-odd-fredholm-condition}
 \int_0^{T_\eps}\int_0^L
   \omega F_\eps(2\partial_tE_\eps+dF_\eps)\dd\beta\,\dd t
 =T_\eps\Omega(F_\eps,G_\eps)=0.
\end{equation}

We now apply Proposition~\ref{prop:periodic-robin} with
right-hand side $2\partial_tE_\eps+dF_\eps$ and zero boundary
data. The right-hand side belongs to
$C^{0,\gamma}_{T_\eps,\mathrm{sym,odd}}(\R\times\Sph^{n-1}_+)$.
The compatibility condition~\eqref{eq:periodic-robin-compatibility}
requires its $L^2([0,T_\eps]\times\Sph^{n-1}_+)$ pairing with
each odd, cylindrically symmetric $T_\eps$-periodic Jacobi field
to vanish. The term involving the boundary data is zero.
Every such field is a multiple of $F_\eps$ by
\eqref{eq:s4-axial-polynomial-kernel} and the parity of $G_\eps$.
It therefore suffices to verify
\[
 \begin{aligned}
 &\int_0^{T_\eps}\int_{\Sph^{n-1}_+}
   F_\eps(2\partial_tE_\eps+dF_\eps)\dd\theta\,\dd t\\
 &\qquad=|\Sph^{n-2}|
   \int_0^{T_\eps}\int_0^L
   \omega F_\eps(2\partial_tE_\eps+dF_\eps)\dd\beta\,\dd t=0.
 \end{aligned}
\]
The factor $|\Sph^{n-2}|$ comes from integration in $\vartheta$,
since both functions depend only on $t$ and $\beta$. The last
equality is~\eqref{eq:s4-odd-fredholm-condition}. We obtain
$y\in C^{2,\gamma}_{T_\eps,\mathrm{sym,odd}}
(\R\times\Sph^{n-1}_+)$ satisfying
\[
 \mathcal Ay=2\partial_tE_\eps+dF_\eps,
 \qquad \mathcal B_{v_\eps}y=0.
\]
Elliptic regularity gives $y\in C^\infty(\R\times\Sph^{n-1}_+)$.

The product rule gives
\[
 \begin{aligned}
 \mathcal A(tE_\eps)
 &=2dt\,\partial_tF_\eps-2\partial_tE_\eps,\\
 \mathcal A\bigl((d/2)t^2F_\eps\bigr)
 &=-2dt\,\partial_tF_\eps-dF_\eps.
 \end{aligned}
\]
Thus the function
\[
 \Phi=y+tE_\eps+(d/2)t^2F_\eps
\]
satisfies $\mathcal A\Phi=0$. It also satisfies
$\mathcal B_{v_\eps}\Phi=0$, since $\mathcal B_{v_\eps}$
commutes with multiplication by functions of $t$. The three
summands are odd in $t$ and have at most quadratic growth,
so $\Phi\in\mathcal K_{\mathrm{pol}}^{\mathrm{sym}}(v_\eps)$.
By~\eqref{eq:s4-axial-polynomial-kernel} and the parity of
$G_\eps$, we have $\Phi\in\operatorname{span}\{F_\eps\}$.
In particular, $\Phi$ is $T_\eps$-periodic.

On the other hand, $y(T_\eps,L)=y(0,L)$ by periodicity,
$F_\eps(T_\eps,L)=F_\eps(0,L)=0$ by periodicity and oddness,
and $E_\eps(T_\eps,L)=E_\eps(0,L)=G_\eps(0,L)=1$. Hence
\[
 \begin{aligned}
 \Phi(T_\eps,L)-\Phi(0,L)
 &=y(T_\eps,L)-y(0,L)
   +T_\eps E_\eps(T_\eps,L)
   +\frac d2 T_\eps^2F_\eps(T_\eps,L)\\
 &=T_\eps>0.
 \end{aligned}
\]
This contradicts the periodicity of $\Phi$ and proves
$\Omega(F_\eps,G_\eps)\ne0$.
\end{proof}

Set $h(\eps)=H(v_\eps)$ for $0<\eps<1$, with $H$ as in
\eqref{eq:s3-hamiltonian}.

\begin{lemma}\label{s4:energy-first-variation}
Let $\Phi\in C^\infty_{\mathrm{sym}}(\R\times\Sph^{n-1}_+)$
satisfy~\eqref{eq:s4-axial-jacobi-equation}. At each fixed $t$,
the first variation of $H$ at $v_\eps$ in the direction $\Phi$ is
\begin{equation}\label{eq:s4-energy-first-variation}
 D H_{v_\eps}[\Phi](t)=\Omega(F_\eps,\Phi).
\end{equation}
In particular,
\begin{equation}\label{eq:s4-energy-neck-derivative}
 h'(\eps)=\Omega(F_\eps,G_\eps)\ne0.
\end{equation}
\end{lemma}

\begin{proof}
Fix $t\in\R$. At $(t,L)$, the boundary value of
$v_\eps+s\Phi$ is $\tilde v_\eps(t)+s\tilde\Phi(t)$.
Since $\tilde v_\eps(t)>0$, this value remains positive
for sufficiently small $|s|$. Holding $t$ fixed, the
chain rule gives the derivative of the boundary term
in~\eqref{eq:s3-hamiltonian}:
\[
 \left.\frac{\dd}{\dd s}\right|_{s=0}
 \left[
   \frac{(n-2)c}{2(n-1)}
   \bigl(\tilde v_\eps(t)+s\tilde\Phi(t)\bigr)^{\frac{2(n-1)}{n-2}}
 \right]
 =c\,\tilde v_\eps(t)^{n/(n-2)}\tilde\Phi(t).
\]
Differentiating the full expression gives
\[
 \begin{aligned}
 D H_{v_\eps}[\Phi](t)
 &=\int_0^L\omega\left(
     \partial_t v_\eps\,\partial_t\Phi
     -\partial_\beta v_\eps\,\partial_\beta\Phi
     -\frac{(n-2)^2}{4}v_\eps\Phi\right)\dd\beta\\
 &\quad+c\,\tilde v_\eps(t)^{n/(n-2)}\tilde\Phi(t).
 \end{aligned}
\]
Integration by parts in $\beta$ gives
\[
 -\int_0^L\omega\,\partial_\beta v_\eps\,
                   \partial_\beta\Phi\dd\beta
 =-\bigl[\omega\,\partial_\beta v_\eps\,\Phi\bigr]_0^L
   +\int_0^L(\omega\,\partial_\beta v_\eps)_\beta\Phi\dd\beta.
\]
The term at $\beta=0$ vanishes because $\omega(0)=0$.
At $\beta=L$, we have $\omega(L)=1$ and
$\partial_\beta v_\eps(t,L)=c\,\tilde v_\eps(t)^{n/(n-2)}$.
The endpoint term is therefore
$-c\,\tilde v_\eps(t)^{n/(n-2)}\tilde\Phi(t)$, which cancels
the contribution from differentiating the boundary term
in~\eqref{eq:s3-hamiltonian}. The interior equation for $v_\eps$ gives
\[
 (\omega\,\partial_\beta v_\eps)_\beta
   -\frac{(n-2)^2}{4}\omega v_\eps
 =-\omega\,\partial_t^2v_\eps .
\]
Substitution into the formula for the first variation yields
\[
 D H_{v_\eps}[\Phi](t)
 =\int_0^L\omega
    (\partial_t v_\eps\,\partial_t\Phi-\partial_t^2v_\eps\,\Phi)
    \dd\beta
 =\Omega(F_\eps,\Phi).
\]
Taking $\Phi=G_\eps=\partial_\eps v_\eps$ gives the derivative
of $h$. Its nonvanishing follows from
Proposition~\ref{s4:axial-polynomial-kernel}.
\end{proof}

For $\Phi\in\mathcal K_{\mathrm{pol}}^{\mathrm{sym}}(v_\eps)$,
write $\Phi=aF_\eps+bG_\eps$ with $a,b\in\R$. By
\eqref{eq:s4-energy-first-variation} and the skew-symmetry
of $\Omega$,
\begin{equation}\label{eq:s4-fixed-energy-axial-kernel}
 D H_{v_\eps}[\Phi]=b\,h'(\eps).
\end{equation}
Thus the first variation vanishes if and only if
$\Phi\in\operatorname{span}\{F_\eps\}$.

\subsection{The range of the Hamiltonian and classification by energy}
\label{subsec:s4-global-energy}
\label{subsec:s4-physical-energy}

We compute the limits of $h$ as $\eps\to0$ and $\eps\to1$.
Together with~\eqref{eq:s4-energy-neck-derivative}, they show
that the energy determines the necksize. We use the energy
as a parameter to describe the bounded Jacobi fields and
then return to the original normalization of the half-space
equation.

\begin{lemma}\label{s4:endpoint-energies}
The centered family satisfies
\begin{equation}\label{eq:s4-energy-endpoints}
 \lim_{\eps\to0}h(\eps)=0,\qquad
 \lim_{\eps\to1}h(\eps)=H_c .
\end{equation}
Moreover,
\begin{equation}\label{eq:s4-zero-neck-local-limit}
 v_\eps\to0\quad\text{in }
 C^1_{\mathrm{loc}}(\R\times\Sph^{n-1}_+)
 \quad\text{as }\eps\to0.
\end{equation}
\end{lemma}

\begin{proof}
By~\cite[Proposition~3.7(i)]{PartI}, $v_\eps\to C_n$ in
$C^{1,\gamma}_{\mathrm{uloc}}(\R\times\Sph^{n-1}_+)$ for
$0<\gamma<1$ as $\eps\to1$. Evaluating
\eqref{eq:s3-hamiltonian} at $t=0$ gives $h(\eps)\to H_c$.

For the other limit, let $(\eps_j)$ be a sequence in $(0,1)$
tending to zero. The bound $v_{\eps_j}\leq M_nC_n$ in
Section~\ref{subsec:companion-family} and the local estimates
of Section~\ref{subsec:linearized-spaces} give uniform
$C^{1,\delta_0}([-R,R]\times\Sph^{n-1}_+)$ bounds for each
$R>0$, for some $0<\delta_0<1$. A diagonal subsequence
therefore converges in
$C^1_{\mathrm{loc}}(\R\times\Sph^{n-1}_+)$ to a nonnegative
function $v_\infty$. Passing to the limit in the interior
equation in the sense of distributions, and in the boundary
condition by $C^1$ convergence, gives the normalized equation
for $v_\infty$. Interior elliptic regularity makes the limit
smooth in $\R\times\operatorname{int}\Sph^{n-1}_+$.

At $(0,L)$, the trace and the outward conormal derivative of
$v_\infty$ both vanish, since their values for $v_{\eps_j}$ are
$\eps_j$ and $c\eps_j^{n/(n-2)}$. If $v_\infty$ were nonzero,
the strong maximum principle would give positivity in the
interior. The boundary point lemma at $(0,L)$ would then
give a strictly negative outward conormal derivative there.
It follows that $v_\infty=0$. Since every sequence tending
to zero admits such a subsequence, this proves
\eqref{eq:s4-zero-neck-local-limit}. Evaluating
\eqref{eq:s3-hamiltonian} at $t=0$ and using this convergence
for the first derivatives and the trace gives
$h(\eps)\to0$ as $\eps\to0$.
\end{proof}

\setcounter{restatedtheorem}{0}
\stepcounter{theorem}
\begin{restatedtheorem}\label{s4:global-energy-coordinate}
The map
\begin{equation}\label{eq:s4-energy-diffeomorphism}
 h:(0,1)\longrightarrow(H_c,0),\qquad h(\eps)=H(v_\eps),
\end{equation}
is a smooth diffeomorphism with $h'(\eps)<0$. Two
noncylindrical solutions in $\mathcal B_n$ have the same
Hamiltonian if and only if they differ by a translation in
$t$. The only solution in $\mathcal B_n$ with Hamiltonian
$H_c$ is $C_n$.

Let $\eps(h)$ denote the inverse of~\eqref{eq:s4-energy-diffeomorphism}
and set
\begin{equation}\label{eq:s4-energy-family}
 v_h=v_{\eps(h)},\qquad T(h)=T_{\eps(h)},\qquad
 \Psi_h=\partial_h v_h,
\end{equation}
where the derivative is taken at fixed $t$. The map
$h\mapsto(T(h),v_h^{(1)})$, where
$v_h^{(1)}(s,\theta)=v_h(T(h)s,\theta)$, is smooth from $(H_c,0)$
into $(0,\infty)\times C^{k,\gamma}_{1,\mathrm{sym,even}}
(\R\times\Sph^{n-1}_+)$ for every integer $k\geq0$ and
$0<\gamma<1$. The field $\Psi_h$ is even in $t$, satisfies
\eqref{eq:jacobi-pointwise-growth} with $q=1$, and obeys
\begin{equation}\label{eq:s4-energy-pairing-normalization}
 \Omega(\partial_t v_h,\Psi_h)=1,\qquad
 \mathcal K_{\mathrm{pol}}^{\mathrm{sym}}(v_h)
 =\operatorname{span}\{\partial_t v_h,\Psi_h\}.
\end{equation}
Moreover,
\begin{equation}\label{eq:s4-energy-monodromy}
 \Psi_h(t+T(h),\theta)
 =\Psi_h(t,\theta)-T'(h)\partial_t v_h(t,\theta).
\end{equation}
\end{restatedtheorem}

\begin{proof}
The smooth dependence of the family on $\eps$ gives
smoothness of $h$, and~\eqref{eq:s4-energy-neck-derivative}
gives $h'(\eps)\ne0$ for every $\eps\in(0,1)$. If $h'$ took
both positive and negative values, its continuity would give
a zero between them. Thus $h'$ has the same sign throughout
$(0,1)$. If it were positive, $h$ would be strictly increasing,
and the endpoint limits in Lemma~\ref{s4:endpoint-energies}
would give $0\leq H_c$. Since $H_c<0$, we have $h'(\eps)<0$.
Strict monotonicity and the endpoint limits give a bijection
onto $(H_c,0)$. The inverse function theorem gives smoothness
of $\eps(h)$.

By Theorem~A of~\cite{PartI}, every noncylindrical member of
$\mathcal B_n$ is a translate of a centered solution with a
unique necksize. Translating in $t$ leaves $H$ unchanged,
so two such solutions have the same energy if and only if
their centered representatives have the same necksize.
Their energies lie in $(H_c,0)$, and $C_n$ has energy $H_c$.
Thus $C_n$ is the only member of $\mathcal B_n$ with that
energy.

The period-one family parametrized by $h$ is the composition
\[
 h\longmapsto\eps(h)
 \longmapsto\bigl(T_{\eps(h)},v_{\eps(h)}^{(1)}\bigr)
 =\bigl(T(h),v_h^{(1)}\bigr).
\]
Both maps in the composition are smooth. At fixed $t$,
the chain rule gives
\[
 \Psi_h=\frac{G_{\eps(h)}}{h'(\eps(h))}.
\]
By Proposition~\ref{s4:axial-polynomial-kernel}, $\Psi_h$
is even in $t$ and has at most linear growth. Since
$h'(\eps(h))\ne0$ and $\partial_t v_h=F_{\eps(h)}$, the same
proposition gives
$\mathcal K_{\mathrm{pol}}^{\mathrm{sym}}(v_h)
 =\operatorname{span}\{\partial_t v_h,\Psi_h\}$.
By~\eqref{eq:s4-energy-neck-derivative},
\[
 \Omega(\partial_t v_h,\Psi_h)
 =\frac{\Omega(F_{\eps(h)},G_{\eps(h)})}{h'(\eps(h))}=1.
\]
Finally, differentiating $v_h(t+T(h),\theta)=v_h(t,\theta)$
in $h$ at fixed $t$ gives
\[
 \Psi_h(t+T(h),\theta)
 +T'(h)\partial_t v_h(t+T(h),\theta)
 =\Psi_h(t,\theta).
\]
The periodicity of $\partial_t v_h$ gives
\eqref{eq:s4-energy-monodromy}.
\end{proof}

For $H_c<h<0$, define
\[
 \mathcal K_{\mathrm{bd}}^{\mathrm{sym}}(v_h)
 =
 \left\{
  \Phi\in\mathcal K_{\mathrm{pol}}^{\mathrm{sym}}(v_h)
  :\sup_{(t,\theta)\in\R\times\Sph^{n-1}_+}|\Phi(t,\theta)|<\infty
 \right\}.
\]

\begin{corollary}\label{s4:bounded-axial-kernel}
For $H_c<h<0$,
\[
 \mathcal K_{\mathrm{bd}}^{\mathrm{sym}}(v_h)
 =\begin{cases}
 \operatorname{span}\{\partial_t v_h\},& T'(h)\ne0,\\
 \operatorname{span}\{\partial_t v_h,\Psi_h\},& T'(h)=0.
 \end{cases}
\]
\end{corollary}

\begin{proof}
Let $\Phi\in\mathcal K_{\mathrm{bd}}^{\mathrm{sym}}(v_h)$.
By~\eqref{eq:s4-energy-pairing-normalization}, write
$\Phi=a\partial_t v_h+b\Psi_h$ with $a,b\in\R$.
If $T'(h)=0$, \eqref{eq:s4-energy-monodromy} gives
$\Psi_h(t+T(h),\theta)=\Psi_h(t,\theta)$. Thus both
$\partial_t v_h$ and $\Psi_h$ are bounded, since each is
continuous on $[0,T(h)]\times\Sph^{n-1}_+$ and
$T(h)$-periodic in $t$.

Suppose that $T'(h)\ne0$. Applying
\eqref{eq:s4-energy-monodromy} at $t,t+T(h),\dots,t+(k-1)T(h)$
and using the periodicity of $\partial_t v_h$ gives
\[
 \Psi_h(t+kT(h),\theta)
 =\Psi_h(t,\theta)-kT'(h)\partial_t v_h(t,\theta)
\]
for every integer $k\geq1$. It follows that
\[
 \Phi(t+kT(h),\theta)-\Phi(t,\theta)
 =-kbT'(h)\partial_t v_h(t,\theta).
\]
Since $v_h$ is noncylindrical, choose $(t_0,\theta_0)$ with
$\partial_t v_h(t_0,\theta_0)\ne0$. The boundedness of $\Phi$
gives
\[
 \begin{aligned}
 k|bT'(h)\partial_t v_h(t_0,\theta_0)|
 &=|\Phi(t_0+kT(h),\theta_0)-\Phi(t_0,\theta_0)|\\
 &\leq 2\sup_{(t,\theta)\in\R\times\Sph^{n-1}_+}|\Phi(t,\theta)|.
 \end{aligned}
\]
The right-hand side is finite and independent of $k$, so
$b=0$. Every multiple of $\partial_t v_h$ is bounded by
periodicity.
\end{proof}

For a positive classical solution $U$ of~\eqref{eq:flat-boundary}
on $\R^n_+\setminus\{0\}$, set
$v_U(t,\theta)=e^{-(n-2)t/2}U(e^{-t}\theta)$ for
$(t,\theta)\in\R\times\Sph^{n-1}_+$. Formula
\eqref{eq:proof-hamiltonian} with $a=(n-2)/2$ gives
\begin{equation}\label{eq:s4-physical-hamiltonian}
 \begin{aligned}
 &|\Sph^{n-2}|H_{(n-2)/2}(v_U)(t)\\
 &\qquad=\frac12\int_{\Sph^{n-1}_+}
 \left((\partial_t v_U)^2-|\nabla_{\Sph^{n-1}_+}v_U|^2
            -\frac{(n-2)^2}{4}v_U^2\right)\dd\theta\\
 &\qquad\quad+\frac{(n-2)^2}{4(n-1)}
 \int_{\partial\Sph^{n-1}_+}
             \tilde v_U^{\frac{2(n-1)}{n-2}}\dd\sigma .
 \end{aligned}
\end{equation}
Since $v_U$ solves~\eqref{eq:flat-cylinder-unnormalized},
Section~\ref{subsec:flat-model} shows that
$H_{(n-2)/2}(v_U)(t)$ is independent of $t$. We denote its
value by $H_{(n-2)/2}(v_U)$. If $v_U=b_nv$ for a solution $v$
of~\eqref{eq:s3-normalized-model}, then
$((n-2)/2)b_n^{2/(n-2)}=c$ gives
\begin{equation}\label{eq:s4-energy-normalization}
 H_{(n-2)/2}(v_U)=b_n^2H(v).
\end{equation}

\begin{corollary}\label{s4:physical-energy-classification}
For every positive classical solution $U$ of~\eqref{eq:flat-boundary}
on $\R^n_+\setminus\{0\}$,
\[
 b_n^2H_c\leq H_{(n-2)/2}(v_U)\leq0.
\]
\begin{enumerate}[label=\textup{(\roman*)}]
 \item $H_{(n-2)/2}(v_U)=0$ if and only if the trace singularity
 at $0$ is removable.
 \item $H_{(n-2)/2}(v_U)=b_n^2H_c$ if and only if
 $U(r\theta)=b_nr^{-(n-2)/2}C_n(\theta)$.
 \item Each value in $[b_n^2H_c,0)$ is attained by a singular
 solution. Two singular solutions have the same value if and
 only if they differ by a dilation
 $U(x)\mapsto\lambda^{(n-2)/2}U(\lambda x)$ with $\lambda>0$.
\end{enumerate}
\end{corollary}

\begin{proof}
By Theorem~B of~\cite{PartI}, a solution with removable trace
singularity is a standard bubble and extends smoothly across
$0$. The definition of $v_U$ and the chain rule give
\[
 |v_U|+|\partial_t v_U|+|\nabla_{\Sph^{n-1}_+}v_U|
 =O(e^{-(n-2)t/2})
\]
uniformly on $\Sph^{n-1}_+$ as $t\to\infty$. Each integral
in~\eqref{eq:s4-physical-hamiltonian} therefore tends to zero.
Since $H_{(n-2)/2}(v_U)(t)$ is independent of $t$, we obtain
$H_{(n-2)/2}(v_U)=0$.

If the trace singularity is nonremovable, Theorem~B
of~\cite{PartI} gives $v_U(t,\theta)=b_nv_\eps(t+\tau,\theta)$
for a unique $\eps\in(0,1]$ and some $\tau\in\R$, with
$v_1=C_n$. By~\eqref{eq:s4-energy-normalization}, its energy
is $b_n^2H(v_\eps)$. Theorem~\ref{s4:global-energy-coordinate}
gives $b_n^2H_c<b_n^2H(v_\eps)<0$ for $0<\eps<1$, while
$\eps=1$ gives $b_n^2H_c$. This proves~\textup{(i)} and
\textup{(ii)}, since the representative $b_nC_n$ corresponds
to $U(r\theta)=b_nr^{-(n-2)/2}C_n(\theta)$.

For an energy in $(b_n^2H_c,0)$, division by $b_n^2$ gives a
number in $(H_c,0)$. The inverse of
\eqref{eq:s4-energy-diffeomorphism} then gives exactly one
$\eps\in(0,1)$. At energy $b_n^2H_c$, the only choice is
$\eps=1$. For each such $\eps$, the representative $b_nv_\eps$
gives a singular solution with the prescribed energy by
Theorem~B of~\cite{PartI}.

If two singular solutions $U_1$ and $U_2$ have the same
energy, their representatives therefore have the form
$b_nv_\eps(t+\tau_1,\theta)$ and $b_nv_\eps(t+\tau_2,\theta)$
with the same $\eps$. By~\eqref{eq:flat-dilation-translation},
the representative of $\lambda^{(n-2)/2}U_1(\lambda x)$ is
$v_{U_1}(t-\log\lambda,\theta)$. Taking
$\lambda=e^{\tau_1-\tau_2}$ gives $v_{U_2}(t,\theta)$ and
hence $U_2(x)=\lambda^{(n-2)/2}U_1(\lambda x)$. Conversely,
such a dilation replaces $t$ by $t-\log\lambda$ in
\eqref{eq:s4-physical-hamiltonian}. Since
$H_{(n-2)/2}(v_{U_1})(t)$ is independent of $t$, the two
solutions have the same energy.
\end{proof}

\section{Jacobi fields without a symmetry assumption}
\label{sec:full-jacobi}

Let $n\geq3$. We use Jacobi fields obtained from tangential
translations and the Kelvin transform to prove that every field
with polynomial growth on $\R\times\Sph^{n-1}_+$ is cylindrically
symmetric. The classification along the noncylindrical solutions
then follows from Section~\ref{sec:energy-coordinate}. For the
cylindrical solution, we solve the eigenvalue problem on
$\Sph^{n-1}_+$ and compute the quadratic term of the Hamiltonian.
We also determine the $L^2$ kernels of the operators on
$\R\times\Sph^{n-1}_+$ and $\R\times\partial\Sph^{n-1}_+$.

Throughout this section, $v$ is either $v_h$ from
\eqref{eq:s4-energy-family}, with $H_c<h<0$, or the cylindrical
solution $C_n$ of~\eqref{eq:s3-normalized-model}. In the first
case let $T=T(h)$. In the second case choose any $T>0$. Then
$v(t+T,\theta)=v(t,\theta)$ on $\R\times\Sph^{n-1}_+$.
By~\eqref{eq:jacobi-linearized-equation}, if $\Phi$ is a Jacobi
field along $v$, then $\Phi(t+\tau,\theta)$ is a Jacobi field
along $v(t+\tau,\theta)$ for any $\tau\in\R$.

We use the spaces $\mathcal K_{\mathrm{pol}}(v)$ and
$\mathcal K_{\mathrm{bd}}(v)$ defined in the introduction.
Lemma~\ref{lem:jacobi-growth-equivalence} characterizes
polynomial growth by $L^2$ bounds on translated intervals.
The spaces $\mathcal K_{\mathrm{pol}}^{\mathrm{sym}}(v)$ and
$\mathcal K_{\mathrm{bd}}^{\mathrm{sym}}(v)$ from
Section~\ref{sec:energy-coordinate} are their respective
intersections with $C^\infty_{\mathrm{sym}}(\R\times\Sph^{n-1}_+)$.

\subsection{Angular estimates and the noncylindrical classification}
\label{subsec:s5-positive-pair}
\label{subsec:s5-nonaxial-coercivity}
\label{subsec:s5-polynomial}

We prove a positive lower bound for the quadratic form associated
with each spherical harmonic of positive degree. The bound follows
from two positive functions obtained by differentiating families
of solutions of~\eqref{eq:flat-boundary} generated by tangential
translations and the Kelvin transform. We then apply the bound to the spherical harmonic coefficients of a
Jacobi field with polynomial growth.

Fix a real spherical harmonic $S_\ell\in C^\infty(\Sph^{n-2})$
of degree $\ell$, normalized by
$\|S_\ell\|_{L^2(\Sph^{n-2})}=1$. Its eigenvalue is
\[
 -\Delta_{\Sph^{n-2}}S_\ell=\lambda_\ell S_\ell,
 \qquad \lambda_\ell=\ell(\ell+n-3).
\]
For a Jacobi field $\Phi$ along $v$ as in
\eqref{eq:jacobi-linearized-equation}, define
\[
 \begin{aligned}
 y_\ell(t,\beta)
   &=\int_{\Sph^{n-2}}\Phi(t,\beta,\vartheta)S_\ell(\vartheta)
       \,\dd\sigma(\vartheta),\\
 Z_\ell(t,\beta)
   &=\sin^{(n-2)/2}\beta\,y_\ell(t,\beta).
 \end{aligned}
\]
The formula for $\Delta_{\Sph^{n-1}_+}$ in
Section~\ref{subsec:flat-model} gives the equation for
$y_\ell$. Multiplying $y_\ell$ by $\sin^{(n-2)/2}\beta$
removes the first derivative in $\beta$ from that equation.
At $\beta=L$, the multiplier equals $1$ and its
$\beta$-derivative vanishes. Thus the boundary condition in
\eqref{eq:jacobi-linearized-equation} is preserved under
multiplication. The first line of the system below holds on
$\R\times(0,L)$, and the second holds at $\beta=L$ for
every $t\in\R$.
\begin{equation}\label{eq:s5-mode-operator}
 \begin{cases}
 \mathcal L_\ell Z_\ell=0,\\[2pt]
 \partial_\beta Z_\ell(t,L)
 =\dfrac{nc}{n-2}\tilde v(t)^{2/(n-2)} Z_\ell(t,L),
 \end{cases}
\end{equation}
where
\[
 \mathcal L_\ell
 =-\partial_t^2-\partial_\beta^2
   +\left(\frac{(n-2)(n-4)}4+\lambda_\ell\right)\csc^2\beta .
\]

For $\ell\geq1$, we describe the behavior as $\beta\to0$.
The value $\beta=0$ corresponds to the point
$(0,\ldots,0,1)\in\Sph^{n-1}_+$. In geodesic normal
coordinates centered at this point, the vector
$\beta\vartheta\in\R^{n-1}$ represents
$\theta=(\sin\beta\,\vartheta,\cos\beta)$.
The restriction to $\Sph^{n-2}$ of any polynomial of degree
less than $\ell$ is orthogonal to $S_\ell$ in
$L^2(\Sph^{n-2})$. Thus, after expanding $\Phi$ in these
coordinates through degree $\ell-1$, only the remainder
contributes to $y_\ell$. The remainder is $O(\beta^\ell)$,
and its derivative along each ray $\beta\vartheta$ is
$O(\beta^{\ell-1})$. Applying the same expansion to
$\partial_t\Phi$ gives
\[
 \begin{aligned}
 y_\ell(t,\beta)&=O(\beta^\ell),\\
 \partial_t y_\ell(t,\beta)&=O(\beta^\ell),\\
 \partial_\beta y_\ell(t,\beta)&=O(\beta^{\ell-1}).
 \end{aligned}
\]
These estimates are uniform on compact intervals in $t$
as $\beta\to0$.

We next construct two positive solutions of the system
\eqref{eq:s5-mode-operator} with $\ell=1$. The function
$c^{(n-2)/2}v$ solves the equation with boundary coefficient
$1$ and has trace bounded above and below by positive
constants on $\R$. Proposition~A.2 and equation~(A.2)
of~\cite{PartI} therefore give
\[
 \left(\frac{\tilde v'(t)}{\tilde v(t)}\right)^2
   +c^2\tilde v(t)^{4/(n-2)}
 \leq\frac{(n-2)^2}{4}.
\]
The second term on the left is positive. It follows that,
for every $t\in\R$,
\begin{equation}\label{eq:s5-trace-slope}
 \left|\frac{\tilde v'(t)}{\tilde v(t)}\right|
 <\frac{n-2}{2}.
\end{equation}

For $(t,\beta)\in\R\times(0,L]$, define
\begin{equation}\label{eq:s5-conformal-amplitudes}
 \begin{aligned}
 h_\pm(t,\beta)
  &=\frac{n-2}{2}v(t,\beta)
     -\cot\beta\,\partial_\beta v(t,\beta)
     \pm\partial_t v(t,\beta),\\
 z_\pm(t,\beta)
  &=e^{\pm t}\sin^{n/2}\beta\,h_\pm(t,\beta).
 \end{aligned}
\end{equation}

\begin{proposition}\label{s5:positive-pair}
The functions $z_\pm$ in~\eqref{eq:s5-conformal-amplitudes}
satisfy the system~\eqref{eq:s5-mode-operator} with $\ell=1$.
Each $h_\pm$ extends to a $T$-periodic function in
$C^\infty_{\mathrm{sym}}(\R\times\Sph^{n+1}_+)$.
There are constants $C_v, c_v>0$ such that, for every
$(t,\beta)\in\R\times[0,L]$,
\begin{equation}\label{eq:s5-lifted-positive-pair}
 \begin{aligned}
 c_v\leq h_+(t,\beta)&\leq C_v,\\
 c_v\leq h_-(t,\beta)&\leq C_v .
 \end{aligned}
\end{equation}
In particular, $z_\pm>0$ on $\R\times(0,L]$.
\end{proposition}

\begin{proof}
Define $U$ on $\R^n_+\setminus\{0\}$ by
$U(r\theta)=b_nr^{-(n-2)/2}v(-\log r,\theta)$. By
\eqref{eq:flat-cylinder-unnormalized} and
\eqref{eq:s3-normalized-model}, $U$ solves
\eqref{eq:flat-boundary}. For $i\in\{1,\ldots,n-1\}$, the
chain rule at $x=r\theta$ gives
\[
 b_n^{-1}r^{(n-2)/2}\partial_{x_i}U(r\theta)
   =-e^t\sin\beta\,h_+(t,\beta)\vartheta_i .
\]
Indeed, $\partial_{x_i}r=\sin\beta\,\vartheta_i$ and
$\partial_{x_i}\beta=r^{-1}\cos\beta\,\vartheta_i$, while
$\partial_{x_i}t=-r^{-1}\sin\beta\,\vartheta_i$. Cylindrical
symmetry eliminates derivatives in $\vartheta$.

For each $x\ne0$, the translates $U(x+a e_i)$ satisfy
\eqref{eq:flat-boundary} in a neighborhood of $x$ for
sufficiently small $|a|$. Differentiating at $a=0$ and
passing to Emden--Fowler coordinates therefore gives a
Jacobi field along $v$. The factor $\vartheta_i$ is the
restriction to $\Sph^{n-2}$ of the linear harmonic function
$\overline x\mapsto x_i$ on $\R^{n-1}$. Thus it is a
spherical harmonic with $\ell=1$ and $\lambda_1=n-2$.
Multiplying its coefficient by $\sin^{(n-2)/2}\beta$
gives $-z_+$, which consequently satisfies
\eqref{eq:s5-mode-operator}.

By~\eqref{eq:flat-kelvin}, the Kelvin transform centered
at $0$ with radius $1$ replaces $v(t,\beta)$ by
$v(-t,\beta)$. Apply the translation calculation to this
solution and then replace $t$ by $-t$ in the resulting
field. This changes $e^t$ to $e^{-t}$ and reverses the sign
of $\partial_t v$. The coefficient of $\vartheta_i$,
multiplied by $\sin^{(n-2)/2}\beta$, is $-z_-$. The
equations for both coefficients also follow from
\cite[Eqs.~(A.9) and~(A.10)]{PartI} by multiplication
by $\sin^{(n-2)/2}\beta$.

To verify smoothness of $h_\pm$ at $\beta=0$, fix
$\vartheta\in\Sph^{n-2}$ and consider the map
\[
 \beta\longmapsto v(t,(\sin\beta\,\vartheta,\cos\beta))
\]
for $\beta$ near zero, including negative $\beta$.
Replacing $\beta$ by $-\beta$ has the same effect as
replacing $\vartheta$ by $-\vartheta$. Cylindrical symmetry
therefore makes this smooth function even in $\beta$.
In particular, $\partial_\beta v(t,0)=0$, and
\[
 \cot\beta\,\partial_\beta v(t,\beta)
  =(\beta\cot\beta)
       \int_0^1\partial_\beta^2v(t,a\beta)\,\dd a.
\]
Both factors on the right are smooth and even in $\beta$
near zero, with $\beta\cot\beta\to1$. Hence $h_\pm$ extends
smoothly and evenly in $\beta$. In geodesic normal coordinates $\overline x\in\R^{n+1}$
centered at $(0,\ldots,0,1)\in\Sph^{n+1}_+$, we have
$\beta=|\overline x|$. The extension is therefore given by
$(t,\overline x)\mapsto h_\pm(t,|\overline x|)$.
It is smooth because the Taylor expansion of $h_\pm$ in
$\beta$ to any order contains only even powers, which
become polynomials in the coordinates of $\overline x$ after setting
$\beta=|\overline x|$. Periodicity follows from the formula for $h_\pm$
and the periodicity of $v$.

For $0<\beta<L$, substitute
$z_\pm=e^{\pm t}\sin^{n/2}\beta\,h_\pm$ into
$\mathcal L_1z_\pm=0$. Differentiating first the factor
$e^{\pm t}$ gives
\begin{equation}\label{eq:s5-conformal-intertwining}
 (\mathcal L_1-1)(\sin^{n/2}\beta\,h_\pm)
    =\pm2\sin^{n/2}\beta\,\partial_th_\pm .
\end{equation}
Since
\[
 \frac{(\sin^{n/2}\beta)''}{\sin^{n/2}\beta}
  =\frac{n(n-2)}4\csc^2\beta-\frac{n^2}{4},
\]
the term with $\csc^2\beta$ cancels. On a function depending
only on $\beta$, the spherical Laplacian on $\Sph^{n+1}_+$
is $\partial_\beta^2+n\cot\beta\,\partial_\beta$.
For $t\in\R$ and $0<\beta<L$, the resulting equation is
\begin{equation}\label{eq:s5-lifted-equation}
 -\partial_t^2h_\pm\mp2\partial_th_\pm
 -\Delta_{\Sph^{n+1}_+}h_\pm
 +\left(\frac{n^2}{4}-1\right)h_\pm=0.
\end{equation}
The smooth extension just proved makes the equation hold
at $\beta=0$ as well. This is the equation in
\cite[Eq.~(A.8)]{PartI}.

At $\beta=L$, the term with $\cot\beta$ in
\eqref{eq:s5-conformal-amplitudes} vanishes. By
\eqref{eq:s5-trace-slope},
\[
 h_\pm(t,L)=\frac{n-2}{2}\tilde v(t)\pm\tilde v'(t)>0 .
\]
For each sign, $h_\pm$ attains a minimum on the compact
set $[0,T]\times\Sph^{n+1}_+$. Periodicity makes it a
minimum on $\R\times\Sph^{n+1}_+$, including when the
minimizing point has $t=0$ or $t=T$. If the minimum were
nonpositive, its angular coordinate would be in
$\operatorname{int}\Sph^{n+1}_+$ because $h_\pm(t,L)>0$.
At a negative minimum, the first derivative in $t$
vanishes and $\partial_t^2h_\pm+\Delta_{\Sph^{n+1}_+}h_\pm
\geq0$. Since $n^2/4-1>0$, the left-hand side of
\eqref{eq:s5-lifted-equation} would be strictly negative.
If the minimum were zero, the strong maximum principle
applied to $-h_\pm\leq0$ in the connected interior would
give $h_\pm\equiv0$, contrary to $h_\pm(t,L)>0$.

For $c_v$, take the smaller of the two positive minima
on $[0,T]\times\Sph^{n+1}_+$. For $C_v$, take the larger of
the two maxima there. Periodicity gives the same bounds
for all $t\in\R$, proving~\eqref{eq:s5-lifted-positive-pair}.
\end{proof}

For $\ell\geq1$, $\lambda_\ell\geq n-2$, so the coefficient
of $\csc^2\beta$ in $\mathcal L_\ell$ is at least
$n(n-2)/4>0$. For $I=\R$ or $I=[0,T]$, define
\begin{equation}\label{eq:s5-sector-form}
 \begin{aligned}
 Q_{\ell,I}[Z]
 &=\int_I\int_0^L
       \bigl(|\partial_t Z|^2+|\partial_\beta Z|^2\bigr)
       \,\dd\beta\,\dd t\\
 &\quad+\left(\frac{(n-2)(n-4)}4+\lambda_\ell\right)
       \int_I\int_0^L\frac{Z^2}{\sin^2\beta}\,\dd\beta\,\dd t\\
 &\quad-\frac{nc}{n-2}\int_I
          \tilde v(t)^{2/(n-2)}Z(t,L)^2\,\dd t .
 \end{aligned}
\end{equation}
We consider this quadratic form on
\[
 \mathcal D_I
  =\left\{Z\in H^1(I\times(0,L)):
       Z|_{\beta=0}=0,\quad
       \frac{Z}{\sin\beta}\in L^2(I\times(0,L))\right\}.
\]
When $I=[0,T]$, we require in addition that the boundary
values at $t=0$ and $t=T$ agree in the Sobolev trace sense.
The trace theorem on $I\times(0,L)$ bounds the
$L^2(I)$ norm of $Z(\cdot,L)$ by the
$H^1(I\times(0,L))$ norm of $Z$.

Every $Z\in\mathcal D_I$ can be approximated by smooth
functions that vanish near $\beta=0$. On $I=\R$, the
approximating functions also have compact support in $t$.
On $I=[0,T]$, they are restrictions of smooth functions
that are $T$-periodic in $t$. The convergence holds in
$H^1(I\times(0,L))$ and, after division by $\sin\beta$,
in $L^2(I\times(0,L))$.

To see this, let $\chi_\delta(\beta)$ be zero for
$\beta\leq\delta$, equal to $1$ for $\beta\geq2\delta$,
and satisfy $0\leq\chi_\delta\leq1$ and
$|\chi_\delta'|\leq C/\delta$. The differences
$Z-\chi_\delta Z$, $(1-\chi_\delta)\partial_t Z$,
$(1-\chi_\delta)\partial_\beta Z$ and
$(1-\chi_\delta)Z/\sin\beta$ tend to zero in
$L^2(I\times(0,L))$ by dominated convergence. The other term
in the derivative of $Z-\chi_\delta Z$ is
$-\chi_\delta' Z$. Its squared $L^2(I\times(0,L))$ norm is
bounded by
\[
 C\int_I\int_\delta^{2\delta}
           \frac{Z^2}{\sin^2\beta}\,\dd\beta\,\dd t,
\]
since $\sin\beta\leq2\delta$ on the support of
$\chi_\delta'$. The integrand is integrable and the
integration region shrinks to $\beta=0$, so the bound tends
to zero.

On $I=\R$, multiplication by cutoffs in $t$ that equal
$1$ on $[-R,R]$, vanish for $|t|\geq2R$, and have derivatives
bounded by $C/R$ gives convergence in the same spaces as
$R\to\infty$. On $I=[0,T]$, equality of the traces at $0$
and $T$ gives a periodic extension in
$H^1_{\mathrm{loc}}(\R\times(0,L))$. In both cases,
extension by reflection across $\beta=L$ and smoothing
by convolution give smooth approximations in
$H^1(I\times(0,L))$. Choose the mollifier with support small
enough to keep the approximations zero near $\beta=0$.
With $\delta$ fixed, $1/\sin\beta$ is bounded
on their support, so the convergence after division by
$\sin\beta$ follows as well.

\begin{proposition}\label{s5:nonaxial-gap}
For the quadratic forms in~\eqref{eq:s5-sector-form},
there is $\lambda_v>0$ independent of $\ell$ such that,
for every $\ell\geq1$ and $Z\in\mathcal D_I$,
\begin{equation}\label{eq:s5-nonaxial-gap}
 Q_{\ell,I}[Z]\geq
       \lambda_v\|Z\|_{L^2(I\times(0,L))}^2.
\end{equation}
\end{proposition}

\begin{proof}
Let $\phi_0=\sqrt{z_+z_-}=\sin^{n/2}\beta\,\sqrt{h_+h_-}$.
It is positive on $\R\times(0,L]$ and $T$-periodic in $t$
by Proposition~\ref{s5:positive-pair}. Differentiating the
product at $\beta=L$ and using~\eqref{eq:s5-mode-operator}
gives $\partial_\beta\phi_0(t,L)=(nc/(n-2))
\tilde v(t)^{2/(n-2)}\phi_0(t,L)$.

On $\R\times(0,L)$, the identity
$\log\phi_0=(\log z_++\log z_-)/2$ and the equations
$\mathcal L_1 z_\pm=0$ give
\begin{equation}\label{eq:s5-positive-supersolution}
 \mathcal L_1\phi_0
   =\frac{\phi_0}{4}
        \left|\nabla_{t,\beta}\log\frac{z_+}{z_-}\right|^2 .
\end{equation}
Here $\nabla_{t,\beta}=(\partial_t,\partial_\beta)$.
Indeed, for each positive smooth function $z$,
\[
 \frac{(\partial_t^2+\partial_\beta^2)z}{z}
 =(\partial_t^2+\partial_\beta^2)\log z
       +|\nabla_{t,\beta}\log z|^2 .
\]
Applying this first to $z_\pm$ and then to $\phi_0$,
the terms linear in $\log z_\pm$ cancel and the
remaining difference is
$|\nabla_{t,\beta}(\log z_+-\log z_-)|^2/4$.

Take a smooth $Z$ from the approximating functions
described after~\eqref{eq:s5-sector-form} and put
$w=Z/\phi_0$. Expanding the derivative of $Z=\phi_0w$
gives
\[
 |\nabla_{t,\beta}Z|^2-\phi_0^2|\nabla_{t,\beta}w|^2
   =\nabla_{t,\beta}\phi_0\cdot
                \nabla_{t,\beta}(\phi_0w^2).
\]
Integrating by parts and using
\eqref{eq:s5-positive-supersolution} gives Picone's identity
\begin{equation}\label{eq:s5-picone}
 \begin{aligned}
 Q_{1,I}[Z]
 &=\int_I\int_0^L\phi_0^2|\nabla_{t,\beta}w|^2
                 \,\dd\beta\,\dd t\\
 &\quad+\frac14\int_I\int_0^L
        \left|\nabla_{t,\beta}\log\frac{z_+}{z_-}\right|^2
             Z^2\,\dd\beta\,\dd t .
 \end{aligned}
\end{equation}
At $\beta=L$, the term from integration by parts is
\[
 \int_I\frac{\partial_\beta\phi_0(t,L)}{\phi_0(t,L)}
                         Z(t,L)^2\,\dd t .
\]
By the boundary condition for $\phi_0$, it cancels the
last term in~\eqref{eq:s5-sector-form}. The vanishing of
$Z$ near $\beta=0$ eliminates the term there. The terms at
$t=0$ and $t=T$ cancel when $I=[0,T]$ by periodicity. When
$I=\R$, $Z$ has compact support in $t$.

Fix $0<\beta\leq L$ and an interval $J=[a,a+T]$. By
\eqref{eq:s5-conformal-amplitudes},
$z_+/z_-=e^{2t}h_+/h_-$, and $h_+/h_-$ is
$T$-periodic. Thus
\[
 \int_J\partial_t\log(z_+/z_-)\,\dd t=2T.
\]
Cauchy--Schwarz gives
$\int_J|\nabla_{t,\beta}\log(z_+/z_-)|^2\,\dd t\geq4T$.
For $t,s\in J$, the fundamental theorem of calculus and
Cauchy--Schwarz give
\[
 |w(t,\beta)|^2
 \leq 2|w(s,\beta)|^2+
       2T\int_J|\partial_tw(r,\beta)|^2\,\dd r .
\]
Multiply by $|\nabla_{s,\beta}\log(z_+/z_-)(s,\beta)|^2$
and integrate in $s$ on $J$. Divide by the integral of
this weight, which is at least $4T$, and integrate in $t$
on $J$. The resulting inequality is
\[
 \begin{aligned}
 \int_J |w(t,\beta)|^2\,\dd t
 &\leq\frac12\int_J
     \left|\nabla_{t,\beta}\log\frac{z_+}{z_-}\right|^2
       |w(t,\beta)|^2\,\dd t\\
 &\quad+2T^2\int_J|\partial_tw(t,\beta)|^2\,\dd t .
 \end{aligned}
\]
By~\eqref{eq:s5-lifted-positive-pair}, the ratio
$\phi_0^2/\sin^n\beta$ lies between $c_v^2$ and $C_v^2$.
Multiplying the last inequality by $\sin^n\beta$ and
integrating in $\beta$ therefore gives
\[
 \begin{aligned}
 \int_J\int_0^L Z^2\,\dd\beta\,\dd t
 &\leq C\int_J\int_0^L
          \phi_0^2|\partial_t w|^2\,\dd\beta\,\dd t\\
 &\quad+C\int_J\int_0^L
          \left|\nabla_{t,\beta}\log\frac{z_+}{z_-}\right|^2
              Z^2\,\dd\beta\,\dd t ,
 \end{aligned}
\]
where $C$ depends on $c_v,C_v,T$. For $I=[0,T]$, take
$J=[0,T]$. For $I=\R$, sum over $J=[jT,(j+1)T]$ for
$j\in\mathbb Z$. In either case,~\eqref{eq:s5-picone}
gives~\eqref{eq:s5-nonaxial-gap} with $\ell=1$
for the smooth approximating functions.

For $Z\in\mathcal D_I$, the approximation following
\eqref{eq:s5-sector-form} gives convergence of the
derivatives in $L^2(I\times(0,L))$, of $Z/\sin\beta$ in the
same space, and of the boundary values at $\beta=L$ in
$L^2(I)$. Since $\tilde v$ is bounded, every term in
\eqref{eq:s5-sector-form} converges along this
approximation. This proves the inequality for
$Z\in\mathcal D_I$ when $\ell=1$. Finally, for
$\ell\geq1$, the boundary terms for $Q_{\ell,I}$ and
$Q_{1,I}$ agree, and
\[
 Q_{\ell,I}[Z]-Q_{1,I}[Z]
   =(\lambda_\ell-(n-2))
         \int_I\int_0^L\frac{Z^2}{\sin^2\beta}\,\dd\beta\,\dd t
   \geq0 .
\]
\end{proof}

\begin{proposition}\label{s5:angular-liouville}
Every $\Phi\in\mathcal K_{\mathrm{pol}}(v)$ is
cylindrically symmetric.
\end{proposition}

\begin{proof}
Fix a spherical harmonic of degree $\ell\geq1$ and
consider its coefficient $Z_\ell$ in
\eqref{eq:s5-mode-operator}. Choose $\chi\in C_c^\infty(\R)$
with $0\leq\chi\leq1$ that equals $1$ on $[-1,1]$ and
vanishes for $|t|\geq2$. For $R\geq1$, write
$\chi_R(t)=\chi(t/R)$. Near $\beta=0$, the estimates for
$y_\ell$ and the product rule for
$Z_\ell=\sin^{(n-2)/2}\beta\,y_\ell$ give
\[
 \begin{aligned}
 Z_\ell(t,\beta)&=O\bigl(\beta^{(n-2)/2+\ell}\bigr),\\
 \partial_t Z_\ell(t,\beta)&=O\bigl(\beta^{(n-2)/2+\ell}\bigr),\\
 \partial_\beta Z_\ell(t,\beta)&=O\bigl(\beta^{(n-4)/2+\ell}\bigr).
 \end{aligned}
\]
The constants are uniform for $t\in[-2R,2R]$. In the last
line, the two terms come from
$(\sin^{(n-2)/2}\beta)'=O(\beta^{(n-4)/2})$
and $\partial_\beta y_\ell=O(\beta^{\ell-1})$.
Since $(n-4)/2+\ell\geq1/2$, both
$\chi_R\partial_\beta Z_\ell$ and
$\chi_RZ_\ell/\sin\beta$ are square integrable near
$\beta=0$, and $\chi_RZ_\ell$ vanishes there.
Thus $\chi_R Z_\ell\in\mathcal D_\R$.

Multiply $\mathcal L_\ell Z_\ell=0$ by $\chi_R^2Z_\ell$
and integrate on $\R\times[\delta,L]$, with
$0<\delta<L$. The integration by parts in $\beta$
gives the terms
\[
 \begin{aligned}
 &\int_\R\chi_R^2 Z_\ell(t,\delta)
                       \partial_\beta Z_\ell(t,\delta)\,\dd t\\
 &\qquad-\int_\R\chi_R^2 Z_\ell(t,L)
                       \partial_\beta Z_\ell(t,L)\,\dd t .
 \end{aligned}
\]
The first is $O(\delta^{n+2\ell-3})$ for fixed $R$,
by the estimates just proved and the support of
$\chi_R$. It tends to zero as $\delta\to0$.
By~\eqref{eq:s5-mode-operator}, the second equals the
boundary term in~\eqref{eq:s5-sector-form} with
$Z=\chi_R Z_\ell$. Expanding
$|\partial_t(\chi_RZ_\ell)|^2$ and letting
$\delta\to0$ now gives
\begin{equation}\label{eq:s5-cutoff-identity}
 Q_{\ell,\R}[\chi_RZ_\ell]
   =\int_\R\int_0^L|\chi_R'|^2 Z_\ell^2\,\dd\beta\,\dd t .
\end{equation}
By~\eqref{eq:s5-nonaxial-gap},
\begin{equation}\label{eq:s5-cutoff-iteration}
 \begin{aligned}
 \int_{-R}^R\int_0^L Z_\ell^2\,\dd\beta\,\dd t
 &\leq \lambda_v^{-1} Q_{\ell,\R}[\chi_RZ_\ell]\\
 &\leq \frac{\|\chi'\|_{L^\infty(\R)}^2}{\lambda_v R^2}
             \int_{-2R}^{2R}\int_0^L Z_\ell^2\,\dd\beta\,\dd t .
 \end{aligned}
\end{equation}

For each $t$, the normalization of $S_\ell$ in
$L^2(\Sph^{n-2})$ and Cauchy--Schwarz in $\vartheta$
give
\[
 \begin{aligned}
 \int_0^L Z_\ell(t,\beta)^2\,\dd\beta
 &=\int_0^L\omega(\beta)y_\ell(t,\beta)^2\,\dd\beta\\
 &\leq\int_{\Sph^{n-1}_+}|\Phi(t,\theta)|^2\,\dd\theta .
 \end{aligned}
\]
The polynomial growth bound for $\Phi$ and
Lemma~\ref{lem:jacobi-growth-equivalence} therefore give
\[
 \int_{-R}^R\int_0^L Z_\ell^2\,\dd\beta\,\dd t
 \leq C R^{2N+1}
\]
for $R\geq1$ and some integer $N\geq0$. Indeed, $[-R,R]$
is covered by at most $2R+2$ intervals of length $2$,
and the squared local $L^2$ norm of $\Phi$ on the
corresponding products with $\Sph^{n-1}_+$ is bounded by
$CR^{2N}$. Applying~\eqref{eq:s5-cutoff-iteration} at
$R,2R,\ldots,2^{j-1}R$ gives
\[
 \begin{aligned}
 \int_{-R}^R\int_0^L Z_\ell^2\,\dd\beta\,\dd t
 &\leq C_jR^{-2j}\int_{-2^jR}^{2^jR}\int_0^L
                               Z_\ell^2\,\dd\beta\,\dd t\\
 &\leq C_j'R^{2N+1-2j}.
 \end{aligned}
\]
Choose an integer $j>N+1/2$ and let $R\to\infty$. The
integral on the left is nonnegative and nondecreasing in
$R$, so it must vanish for every $R$. This proves
$Z_\ell=0$. Since the argument applies to every spherical
harmonic of degree $\ell\geq1$, completeness in
$L^2(\Sph^{n-2})$ and smoothness imply that
$\Phi(t,\beta,\vartheta)$ is independent of $\vartheta$.
\end{proof}

For $y\in C^\infty(\R\times\Sph^{n-1}_+)$ and fixed $t\in\R$,
the Hamiltonian in~\eqref{eq:proof-hamiltonian} satisfies
\begin{equation}\label{eq:s5-full-first-variation}
 |\Sph^{n-2}|D H_v[y](t)
 =\int_{\Sph^{n-1}_+}
        \bigl(v_t y_t-v_{tt}y\bigr)\dd\theta .
\end{equation}
This is the same integration by parts as in the proof of
Lemma~\ref{s4:energy-first-variation}, now on $\Sph^{n-1}_+$.
The boundary term cancels by~\eqref{eq:s3-normalized-model},
and the interior equation replaces
$\Delta_{\Sph^{n-1}_+}v-(n-2)^2v/4$ by $-v_{tt}$.
The differentiation at $s=0$ is valid because the trace of
$v+s y$ stays positive near the fixed section for small $|s|$.
Indeed, $\tilde v$ has a positive lower bound by the definition
of $\mathcal B_n$ in Section~\ref{subsec:companion-family}, and
$\tilde y$ is bounded on each compact product of a closed
interval with $\partial\Sph^{n-1}_+$.

If $y=\Phi$ satisfies~\eqref{eq:jacobi-linearized-equation},
the right-hand side of~\eqref{eq:s5-full-first-variation}
is independent of $t$. Its derivative, by Green's identity,
is the integral over $\partial\Sph^{n-1}_+$ of
\[
 \tilde\Phi\,\partial_{\nu_{\Sph^{n-1}_+}}v_t
 -\partial_t\tilde v\,\partial_{\nu_{\Sph^{n-1}_+}}\Phi .
\]
This integrand vanishes because $v_t$ and $\Phi$ satisfy the
same boundary condition in~\eqref{eq:jacobi-linearized-equation}.

\setcounter{restatedtheorem}{1}
\stepcounter{theorem}
\begin{restatedtheorem}\label{s5:full-polynomial-kernel}
Let $H_c<h<0$, and let $\Psi_h=\partial_hv_h$ be the field in
\eqref{eq:s4-energy-family}. Then
\begin{equation}\label{eq:s5-full-polynomial-kernel}
 \mathcal K_{\mathrm{pol}}(v_h)
 =\operatorname{span}\{\partial_t v_h,\Psi_h\}.
\end{equation}
For $\Phi=a\partial_t v_h+d\Psi_h$, with $a,d\in\R$, the
derivative of $H$ at $v_h$ in the direction $\Phi$ satisfies
\begin{equation}\label{eq:s5-full-linear-energy}
 D H_{v_h}[\Phi](t)=d
\end{equation}
for every $t\in\R$. In particular, the derivative vanishes if and
only if $\Phi\in\operatorname{span}\{\partial_t v_h\}$.
\end{restatedtheorem}

\begin{proof}
By Proposition~\ref{s5:angular-liouville},
$\mathcal K_{\mathrm{pol}}(v_h)=
\mathcal K_{\mathrm{pol}}^{\mathrm{sym}}(v_h)$.
The description of the latter space in
\eqref{eq:s4-energy-pairing-normalization} gives
\eqref{eq:s5-full-polynomial-kernel}. By
Lemma~\ref{s4:energy-first-variation} and
\eqref{eq:s4-energy-pairing-normalization},
\[
 D H_{v_h}[\Phi](t)
 =a\,\Omega(\partial_t v_h,\partial_t v_h)
   +d\,\Omega(\partial_t v_h,\Psi_h)
 =d .
\]
The first term is zero by skew-symmetry, and the pairing in
the second term is equal to $1$.
\end{proof}

\begin{corollary}\label{s5:full-bounded-kernel}
For $H_c<h<0$,
\begin{equation}\label{eq:s5-full-bounded-kernel}
 \mathcal K_{\mathrm{bd}}(v_h)
 =\begin{cases}
 \operatorname{span}\{\partial_t v_h\},& T'(h)\ne0,\\
 \operatorname{span}\{\partial_t v_h,\Psi_h\},& T'(h)=0 .
 \end{cases}
\end{equation}
Every field in $\mathcal K_{\mathrm{bd}}(v_h)$ is $T(h)$-periodic
in $t$.
\end{corollary}

\begin{proof}
Proposition~\ref{s5:angular-liouville} gives
$\mathcal K_{\mathrm{bd}}(v_h)=
\mathcal K_{\mathrm{bd}}^{\mathrm{sym}}(v_h)$.
Corollary~\ref{s4:bounded-axial-kernel} then gives
\eqref{eq:s5-full-bounded-kernel}. The field $\partial_t v_h$
is $T(h)$-periodic, and when $T'(h)=0$,
\eqref{eq:s4-energy-monodromy} gives the same periodicity for
$\Psi_h$.
\end{proof}

\subsection{The cylindrical solution}
\label{subsec:s5-stationary-spectrum}
\label{subsec:s5-stationary-kernel}

We use the eigenvalue problem on $\Sph^{n-1}_+$ with the
boundary condition obtained by linearizing
\eqref{eq:s3-normalized-model} at $C_n$. Its negative
eigenvalue gives two bounded Jacobi fields. The uniqueness
results for even and odd fields then give the classification
of $\mathcal K_{\mathrm{pol}}(C_n)$. We also compute the
quadratic term in the expansion of the Hamiltonian at $C_n$.

On $L^2((0,L),\omega(\beta)\,\dd\beta)$, let
$\mathcal A_{nc/(n-2)}$ be the self-adjoint operator
associated with the following quadratic form, defined on
$H^1_{\mathrm{sym}}(\Sph^{n-1}_+)$,
\[
 \int_0^L\omega(\beta)
   \left(|f'(\beta)|^2+\frac{(n-2)^2}{4}|f(\beta)|^2\right)\dd\beta
 -\frac{nc}{n-2}|f(L)|^2 .
\]
This is the operator of~\cite[Eqs.~(A.21)--(A.23)]{PartI}
with the boundary coefficient $nc/(n-2)$ from
\eqref{eq:jacobi-linearized-equation} at $C_n$.
It has compact resolvent, and its eigenfunctions are smooth
on $\Sph^{n-1}_+$. The operator acts on $0<\beta<L$ by
the first line below, and its domain imposes the boundary
condition on the second line.
\begin{equation}\label{eq:s5-stationary-angular-operator}
 \begin{aligned}
 \mathcal A_{nc/(n-2)}f
   &=-\frac{1}{\omega}(\omega f')'+\frac{(n-2)^2}{4}f,\\
 f'(L)&=\frac{nc}{n-2}f(L).
 \end{aligned}
\end{equation}
At the point $(0,\ldots,0,1)$, corresponding to $\beta=0$,
we use the $H^1_{\mathrm{sym}}(\Sph^{n-1}_+)$ condition on
the hemisphere itself. The weighted integral includes smooth
functions with $f(0)\ne0$ also when $n=3$.

\begin{proposition}\label{s5:stationary-spectrum}
Let $\mu_0$ and $\mu_1$ be the first two eigenvalues of
the operator in~\eqref{eq:s5-stationary-angular-operator}.
Then
\begin{equation}\label{eq:s5-stationary-spectral-gap}
 \mu_0=-\xi_*^2<0<\frac{n^2}{4}<\mu_1,
\end{equation}
where $\xi_*>0$ is the solution of
$\Theta_n(\xi_*)=nc/(n-2)$ in~\cite[Eq.~(3.51)]{PartI}.
Here $\Theta_n(\xi)=A_\xi'(L)$, with $A_\xi$ defined
in~\cite[Eq.~(3.47)]{PartI}. The eigenfunction for $\mu_0$
with value $1$ at $\beta=L$ is $A_{\xi_*}$, which is
positive on $\Sph^{n-1}_+$. Its equation on $0<\beta<L$
and its boundary values at $\beta=L$ are
\begin{equation}\label{eq:s5-stationary-ground-state}
 \begin{aligned}
 -\frac{1}{\omega}(\omega A_{\xi_*}')'
       +\frac{(n-2)^2}{4}A_{\xi_*}
   &=-\xi_*^2 A_{\xi_*},\\
 A_{\xi_*}(L)&=1,\\
 A_{\xi_*}'(L)&=\frac{nc}{n-2}A_{\xi_*}(L).
 \end{aligned}
\end{equation}
\end{proposition}

\begin{proof}
Equations~(3.47) and~(3.51) of~\cite{PartI} give
\eqref{eq:s5-stationary-ground-state}. The Dirichlet problem
for $A_{\xi_*}$ on $\Sph^{n-1}_+$ has boundary value
$1$ and zeroth-order coefficient
$(n-2)^2/4+\xi_*^2>0$. The maximum principle therefore
gives $A_{\xi_*}>0$. By~\cite[Lemma~A.6]{PartI},
$\mu_1>n^2/4$ for $n\geq3$. Thus $-\xi_*^2$ is the
only negative eigenvalue and equals $\mu_0$.
\end{proof}

\setcounter{restatedtheorem}{2}
\stepcounter{theorem}
\begin{restatedtheorem}\label{s5:stationary-kernel}
Let $A_{\xi_*}$ be as in
\eqref{eq:s5-stationary-ground-state}. For each
$\Phi\in\mathcal K_{\mathrm{pol}}(C_n)$ there are
$a,d\in\R$ such that, on $\R\times\Sph^{n-1}_+$,
\begin{equation}\label{eq:s5-stationary-polynomial-kernel}
 \Phi(t,\theta)
   =A_{\xi_*}(\beta)\bigl(a\cos(\xi_*t)+d\sin(\xi_*t)\bigr).
\end{equation}
In particular,
\[
 \begin{aligned}
 \mathcal K_{\mathrm{pol}}(C_n)
 &=\mathcal K_{\mathrm{bd}}(C_n)\\
 &=\operatorname{span}\bigl\{
     A_{\xi_*}(\beta)\cos(\xi_*t),
     A_{\xi_*}(\beta)\sin(\xi_*t)\bigr\}.
 \end{aligned}
\]
\end{restatedtheorem}

\begin{proof}
By~\eqref{eq:s5-stationary-ground-state}, both functions in
the span satisfy~\eqref{eq:jacobi-linearized-equation} at
$C_n$ and are bounded (the cosine field is constructed
after equation~(3.52) in~\cite{PartI}).

By Proposition~\ref{s5:angular-liouville}, each
$\Phi\in\mathcal K_{\mathrm{pol}}(C_n)$ is cylindrically
symmetric. Its even and odd parts in $t$, denoted by
$\Phi_{\mathrm{even}}$ and $\Phi_{\mathrm{odd}}$, belong
to $\mathcal K_{\mathrm{pol}}(C_n)$ because the coefficients
in~\eqref{eq:jacobi-linearized-equation} at $C_n$ are
independent of $t$. Since $A_{\xi_*}(L)=1$, the even field
\[
 \Phi_{\mathrm{even}}(t,\theta)
 -\Phi(0,L)A_{\xi_*}(\beta)\cos(\xi_*t)
\]
vanishes at $(0,L)$, and the odd field
\[
 \Phi_{\mathrm{odd}}(t,\theta)
 -\frac{\partial_t\Phi(0,L)}{\xi_*}A_{\xi_*}(\beta)\sin(\xi_*t)
\]
has zero $t$-derivative at $(0,L)$. The cylindrical
solution $C_n$ is even in $t$ and has trace $1$.
Thus~\cite[Proposition~3.1(iii)]{PartI} applies to both fields
and forces them to vanish. This proves
\eqref{eq:s5-stationary-polynomial-kernel}.
\end{proof}

Fix $y\in C^\infty(\R\times\Sph^{n-1}_+)$ and $t\in\R$.
The trace of $C_n+s y$ on
$\{t\}\times\partial\Sph^{n-1}_+$ is
$1+s\tilde y(t,\vartheta)$, which is positive for small
$|s|$ because $\tilde y(t,\cdot)$ is bounded on
$\partial\Sph^{n-1}_+$. By
\eqref{eq:s5-full-first-variation} and $\partial_t C_n=0$,
\begin{equation}\label{eq:s5-stationary-first-variation}
 D H_{C_n}[y](t)=0 .
\end{equation}
Define the quadratic term for $|\Sph^{n-2}|H$ by
\begin{equation}\label{eq:s5-stationary-quadratic-form}
 \begin{aligned}
 \mathcal Q_c[y](t)
 &=\frac12\int_{\Sph^{n-1}_+}
    \left((\partial_t y)^2-|\nabla_{\Sph^{n-1}_+}y|^2
                -\frac{(n-2)^2}{4}y^2\right)\dd\theta\\
 &\quad+\frac{nc}{2(n-2)}
             \int_{\partial\Sph^{n-1}_+}\tilde y^2\,\dd\sigma .
 \end{aligned}
\end{equation}

\begin{proposition}\label{s5:stationary-energy}
Let $\Phi\in\mathcal K_{\mathrm{pol}}(C_n)$, and let $a,d$
be its coefficients in~\eqref{eq:s5-stationary-polynomial-kernel}.
For every $t\in\R$,
\begin{equation}\label{eq:s5-stationary-positive-energy}
 \mathcal Q_c[\Phi](t)=
 \frac{\xi_*^2}{2}(a^2+d^2)
   \int_{\Sph^{n-1}_+} A_{\xi_*}^2\,\dd\theta .
\end{equation}
In particular, $\mathcal Q_c[\Phi](t)=0$ at any $t$ if and
only if $\Phi=0$.
\end{proposition}

\begin{proof}
Multiply the equation for $A_{\xi_*}$ in
\eqref{eq:s5-stationary-ground-state} by $A_{\xi_*}$ and
integrate by parts on $\Sph^{n-1}_+$. Using the conormal
derivative at $\beta=L$ from the same equation gives
\[
 \begin{aligned}
 &\int_{\Sph^{n-1}_+}
   \left(|\nabla_{\Sph^{n-1}_+}A_{\xi_*}|^2
         +\frac{(n-2)^2}{4}A_{\xi_*}^2\right)\dd\theta\\
 &\quad-\frac{nc}{n-2}\int_{\partial\Sph^{n-1}_+}
                    A_{\xi_*}^2\,\dd\sigma
   =-\xi_*^2\int_{\Sph^{n-1}_+}A_{\xi_*}^2\,\dd\theta .
 \end{aligned}
\]
Substitute~\eqref{eq:s5-stationary-polynomial-kernel} in
\eqref{eq:s5-stationary-quadratic-form} and use this identity.
The result is
\[
 \begin{aligned}
 \mathcal Q_c[\Phi](t)
  &={}\frac{\xi_*^2}{2}\int_{\Sph^{n-1}_+}A_{\xi_*}^2\,\dd\theta\\
    &\quad\times\Bigl[
      \bigl(-a\sin(\xi_*t)+d\cos(\xi_*t)\bigr)^2
       +\bigl(a\cos(\xi_*t)+d\sin(\xi_*t)\bigr)^2
     \Bigr].
 \end{aligned}
\]
The mixed terms cancel, and $\cos^2(\xi_*t)+\sin^2(\xi_*t)=1$
gives~\eqref{eq:s5-stationary-positive-energy}.
\end{proof}

To expand $H(C_n+s y)(t)$, substitute $C_n+s y$ in
\eqref{eq:proof-hamiltonian} with $a=c$. The integral over
$\Sph^{n-1}_+$ is a polynomial in $s$ of degree at most
two. The boundary integral contains the factor
$(1+s\tilde y(t,\vartheta))^{2(n-1)/(n-2)}$.
For $|s|\sup_{\vartheta}|\tilde y(t,\vartheta)|\leq1/2$,
its argument lies in $[1/2,3/2]$ at every $\vartheta$.
The function $x\mapsto x^{2(n-1)/(n-2)}$ has bounded
third derivative on this interval. Taylor's formula at
$x=1$ therefore bounds the error in the boundary integral
by
\[
 C_{n,c}|s|^3
   \int_{\partial\Sph^{n-1}_+}|\tilde y(t,\vartheta)|^3\,\dd\sigma .
\]
The linear term vanishes by
\eqref{eq:s5-stationary-first-variation}. Combining the
quadratic terms gives~\eqref{eq:s5-stationary-quadratic-form}.
Thus, at a fixed $t$, as $s\to0$,
\begin{equation}\label{eq:s5-stationary-energy-taylor}
 \begin{aligned}
 |\Sph^{n-2}|H(C_n+s y)(t)
 &=|\Sph^{n-2}|H(C_n)\\
 &\quad+s^2\mathcal Q_c[y](t)+O(|s|^3).
 \end{aligned}
\end{equation}
If $\tilde y$ is bounded on
$\R\times\partial\Sph^{n-1}_+$, the same choice of
$|s|$ works for all $t$, and the bound on the remainder is
uniform in $t$. This applies to the Jacobi fields in
\eqref{eq:s5-stationary-polynomial-kernel}, whose traces are
linear combinations of $\cos(\xi_*t)$ and $\sin(\xi_*t)$.

\subsection{The \texorpdfstring{$L^2$}{L2} kernels on the cylinder and its boundary}
\label{subsec:s5-l2}

We define the self-adjoint realizations of the linearized
equations on $\R\times\Sph^{n-1}_+$ and on
$\R\times\partial\Sph^{n-1}_+$. The periodicity of bounded Jacobi
fields implies that both $L^2$ kernels are trivial. A sequence
of cutoffs of a periodic Jacobi field then places zero in the
essential spectrum of the operator on $\R\times\Sph^{n-1}_+$.

For $y,z\in H^1(\R\times\Sph^{n-1}_+)$, we define $Q_v(y,z)$
using the integrands in~\eqref{eq:jacobi-periodic-form},
with integration in $t$ over $\R$.
\begin{equation}\label{eq:s5-bulk-form}
 \begin{aligned}
 Q_v(y,z)
 &=\int_\R\int_{\Sph^{n-1}_+}
   \left(y_tz_t+
   \langle\nabla_{\Sph^{n-1}_+}y,\nabla_{\Sph^{n-1}_+}z\rangle
      +\frac{(n-2)^2}{4}yz\right)\dd\theta\,\dd t\\
 &\quad-\frac{nc}{n-2}
   \int_\R\int_{\partial\Sph^{n-1}_+}
   \tilde v^{2/(n-2)}\tilde y\,\tilde z\,\dd\sigma\,\dd t .
 \end{aligned}
\end{equation}
Write $Q_0(y,z)$ for the integral over $\R\times\Sph^{n-1}_+$
in~\eqref{eq:s5-bulk-form}. We write the forms for real
functions and take their sesquilinear extensions on complex
functions when using the spectral theorem.

The constants in~\eqref{eq:jacobi-trace-inequality} are
independent of the interval $I$. Apply that inequality on
$[-R,R]$ and let $R\to\infty$. The resulting inequality
holds for all $y\in H^1(\R\times\Sph^{n-1}_+)$ by the
continuity of the trace map in
Section~\ref{subsec:linearized-spaces}. The trace $\tilde v$
is bounded by periodicity. Choose the coefficient of the
angular derivative term in~\eqref{eq:jacobi-trace-inequality}
small enough that the absolute value of the boundary term
in $Q_v(y,y)$ is bounded by
\[
 \frac12\int_\R\int_{\Sph^{n-1}_+}
             |\nabla_{\Sph^{n-1}_+}y|^2\,\dd\theta\,\dd t
       +C\|y\|_{L^2(\R\times\Sph^{n-1}_+)}^2 .
\]
For a sufficiently large $\lambda>0$, the square root of
\[
 Q_v(y,y)+\lambda\|y\|_{L^2(\R\times\Sph^{n-1}_+)}^2
\]
is thus equivalent to the norm on
$H^1(\R\times\Sph^{n-1}_+)$. In particular, $Q_v$ is bounded
below, and $H^1(\R\times\Sph^{n-1}_+)$ is complete in this
equivalent norm. This proves that $Q_v$ is closed.

Let $\mathcal A_v$ be the self-adjoint operator on
$L^2(\R\times\Sph^{n-1}_+)$ associated with $Q_v$. Its domain
consists of the functions $y\in H^1(\R\times\Sph^{n-1}_+)$
for which there is an $f\in L^2(\R\times\Sph^{n-1}_+)$ such that
\begin{equation}\label{eq:s5-bulk-domain}
 Q_v(y,z)=\langle f,z\rangle_{L^2(\R\times\Sph^{n-1}_+)}
\end{equation}
for every $z\in H^1(\R\times\Sph^{n-1}_+)$. For such $y$,
we set $\mathcal A_vy=f$. By integration by parts,
\eqref{eq:s5-bulk-domain} is the weak formulation of
$\mathcal Ay=f$ with $\mathcal B_vy=0$, with the operators
from~\eqref{eq:jacobi-linearized-equation}.

Let $f\in H^{1/2}(\R\times\partial\Sph^{n-1}_+)$. We denote
by $Ef$ the function in $H^1(\R\times\Sph^{n-1}_+)$ whose
trace on $\R\times\partial\Sph^{n-1}_+$ equals $f$ in the
Sobolev trace sense and that satisfies
\begin{equation}\label{eq:s5-trace-extension}
 Q_0(Ef,z)=0
\end{equation}
for every $z\in H^1(\R\times\Sph^{n-1}_+)$ with
$z|_{\R\times\partial\Sph^{n-1}_+}=0$. In particular, $Ef$
solves $\mathcal A(Ef)=0$ in the weak sense in
$\R\times\operatorname{int}\Sph^{n-1}_+$.

To construct $Ef$, use the right inverse of the trace from
Section~\ref{subsec:linearized-spaces} to choose
$y_f\in H^1(\R\times\Sph^{n-1}_+)$ whose trace equals $f$ and
whose $H^1(\R\times\Sph^{n-1}_+)$ norm is bounded by a constant
times $\|f\|_{H^{1/2}(\R\times\partial\Sph^{n-1}_+)}$.
Since $n\geq3$, the coefficient $(n-2)^2/4$ is positive, so
$Q_0(y,y)$ is bounded below by a positive multiple of
$\|y\|_{H^1(\R\times\Sph^{n-1}_+)}^2$. The functions with trace zero
on $\R\times\partial\Sph^{n-1}_+$ form a closed subspace of
$H^1(\R\times\Sph^{n-1}_+)$. The Lax--Milgram theorem on
this subspace gives a unique $w$ with
$Q_0(w,z)=-Q_0(y_f,z)$ for all $z$ in the same subspace.
Then $Ef=y_f+w$ has trace $f$ and satisfies
\eqref{eq:s5-trace-extension}. Uniqueness also follows
because the difference of two such extensions has trace
zero and its $Q_0$ energy vanishes by
\eqref{eq:s5-trace-extension}.

For any $y\in H^1(\R\times\Sph^{n-1}_+)$ with trace $f$,
the difference $y-Ef$ has trace zero on
$\R\times\partial\Sph^{n-1}_+$. Consequently,
\[
 Q_0(y,y)=Q_0(Ef,Ef)+Q_0(y-Ef,y-Ef)
\]
by~\eqref{eq:s5-trace-extension}. Thus $Ef$ minimizes the
energy among extensions of $f$. For the upper bound below,
use $Q_0(Ef,Ef)\leq Q_0(y_f,y_f)$ and the norm bound for
$y_f$. For the lower bound, continuity of the trace map
bounds the $H^{1/2}(\R\times\partial\Sph^{n-1}_+)$ norm of
$f$ by a constant times the $H^1(\R\times\Sph^{n-1}_+)$ norm
of $Ef$. The squared $H^1(\R\times\Sph^{n-1}_+)$ norm of $Ef$ is
bounded by a constant times $Q_0(Ef,Ef)$. This gives
\begin{equation}\label{eq:s5-dn-form-equivalence}
 \begin{aligned}
 C^{-1}\|f\|_{H^{1/2}(\R\times\partial\Sph^{n-1}_+)}^2
 &\leq Q_0(Ef,Ef)\\
 &\leq C\|f\|_{H^{1/2}(\R\times\partial\Sph^{n-1}_+)}^2 .
 \end{aligned}
\end{equation}

We define the Dirichlet-to-Neumann map
\[
 P:H^{1/2}(\R\times\partial\Sph^{n-1}_+)
      \longrightarrow H^{-1/2}(\R\times\partial\Sph^{n-1}_+)
\]
by $\langle Pf,g\rangle=Q_0(Ef,Eg)$ for every
$g\in H^{1/2}(\R\times\partial\Sph^{n-1}_+)$. Here
$\langle\cdot,\cdot\rangle$ denotes the dual pairing.
By Cauchy--Schwarz for $Q_0$ and
\eqref{eq:s5-dn-form-equivalence},
\[
 |\langle Pf,g\rangle|
 \leq C\|f\|_{H^{1/2}(\R\times\partial\Sph^{n-1}_+)}
       \|g\|_{H^{1/2}(\R\times\partial\Sph^{n-1}_+)} .
\]
To relate this definition to the conormal derivative, take
$z\in H^1(\R\times\Sph^{n-1}_+)$ and write
$z=E(\tilde z)+(z-E(\tilde z))$. The second summand has
trace zero on $\R\times\partial\Sph^{n-1}_+$, so
\eqref{eq:s5-trace-extension} gives
\[
 Q_0(Ef,z)=\langle Pf,\tilde z\rangle .
\]
If $f$ is smooth, integration by parts using
$\mathcal A(Ef)=0$ identifies $Pf$ with
$\partial_{\nu_{\Sph^{n-1}_+}}(Ef)$ on
$\R\times\partial\Sph^{n-1}_+$.

For $f,g\in H^{1/2}(\R\times\partial\Sph^{n-1}_+)$, define
\begin{equation}\label{eq:s5-trace-form}
 \begin{aligned}
 Q_v^\partial(f,g)
 &=\langle Pf,g\rangle-\frac{nc}{n-2}
    \int_\R\int_{\partial\Sph^{n-1}_+}
      \tilde v^{2/(n-2)}fg\,\dd\sigma\,\dd t .
 \end{aligned}
\end{equation}
This form is closed and bounded below on
$L^2(\R\times\partial\Sph^{n-1}_+)$ by
\eqref{eq:s5-dn-form-equivalence} and boundedness of $\tilde v$.
Its self-adjoint operator $\mathcal J_v$ acts by
\begin{equation}\label{eq:s5-trace-domain}
 \mathcal J_v f=Pf-\frac{nc}{n-2}\tilde v^{2/(n-2)}f .
\end{equation}
The right-hand side is first defined in
$H^{-1/2}(\R\times\partial\Sph^{n-1}_+)$. The domain of
$\mathcal J_v$ consists of the functions
$f\in H^{1/2}(\R\times\partial\Sph^{n-1}_+)$ for which this
right-hand side belongs to $L^2(\R\times\partial\Sph^{n-1}_+)$.
Indeed, this is precisely when the functional
$g\mapsto Q_v^\partial(f,g)$ is represented by an element
of $L^2(\R\times\partial\Sph^{n-1}_+)$ under the $L^2$ pairing.

\begin{corollary}\label{s5:l2-kernels}
For the solutions $v$ considered in this section,
\begin{equation}\label{eq:s5-l2-kernels}
 \begin{aligned}
 \ker\mathcal A_v&=\{0\},\\
 \ker\mathcal J_v&=\{0\}.
 \end{aligned}
\end{equation}
The kernels are taken in $L^2(\R\times\Sph^{n-1}_+)$ and
$L^2(\R\times\partial\Sph^{n-1}_+)$, respectively.
\end{corollary}

\begin{proof}
Let $\Phi\in\ker\mathcal A_v$. Equation~\eqref{eq:s5-bulk-domain}
with $f=0$ is the weak formulation of
\eqref{eq:jacobi-linearized-equation}. Interior elliptic
regularity gives smoothness in
$\R\times\operatorname{int}\Sph^{n-1}_+$. Near
$\R\times\partial\Sph^{n-1}_+$, the conormal derivative
is prescribed by $\mathcal B_v\Phi=0$ as a smooth coefficient
times the trace of $\Phi$. Starting from
$\Phi\in H^1(\R\times\Sph^{n-1}_+)$, the trace belongs locally
to $H^{1/2}$ on the boundary, so the local estimate for
the prescribed conormal derivative gives $H^2$ regularity
on smaller products with $\Sph^{n-1}_+$. The improved trace
regularity then allows the estimate to be applied again,
as in the proof of Lemma~\ref{lem:jacobi-growth-equivalence}.
Repeating this step gives
$\Phi\in C^\infty(\R\times\Sph^{n-1}_+)$.

For every $a\in\R$,
\[
 \|\Phi\|_{L^2([a-1,a+1]\times\Sph^{n-1}_+)}
 \leq\|\Phi\|_{L^2(\R\times\Sph^{n-1}_+)} .
\]
By Lemma~\ref{lem:jacobi-growth-equivalence}, this uniform
bound gives $\Phi\in\mathcal K_{\mathrm{bd}}(v)$.
By Corollary~\ref{s5:full-bounded-kernel}, $\Phi$ is
$T(h)$-periodic when $v=v_h$. By
Theorem~\ref{s5:stationary-kernel}, $\Phi$ is
$2\pi/\xi_*$-periodic when $v=C_n$. Denote the corresponding
period by $\tau$. If $\Phi\ne0$, the integral of $|\Phi|^2$ over
$[0,\tau]\times\Sph^{n-1}_+$ is positive. Periodicity gives
the same positive integral over every
$[k\tau,(k+1)\tau]\times\Sph^{n-1}_+$ with $k\in\mathbb Z$,
contrary to $\Phi\in L^2(\R\times\Sph^{n-1}_+)$.
Thus $\Phi=0$.

Let $f\in\ker\mathcal J_v$ and let
$z\in H^1(\R\times\Sph^{n-1}_+)$. Since the trace of $Ef$
equals $f$, the identity $Q_0(Ef,z)=\langle Pf,\tilde z\rangle$
gives
\[
 Q_v(Ef,z)=Q_v^\partial(f,\tilde z)
          =\langle\mathcal J_v f,\tilde z\rangle=0.
\]
Equation~\eqref{eq:s5-bulk-domain} now gives
$Ef\in\ker\mathcal A_v$. Hence $Ef=0$, and its trace is $f=0$.
\end{proof}

We use the Weyl criterion to show that $0$ belongs to the
essential spectrum of $\mathcal A_v$. For a self-adjoint
operator, the criterion states that a real number $\lambda$
belongs to the essential spectrum if and only if there is
a sequence of unit vectors in the operator domain that
converges weakly to zero and whose residuals at $\lambda$
converge to zero in norm. In our case, $\lambda=0$, so it
suffices to find unit vectors $e_j\in D(\mathcal A_v)$ such that
\[
 e_j\rightharpoonup0
 \quad\text{and}\quad
 \|\mathcal A_v e_j\|_{L^2(\R\times\Sph^{n-1}_+)}\to0 .
\]
The first convergence is weak in $L^2(\R\times\Sph^{n-1}_+)$.
It matters because, if $0$ were an isolated eigenvalue of
finite multiplicity, the projections of $e_j$ onto its
eigenspace would converge to zero in norm. The components
orthogonal to the eigenspace would then have norms tending
to $1$. On that orthogonal complement, the spectral gap
around $0$ would give a positive lower bound for
$\|\mathcal A_v e_j\|_{L^2(\R\times\Sph^{n-1}_+)}$.
If $0$ belonged to the resolvent set, the same lower bound
would follow from boundedness of $\mathcal A_v^{-1}$.

Let $\Phi_0=\partial_t v_h$ when $v=v_h$ and
$\Phi_0(t,\theta)=A_{\xi_*}(\beta)\cos(\xi_*t)$ when
$v=C_n$, with $A_{\xi_*}$ as in
\eqref{eq:s5-stationary-ground-state}. These are nonzero
periodic Jacobi fields by Theorems~\ref{s5:full-polynomial-kernel}
and~\ref{s5:stationary-kernel}. Choose a nonzero
$\chi\in C_c^\infty((1,2))$. For sufficiently large $R$, set
\[
 \begin{aligned}
 y_R(t,\theta)&=\chi(t/R)\Phi_0(t,\theta),\\
 e_R&=\frac{y_R}{\|y_R\|_{L^2(\R\times\Sph^{n-1}_+)}} .
 \end{aligned}
\]
Since the cutoff depends only on $t$,
$\mathcal B_vy_R=\chi(t/R)\mathcal B_v\Phi_0=0$.
Integration by parts in~\eqref{eq:s5-bulk-form} shows that
$y_R$ satisfies~\eqref{eq:s5-bulk-domain} with
$f=\mathcal Ay_R$. Thus $y_R\in D(\mathcal A_v)$ and
\[
 \mathcal A_vy_R=-\frac{2}{R}\chi'(t/R)\partial_t\Phi_0
                   -\frac{1}{R^2}\chi''(t/R)\Phi_0 .
\]
The functions $\Phi_0$ and $\partial_t\Phi_0$ are bounded
by periodicity, and the terms on the right are supported in
$[R,2R]\times\Sph^{n-1}_+$. Their
$L^2(\R\times\Sph^{n-1}_+)$ norms are bounded by $CR^{-1/2}$
and $CR^{-3/2}$, respectively.

Choose $[a,b]\subset(1,2)$ and $\eta>0$ such that
$|\chi|\geq\eta$ on $[a,b]$. If $\tau$ is a period of
$\Phi_0$, the interval $[aR,bR]$ contains
$\lfloor(b-a)R/\tau\rfloor$ consecutive intervals of length
$\tau$. Periodicity gives
\[
 \|y_R\|_{L^2(\R\times\Sph^{n-1}_+)}^2
 \geq\eta^2\left\lfloor\frac{(b-a)R}{\tau}\right\rfloor
       \int_0^\tau\int_{\Sph^{n-1}_+}|\Phi_0|^2\,\dd\theta\,\dd t .
\]
The last integral is positive because $\Phi_0\ne0$. For
large $R$, the lower bound for the norm of $y_R$ is therefore
a positive multiple of $R^{1/2}$. Dividing the bounds for
$\mathcal A_vy_R$ by this norm shows that
$\|\mathcal A_ve_R\|_{L^2(\R\times\Sph^{n-1}_+)}=O(R^{-1})$.
For $z\in L^2(\R\times\Sph^{n-1}_+)$, Cauchy--Schwarz and
the support of $e_R$ give
\[
 |\langle e_R,z\rangle_{L^2(\R\times\Sph^{n-1}_+)}|
 \leq\|z\|_{L^2([R,\infty)\times\Sph^{n-1}_+)} .
\]
The right-hand side tends to zero as $R\to\infty$ by square
integrability of $z$. This proves weak convergence to zero.
The Weyl criterion now gives $0$ in the essential spectrum
of $\mathcal A_v$.

\section{Asymptotics at an isolated boundary singularity}
\label{sec:fixed-profile}

Section~\ref{subsec:asymptotic-convergence} proves
Theorem~\ref{intro:asymptotics}, which gives convergence to a fixed
member of the flat family.
Section~\ref{subsec:deformed-boundary-models} studies conformal
deformations of the flat solutions and gives a first-order
obstruction in boundary Fermi coordinates.

\subsection{Asymptotic convergence}
\label{subsec:asymptotic-convergence}

For a noncylindrical limiting solution, we select a translation
parameter on each fixed interval in $t$ by an orthogonality
condition. We use the first variation of the Hamiltonian to
estimate the error and show that the parameter minus $t$ converges.
At the cylindrical solution, the quadratic term of the Hamiltonian
gives the error estimate.

For $0<\eps\leq1$ and $\tau\in\R$, let $U_{\eps,\tau}$ be the
flat solution with boundary coefficient $\kappa$ obtained from
$v_\eps$ in Section~\ref{subsec:companion-family}. For $r>0$
and $\theta\in\Sph^{n-1}_+$, with $A_\kappa$ as in
\eqref{eq:main-change}, we have
\begin{equation}\label{eq:main-comparison}
 U_{\eps,\tau}(r\theta)
 = A_\kappa^{-1}r^{-\frac{n-2}{2}}
             v_\eps(-\log r+\tau,\theta).
\end{equation}

\setcounter{restatedtheorem}{3}
\stepcounter{theorem}
\begin{restatedtheorem}\label{thm:fixed-profile}
Let $n\geq3$ and let $g$ be a smooth metric in boundary Fermi
coordinates centered at $p$, which is represented by $0$.
Let $U$ be a positive smooth solution of~\eqref{eq:main-geometric}
on $B_\rho^+\setminus\{0\}$. Suppose that there are
$c_0,C_0>0$ such that the trace $u(\overline x)=U(\overline x,0)$
satisfies, for $0<|\overline x|<\rho$,
\begin{equation}\label{eq:main-trace-bounds}
 c_0|\overline x|^{-\frac{n-2}{2}}
 \leq u(\overline x)
 \leq C_0|\overline x|^{-\frac{n-2}{2}}.
\end{equation}
Then there are $\delta>0$, $\rho_0\in(0,\rho)$,
$\eps_*\in(0,1]$, $\tau_*\in\R$, and $C>0$ such that
$U_*=U_{\eps_*,\tau_*}$ satisfies
\begin{equation}\label{eq:main-coordinate-estimate}
 \sum_{j=0}^2|x|^{\frac{n-2}{2}+j}
      |D^j(U-U_*)(x)|\leq C|x|^\delta
\end{equation}
for $x\in\R^n_+$ with $0<|x|<\rho_0$. In particular,
$U(x)=(1+O(|x|^\delta))U_*(x)$ as $x\to0$.

The necksize $\eps_*$ is unique. If $\eps_*<1$, the translation
parameter $\tau_*$ is unique modulo $T_{\eps_*}$. If $\eps_*=1$,
all values of $\tau_*$ give the same cylindrical solution.
For the same $U_*$, the estimate holds with every
$\delta'\in(0,\delta]$ after adjusting $C$.
\end{restatedtheorem}

\begin{proof}
Let $w$ be the transform of $U$ in~\eqref{eq:main-change}.
Proposition~\ref{prop:geometric-inputs} gives two-sided positive
bounds and uniform local smooth estimates for $w$, the error
estimate~\eqref{eq:proof-defect-decay}, and the Hamiltonian
limit~\eqref{eq:proof-energy-limit}. We prove that there are
$\eps_*\in(0,1]$, $\tau_*\in\R$, and $\delta>0$ such that,
for all sufficiently large $t$,
\begin{equation}\label{eq:main-cylinder-estimate}
 \left\|w(t+\cdot,\cdot)-v_{\eps_*}(t+\cdot+\tau_*,\cdot)\right\|
       _{C^2([-1,1]\times\Sph^{n-1}_+)}\leq Ce^{-\delta t}.
\end{equation}

\medskip
\noindent\emph{Identification of the limiting necksize.}
Let $t_j\to\infty$. The estimates in
Proposition~\ref{prop:geometric-inputs} give, after passing to a
subsequence, convergence in
$C^k_{\mathrm{loc}}(\R\times\Sph^{n-1}_+)$ for every integer $k\geq0$,
\begin{equation}\label{eq:proof-translated-limit}
 w(t_j+\cdot,\cdot)\longrightarrow v_\infty .
\end{equation}
The convergence holds up to the equator. The positive upper and
lower bounds pass to $v_\infty$. By~\eqref{eq:proof-defect-decay},
the error terms in~\eqref{eq:proof-defect-system} tend to zero
on every compact interval in the translated $t$ variable.
Thus $v_\infty$ solves~\eqref{eq:s3-normalized-model} on the
whole cylinder. For each $s\in\R$, the convergence on the
section at $s$ and~\eqref{eq:proof-energy-limit} give
\begin{equation}\label{eq:proof-common-energy}
 |\Sph^{n-2}|H(v_\infty)(s)=H_0 .
\end{equation}

By the normalization in Section~\ref{subsec:flat-model},
$b_n r^{-(n-2)/2}v_\infty(-\log r,\theta)$ defines a positive
solution of~\eqref{eq:flat-boundary} on
$\R^n_+\setminus\{0\}$. The lower bound for the trace of
$v_\infty$ makes the singularity at $0$ nonremovable.
Theorem~B\textup{(ii)} of~\cite{PartI} therefore gives
$v_\infty=v_\eps(\cdot+\tau,\cdot)$ for some
$\eps\in(0,1]$ and $\tau\in\R$. By
Theorem~\ref{s4:global-energy-coordinate}, $H(v_\eps)$
decreases strictly from $0$ to $H_c$ as $\eps$ increases
from $0$ to $1$, with $v_1=C_n$ and $H(C_n)=H_c$.
It follows from~\eqref{eq:proof-common-energy} that every
translated limit has the same necksize $\eps_*\in(0,1]$,
determined by
\begin{equation}\label{eq:proof-normalized-limit-energy}
 H(v_{\eps_*})=\frac{H_0}{|\Sph^{n-2}|}.
\end{equation}

Write $v_*=v_{\eps_*}$. For every compact interval $I\subset\R$,
as $t\to\infty$ we have
\begin{equation}\label{eq:proof-orbit-entry}
 \inf_{a\in\R}
 \left\|w(t+\cdot,\cdot)-v_*(\cdot+a,\cdot)\right\|
        _{C^2(I\times\Sph^{n-1}_+)}
 \longrightarrow0 .
\end{equation}
Otherwise there is a sequence of translates whose distance
from every $v_*(\cdot+a,\cdot)$ stays bounded below by a positive
constant. A subsequence converges by~\eqref{eq:proof-translated-limit}
and its limit is a translate of $v_*$, a contradiction.

\medskip
\noindent\emph{Choice of the translation parameter.}
Suppose that $0<\eps_*<1$ and put $T=T_{\eps_*}$. Choose
$a_0>\max\{T,1\}$ and set $J=[-a_0,a_0]$. Let
$\chi\in C_c^\infty((-a_0,a_0))$ be nonnegative and positive on
an interval of length greater than $T$. In
$L^2(J\times\Sph^{n-1}_+,\chi(s)\dd s\dd\theta)$, use
\[
 \begin{aligned}
 \langle F,G\rangle_\chi
   &=\int_J\int_{\Sph^{n-1}_+}\chi(s)FG\dd\theta\dd s,\\
 \|F\|_\chi^2&=\langle F,F\rangle_\chi .
 \end{aligned}
\]
Write $\Gamma(a)=v_*(\cdot+a,\cdot)$ for $a\in\R$. This
map is smooth and $T$-periodic. If
$\|\Gamma(a)-\Gamma(b)\|_\chi=0$, the two translates agree
on an interval of length greater than $T$ and hence everywhere
by periodicity. Since $T$ is the least positive period,
$a-b\in T\mathbb Z$. Also $\|\Gamma'(a)\|_\chi>0$ for every
$a\in\R$. If this norm vanished, $\partial_t v_*$ would vanish
on an interval of length greater than $T$ and hence on the
whole cylinder, contrary to $0<\eps_*<1$.

For $F\in L^2(J\times\Sph^{n-1}_+,\chi(s)\dd s\dd\theta)$,
minimize $\|F-\Gamma(a)\|_\chi$ over $a\in[0,T]$.
A minimizer exists by compactness and satisfies
\[
 \langle F-\Gamma(a),\Gamma'(a)\rangle_\chi=0 .
\]
Since $a\mapsto\|F-\Gamma(a)\|_\chi^2$ is $T$-periodic,
a minimizer on $[0,T]$ is also a minimizer on $\R$.
Thus the orthogonality equation holds even when the
minimizer is represented by $0$ or $T$.

Fix $\alpha\in\R$. At $F=\Gamma(\alpha)$ and $a=\alpha$,
the derivative of the orthogonality equation with respect
to $a$, holding $F$ fixed, is
\[
 \left.\frac{\partial}{\partial a}
 \langle F-\Gamma(a),\Gamma'(a)\rangle_\chi
 \right|_{F=\Gamma(\alpha),\,a=\alpha}
 =-\|\Gamma'(\alpha)\|_\chi^2<0.
\]
The implicit function theorem therefore gives, for $F$
in a neighborhood of $\Gamma(\alpha)$, a unique solution
$a$ near $\alpha$. This solution depends smoothly on $F$.

We claim that there is $\eta>0$ such that the minimizing
parameter is unique modulo $T$ whenever
\[
 \inf_{a\in[0,T]}\|F-\Gamma(a)\|_\chi<\eta .
\]
Otherwise there are functions $F_j$ and minimizing
parameters $a_j,b_j\in[0,T]$ such that
$a_j-b_j\notin T\mathbb Z$ and
\[
 \|F_j-\Gamma(a_j)\|_\chi
 =\|F_j-\Gamma(b_j)\|_\chi\longrightarrow0.
\]
By compactness, after passing to a subsequence we have
$a_j\to a_\infty$ and $b_j\to b_\infty$. The triangle
inequality gives $\Gamma(a_\infty)=\Gamma(b_\infty)$.
Hence $b_\infty-a_\infty=kT$ for some $k\in\mathbb Z$.
By periodicity, $b_j-kT$ is also a minimizing parameter,
and both $a_j$ and $b_j-kT$ tend to $a_\infty$.
Since $F_j\to\Gamma(a_\infty)$ in $\|\cdot\|_\chi$, the
local uniqueness given by the implicit function theorem
implies $a_j=b_j-kT$ for all sufficiently large $j$.
This contradicts $a_j-b_j\notin T\mathbb Z$.

By~\eqref{eq:proof-orbit-entry} with $I=J$, and since the
weighted norm is bounded by a constant times the $C^0$ norm
on $J\times\Sph^{n-1}_+$,
\[
 \inf_{a\in[0,T]}
 \|w(t+\cdot,\cdot)-\Gamma(a)\|_\chi\longrightarrow0
 \qquad\text{as }t\to\infty.
\]
Thus, for all sufficiently large $t$, the minimizing
parameter is unique modulo $T$.

At a minimizing parameter $a$, the derivative with respect
to $a$ of the orthogonality equation is
\[
 -\|\Gamma'(a)\|_\chi^2+
 \langle w(t+\cdot,\cdot)-\Gamma(a),\Gamma''(a)\rangle_\chi.
\]
The first term is bounded above by a negative constant,
since $\Gamma'$ never vanishes and is $T$-periodic. The
second tends to zero uniformly over the minimizing parameters,
since $\Gamma''$ is bounded and the minimizing distance tends
to zero. The implicit function theorem therefore gives a
local $C^1$ solution of the orthogonality equation through
each minimizing parameter.

These local solutions give the minimizing parameters for
nearby times. Indeed, if $t_j\to t$, any limit of minimizing
representatives at $t_j$ in $[0,T]$ minimizes the distance
at $t$, by continuity. Uniqueness modulo $T$ then allows
the representatives at $t_j$ to be chosen to converge to
the representative at $t$. Local uniqueness in the implicit
function theorem identifies them with the local $C^1$
solution. Let $a_1$ and $a_2$ be local representatives defined on
intervals $I_1$ and $I_2$ with nonempty intersection.
By uniqueness modulo $T$, we have
$a_2(t)-a_1(t)\in T\mathbb Z$ for every $t\in I_1\cap I_2$.
Since $a_2-a_1$ is continuous on the interval $I_1\cap I_2$
and $T\mathbb Z$ is discrete, there is a fixed $k\in\mathbb Z$
such that $a_2(t)-a_1(t)=kT$ throughout the intersection.
Replacing $a_2$ by $a_2-kT$ leaves the minimizing translate
unchanged and makes the two representatives agree on
$I_1\cap I_2$. Fix a sufficiently large initial time $t_0$ and a real
minimizing representative at $t_0$. By compactness, we can
subdivide each interval $[t_0+k,t_0+k+1]$, with $k$ a
nonnegative integer, into finitely many subintervals, each
contained in an open interval carrying a local $C^1$
minimizing representative. The subdivision points give
a sequence $t_0=s_0<s_1<\cdots$ tending to infinity. For
each $j$, let $a_j$ be a local representative on an open
interval containing $[s_j,s_{j+1}]$.

Add an integer multiple of $T$ to $a_0$ so that $a_0(t_0)$
equals the chosen initial value. Having fixed $a_j$, add
the constant $a_j(s_{j+1})-a_{j+1}(s_{j+1})$ to $a_{j+1}$.
This constant belongs to $T\mathbb Z$ by uniqueness modulo
$T$. The adjusted representatives agree at $s_{j+1}$ and
therefore agree on a neighborhood of $s_{j+1}$, since
their difference is constant on the overlap.

Define $a(t)=a_j(t)$ for $t\in[s_j,s_{j+1}]$, using the
adjusted representatives. Agreement near each junction
makes $a$ a real-valued $C^1$ function on $[t_0,\infty)$,
and $a(t)$ is a minimizing representative at every $t$.

Define $E_t$ on $J\times\Sph^{n-1}_+$. The minimizing property
of $a(t)$ gives the orthogonality to the translation field.
\begin{equation}\label{eq:proof-phase-projection}
 \begin{aligned}
 E_t(s,\theta)&=w(t+s,\theta)-v_*(s+a(t),\theta),\\
 0&=\int_J\int_{\Sph^{n-1}_+}\chi(s)E_t(s,\theta)
          \,\partial_t v_*(s+a(t),\theta)\dd\theta\dd s .
 \end{aligned}
\end{equation}
Define
\begin{equation}\label{eq:proof-distance}
 d(t)=\|E_t\|_{C^2(J\times\Sph^{n-1}_+)} .
\end{equation}
This is a continuous function of $t$. To show that $d(t)\to0$,
choose parameters $a_1(t)$ for which the error in
\eqref{eq:proof-orbit-entry} tends to zero on $J$. The
minimizing property of $a(t)$ and the triangle inequality give
\[
 \|\Gamma(a(t))-\Gamma(a_1(t))\|_\chi
   \leq 2\|w(t+\cdot,\cdot)-\Gamma(a_1(t))\|_\chi
   \longrightarrow0.
\]
By compactness of $[0,T]$ and the characterization of equal
translates, the distance from $a(t)-a_1(t)$ to $T\mathbb Z$
tends to zero. Smoothness and periodicity of $v_*$ then give
convergence of the difference of these translates to zero
in $C^2(J\times\Sph^{n-1}_+)$, so $d(t)\to0$.

Differentiate the second line of~\eqref{eq:proof-phase-projection}
in $t$ and use $\partial_t w(t+s,\theta)=\partial_s w(t+s,\theta)$.
Writing $w(t+s,\theta)=E_t(s,\theta)+v_*(s+a(t),\theta)$ and
integrating the term with $\partial_s E_t$ by parts gives
\[
 \begin{aligned}
 &(a'(t)-1)\int_J\int_{\Sph^{n-1}_+}
    \chi(s)\bigl((\partial_t v_*)^2-E_t\partial_t^2 v_*\bigr)
          \dd\theta\dd s\\
 &\qquad=-\int_J\int_{\Sph^{n-1}_+}
    \chi'(s)E_t\partial_t v_*\dd\theta\dd s .
 \end{aligned}
\]
Here the derivatives of $v_*$ are evaluated at $(s+a(t),\theta)$,
and $\chi$ vanishes near the endpoints of $J$. The integral
of $\chi(s)(\partial_t v_*)^2$ is
$\|\Gamma'(a(t))\|_\chi^2$. The function
$a\mapsto\|\Gamma'(a)\|_\chi^2$ has a positive minimum on
$[0,T]$. The term with $E_t$ in the same integral is
bounded by $Cd(t)$ and tends to zero. The right-hand side is
also bounded by $Cd(t)$. Thus there is $t_0$ such that, for
$t\geq t_0$,
\begin{equation}\label{eq:proof-phase-derivative}
 |a'(t)-1|\leq Cd(t).
\end{equation}

\medskip
\noindent\emph{An estimate for the error along a noncylindrical solution.}\enspace
We compare $d(t)$ with its supremum on an interval centered at
$t$. If the estimate below fails, a subsequence of the rescaled
errors converges to a Jacobi field. The Hamiltonian and the
orthogonality condition then force this field to vanish.

We claim that there are $R,C>0$ and $t_1\geq t_0$ such that,
for $t\geq t_1+R$,
\begin{equation}\label{eq:proof-window}
 d(t)\leq\frac12\sup_{|s-t|\leq R}d(s)+Ce^{-(t-R)}.
\end{equation}
If the estimate failed, taking $R=C=j$ and $t_1=t_0+j$ would
give $t_j\geq t_0+2j$ such that, with
$D_j=\sup_{|s-t_j|\leq j}d(s)$,
\begin{equation}\label{eq:proof-window-failure}
 d(t_j)>\frac12D_j+j e^{-(t_j-j)}.
\end{equation}
Since $t_j-j\to\infty$ and $d(t)\to0$, we have $D_j\to0$.
Moreover, $D_j\geq d(t_j)>0$, so
\begin{equation}\label{eq:proof-window-ratios}
 \frac{d(t_j)}{D_j}>\frac12,
 \qquad
 \frac{e^{-(t_j-j)}}{D_j}<\frac1j.
\end{equation}
Set $a_j=a(t_j)$ and define
\begin{equation}\label{eq:proof-normalized-fields}
 y_j(s,\theta)=
 \frac{w(t_j+s,\theta)-v_*(s+a_j,\theta)}{D_j}
\end{equation}
for $s>-\log\rho-t_j$ and $\theta\in\Sph^{n-1}_+$.
After passing to a subsequence, the representatives of $a_j$
in $[0,T)$ converge to some $a_\infty\in[0,T]$. We continue
to use $a_j=a(t_j)$ when estimating differences of the
real-valued function $a$.

Fix $\xi\in\R$. For $j>|\xi|$, the interval with endpoints
$t_j$ and $t_j+\xi$ lies in $[t_j-j,t_j+j]$. Integrating
\eqref{eq:proof-phase-derivative} on this interval gives
\begin{equation}\label{eq:proof-phase-growth}
 |a(t_j+\xi)-a_j-\xi|\leq C|\xi|D_j.
\end{equation}
On $J\times\Sph^{n-1}_+$, the $C^2$ norm of
$w(t_j+\xi+s,\theta)-v_*(s+a(t_j+\xi),\theta)$ is at most
$D_j$ by~\eqref{eq:proof-distance}. By
\eqref{eq:proof-phase-growth} and the uniform bounds for
derivatives of $v_*$ through order three, the $C^2$ norm of
$v_*(s+a(t_j+\xi),\theta)-v_*(s+\xi+a_j,\theta)$ is at most
$C|\xi|D_j$. Both translation parameters are fixed when
taking derivatives in $(s,\theta)$. Thus
\begin{equation}\label{eq:proof-jacobi-growth}
 \|y_j(\xi+\cdot,\cdot)\|_{C^2(J\times\Sph^{n-1}_+)}
 \leq C(1+|\xi|),
\end{equation}
for $|\xi|<j$, with $C$ independent of $\xi$ and $j$.

Subtract the equation~\eqref{eq:s3-normalized-model} for
$v_*(\cdot+a_j,\cdot)$ from~\eqref{eq:proof-defect-system} for
$w(t_j+\cdot,\cdot)$. With the operators in
\eqref{eq:jacobi-linearized-equation}, we obtain
\begin{align}
 \mathcal A y_j
 &=-D_j^{-1}f_g(t_j+\cdot,\cdot),
 \label{eq:proof-jacobi-bulk}\\
 \mathcal B_{v_*(\cdot+a_j,\cdot)}y_j
 &=D_j^{-1}f_{\partial,g}(t_j+\cdot,\cdot)+\mathcal R_j.
 \label{eq:proof-jacobi-boundary}
\end{align}
on the translated cylinder and its lateral boundary, respectively.
At a point $(s,\vartheta)$ on that boundary, Taylor's formula
for the power in~\eqref{eq:s3-normalized-model} gives
\begin{equation}\label{eq:proof-taylor}
 \begin{split}
 \mathcal R_j(s,\vartheta)
 &=\frac{2nc}{(n-2)^2}D_j\widetilde y_j(s,\vartheta)^2\\
 &\quad\times\int_0^1(1-\lambda)
   \bigl(\tilde v_*(s+a_j)+\lambda D_j\widetilde y_j(s,\vartheta)
   \bigr)^{(4-n)/(n-2)}\dd\lambda .
 \end{split}
\end{equation}
For $0\leq\lambda\leq1$, the expression in parentheses
in~\eqref{eq:proof-taylor} is the convex combination
\[
 (1-\lambda)\tilde v_*(s+a_j)
 +\lambda\tilde w(t_j+s,\vartheta).
\] 
Proposition~\ref{prop:geometric-inputs}
and positivity and periodicity of $v_*$ bound it above and
away from zero, uniformly on compact intervals in $s$.
By~\eqref{eq:proof-jacobi-growth}, the traces $\widetilde y_j$
have bounded $C^2$ norms on these intervals. Differentiating
the products and the power in~\eqref{eq:proof-taylor} therefore
gives a bound $CD_j$ for the $C^{1,\gamma}$ norm of $\mathcal R_j$
on each compact product with $\partial\Sph^{n-1}_+$, for
$0<\gamma<1$.

By~\eqref{eq:proof-defect-decay} and~\eqref{eq:proof-window-ratios},
the geometric errors in~\eqref{eq:proof-jacobi-bulk}
and~\eqref{eq:proof-jacobi-boundary} tend to zero in
$C^{0,\gamma}(I\times\Sph^{n-1}_+)$ and
$C^{1,\gamma}(I\times\partial\Sph^{n-1}_+)$, respectively,
for each compact interval $I\subset\R$. The operator $\mathcal B_{v_*(\cdot+a_j,\cdot)}$
in~\eqref{eq:proof-jacobi-boundary} is the boundary operator
defined in~\eqref{eq:jacobi-linearized-equation}, evaluated
along the translate $v_*(\cdot+a_j,\cdot)$. On
$I\times\partial\Sph^{n-1}_+$ it acts by
\[
 \mathcal B_{v_*(\cdot+a_j,\cdot)}y_j
 =\partial_{\nu_{\Sph^{n-1}_+}}y_j
  -\frac{nc}{n-2}\tilde v_*(s+a_j)^{2/(n-2)}
     \widetilde y_j(s,\vartheta).
\]
The factor multiplying $\widetilde y_j$ after the minus sign
is the derivative of the boundary nonlinearity
$z\mapsto cz^{n/(n-2)}$ in~\eqref{eq:s3-normalized-model},
evaluated at the translated trace $z=\tilde v_*(s+a_j)$.
Since $v_*=v_{\eps_*}$ belongs to the smooth positive periodic
family of Section~\ref{subsec:companion-family}, the translated
traces are bounded away from zero and have uniform bounds
through order two. Hence the coefficient is bounded in
$C^{1,\gamma}(I\times\partial\Sph^{n-1}_+)$ independently of $j$. The local Schauder estimates, applied on
nested intervals as in Section~\ref{subsec:linearized-spaces},
give uniform $C^{2,\gamma}$ bounds on each compact product
with $\Sph^{n-1}_+$. A diagonal subsequence converges to $\Phi$
in $C^2_{\mathrm{loc}}(\R\times\Sph^{n-1}_+)$, including the
boundary. Set $v_\infty(s,\theta)=v_*(s+a_\infty,\theta)$.
The convergence of the translated coefficients and the
vanishing of the right-hand sides give the interior
and boundary equations of~\eqref{eq:jacobi-linearized-equation}
along $v_\infty$.
\begin{equation}\label{eq:proof-limit-jacobi}
 \begin{cases}
 \mathcal A\Phi=0,\\
 \mathcal B_{v_\infty}\Phi=0.
 \end{cases}
\end{equation}
Elliptic regularity makes $\Phi$ smooth up to the boundary,
and~\eqref{eq:proof-jacobi-growth} gives at most linear growth
in both directions of $s$. Thus
$\Phi\in\mathcal K_{\mathrm{pol}}(v_\infty)$.

We evaluate the Hamiltonian at $s=0$ in
\eqref{eq:proof-normalized-fields}. The definition
\eqref{eq:proof-hamiltonian} and Taylor's formula give
\[
 \begin{aligned}
 &|\Sph^{n-2}|\bigl(H(w)(t_j)-H(v_*)\bigr)\\
 &\qquad=D_j|\Sph^{n-2}|D H_{v_*(\cdot+a_j,\cdot)}[y_j](0)
     +O(D_j^2).
 \end{aligned}
\]
The remainder is $O(D_j^2)$ uniformly in $j$. Indeed,
\eqref{eq:proof-jacobi-growth} with $\xi=0$ bounds $y_j$ and
its first derivatives on $\{0\}\times\Sph^{n-1}_+$. Since
the interior integral in~\eqref{eq:proof-hamiltonian} is
quadratic in the function and its first derivatives, its
remainder after the linear term is bounded by $CD_j^2$.

For the boundary integral, fix $\vartheta\in\partial\Sph^{n-1}_+$.
By~\eqref{eq:proof-normalized-fields}, the trace along the
interpolation, for $0\leq\lambda\leq1$, is
\[
 \tilde v_*(a_j)+\lambda D_j\widetilde y_j(0,\vartheta)
 =(1-\lambda)\tilde v_*(a_j)
     +\lambda\tilde w(t_j,\vartheta).
\]
As $\lambda$ varies, this expression runs through the
segment joining $\tilde v_*(a_j)$ to $\tilde w(t_j,\vartheta)$.
By Proposition~\ref{prop:geometric-inputs} and the positivity
and periodicity of $v_*$ from Section~\ref{subsec:companion-family},
these segments lie in a fixed compact subinterval of $(0,\infty)$,
independently of $j$ and $\vartheta$.

On this interval the boundary integrand in
\eqref{eq:proof-hamiltonian}, as a function of the trace, has
bounded second derivative. Taylor's formula therefore bounds
the pointwise remainder by $CD_j^2|\widetilde y_j(0,\vartheta)|^2$.
Integrating over $\partial\Sph^{n-1}_+$ and using the bound
from~\eqref{eq:proof-jacobi-growth} gives the uniform
$O(D_j^2)$ estimate for the boundary integral.
By~\eqref{eq:proof-energy-limit}
and~\eqref{eq:proof-normalized-limit-energy}, the left-hand
side is $O(e^{-t_j})$. Dividing by $D_j$ and using
\eqref{eq:proof-window-ratios} yields
\[
 \left|D H_{v_*(\cdot+a_j,\cdot)}[y_j](0)\right|
 \leq C\left(\frac{e^{-t_j}}{D_j}+D_j\right)
 \leq C\left(\frac{e^{-j}}{j}+D_j\right)\longrightarrow0.
\]
Passing to the limit in~\eqref{eq:s5-full-first-variation}
at $s=0$ gives
\begin{equation}\label{eq:proof-zero-linear-energy}
 D H_{v_\infty}[\Phi](0)=0.
\end{equation}
The first variation in~\eqref{eq:s5-full-first-variation}
is constant in $s$ along a Jacobi field.
Theorem~\ref{s5:full-polynomial-kernel}, applied after
translating by $-a_\infty$, gives
$\Phi=b\,\partial_s v_\infty$ for some $b\in\R$.
Dividing the orthogonality condition
\eqref{eq:proof-phase-projection} at $t_j$ by $D_j$ and passing
to the limit gives
\[
 0=b\int_J\int_{\Sph^{n-1}_+}
       \chi(s)\bigl(\partial_s v_\infty(s,\theta)\bigr)^2
       \dd\theta\dd s .
\]
The integral is positive because $\chi$ is positive on an
interval containing a period and $v_\infty$ depends on $s$.
Hence $\Phi=0$. On the other hand,
\eqref{eq:proof-normalized-fields} and
\eqref{eq:proof-window-ratios} give
$\|y_j\|_{C^2(J\times\Sph^{n-1}_+)}>1/2$, and convergence
in this norm yields $\|\Phi\|_{C^2(J\times\Sph^{n-1}_+)}\geq1/2$.
This contradiction proves~\eqref{eq:proof-window}.

\medskip
\noindent\emph{Rate of convergence.}
We use the following argument here and for the cylindrical
solution. Suppose $d$ is continuous and nonnegative on
$[t_0,\infty)$, tends to zero, and satisfies, for $t\geq t_0+R$,
\begin{equation}\label{eq:proof-general-window}
 d(t)\leq\frac12\sup_{|s-t|\leq R}d(s)
                  +Ce^{-\sigma(t-R)}
\end{equation}
with $R,\sigma>0$. Choose
$0<\delta<\min\{\sigma,(\log 2)/R\}$. For a sufficiently
large $B_0>0$, the function $B(t)=B_0e^{-\delta(t-t_0)}$
satisfies $d(t)\leq B(t)$ on $[t_0,t_0+R]$ and
\[
 \left(1-\frac12e^{\delta R}\right)B(t)
                  \geq Ce^{-\sigma(t-R)}
\]
for $t\geq t_0+R$. Indeed, the factor on the left is positive,
and the ratio of the forcing term to $e^{-\delta(t-t_0)}$
decreases for $t\geq t_0+R$ because $\delta<\sigma$.

If $d-B$ has a positive supremum $M$, it attains this value at
some $t_*>t_0+R$, since $d-B$ is continuous and tends to zero.
For $|s-t_*|\leq R$, we have
$d(s)\leq M+B(s)\leq M+e^{\delta R}B(t_*)$. Substituting
this bound in~\eqref{eq:proof-general-window} gives
\[
 M+B(t_*)\leq\frac12 M+\frac12e^{\delta R}B(t_*)
                          +Ce^{-\sigma(t_*-R)}
             \leq\frac12 M+B(t_*),
\]
which contradicts $M>0$. Thus $d(t)\leq B(t)$.

For the noncylindrical solution, apply this argument to
\eqref{eq:proof-window} with $\sigma=1$ and $R$ as in that
estimate, taking $t_1$ as the initial time. It gives
$d(t)\leq Ce^{-\delta t}$ for each
$0<\delta<\min\{1,(\log 2)/R\}$. By
\eqref{eq:proof-phase-derivative}, the derivative of $a(t)-t$
is bounded in absolute value by $Ce^{-\delta t}$ and is
therefore integrable on $[t_0,\infty)$. There is a real
number $\tau_*$ such that $a(t)-t\to\tau_*$ as $t\to\infty$.
For all sufficiently large $t$,
\begin{equation}\label{eq:proof-phase-limit}
 |a(t)-t-\tau_*|\leq\int_t^\infty |a'(s)-1|\dd s
                       \leq Ce^{-\delta t}.
\end{equation}
By the smoothness and periodicity of $v_*$, the mean value
theorem gives
\[
 \|v_*(\cdot+a(t),\cdot)-v_*(t+\cdot+\tau_*,\cdot)
      \|_{C^2([-1,1]\times\Sph^{n-1}_+)}
 \leq C|a(t)-t-\tau_*|.
\]
Together with~\eqref{eq:proof-distance}, this proves
\eqref{eq:main-cylinder-estimate} when $\eps_*<1$.

\medskip
\noindent\emph{The cylindrical limiting solution.}
Suppose $\eps_*=1$, so that $v_*=C_n$. Fix $a_0>1$ and let
\[
 d(t)=\|w(t+\cdot,\cdot)-C_n\|_
             {C^2([-a_0,a_0]\times\Sph^{n-1}_+)}.
\]
Equation~\eqref{eq:proof-orbit-entry} gives $d(t)\to0$.
We prove that, for some $R,C>0$ and all sufficiently large $t$,
\begin{equation}\label{eq:proof-cylinder-window}
 d(t)\leq\frac12\sup_{|s-t|\leq R}d(s)+Ce^{-(t-R)/2}.
\end{equation}
By~\eqref{eq:s5-stationary-first-variation}, the first variation
of $H$ at $C_n$ vanishes. The argument therefore compares
the quadratic term in~\eqref{eq:s5-stationary-energy-taylor}
with the Hamiltonian error of order $e^{-t}$ in
\eqref{eq:proof-energy-limit}. This accounts for the
exponent $1/2$ in~\eqref{eq:proof-cylinder-window}.

Fix $t_0$ sufficiently large that $d$ is defined on
$[t_0,\infty)$. If~\eqref{eq:proof-cylinder-window} fails,
take $R=C=j$ and choose $t_j\geq t_0+2j$ such that, with
$D_j=\sup_{|s-t_j|\leq j}d(s)$,
\[
 d(t_j)>\frac12D_j+j e^{-(t_j-j)/2}.
\]
Then $D_j>0$ and $D_j\to0$, because $t_j-j\to\infty$. We
also have
\begin{equation}\label{eq:proof-cylinder-ratios}
 \frac{d(t_j)}{D_j}>\frac12,
 \qquad
 \frac{e^{-(t_j-j)}}{D_j^2}<\frac1{j^2}.
\end{equation}
Set $y_j(s,\theta)=(w(t_j+s,\theta)-C_n(\theta))/D_j$.
For each $\xi\in\R$ and all sufficiently large $j$,
\[
 \|y_j(\xi+\cdot,\cdot)\|_
          {C^2([-a_0,a_0]\times\Sph^{n-1}_+)}
     =\frac{d(t_j+\xi)}{D_j}\leq 1.
\]
The geometric errors divided by $D_j$ tend to zero on each
fixed interval in the local H\"older norms of
\eqref{eq:proof-defect-decay}, since
\[
 \frac{e^{-t_j}}{D_j}
     \leq\frac{e^{-j}}{j^2}D_j\longrightarrow0.
\]
The trace of $C_n$ equals $1$. Subtracting its equation
from~\eqref{eq:proof-defect-system} gives
\eqref{eq:proof-jacobi-bulk} and~\eqref{eq:proof-jacobi-boundary}
with $C_n$ as the solution at which the boundary condition
is linearized. The remainder in~\eqref{eq:proof-taylor},
with trace $1$ in place of the translated trace, tends to
zero in $C^{1,\gamma}$ on each fixed product with
$\partial\Sph^{n-1}_+$, since $y_j$ is bounded there in
$C^2$ and $D_j\to0$. The same local Schauder estimates as
in the noncylindrical case yield, after passing to a
subsequence, convergence in
$C^2_{\mathrm{loc}}(\R\times\Sph^{n-1}_+)$ to a Jacobi
field $\Phi$ along $C_n$. The bounds on the translated
intervals give $\Phi\in\mathcal K_{\mathrm{bd}}(C_n)$.
By~\eqref{eq:proof-cylinder-ratios} and convergence on
$[-a_0,a_0]\times\Sph^{n-1}_+$,
\begin{equation}\label{eq:proof-cylinder-nonzero}
 \|\Phi\|_{C^2([-a_0,a_0]\times\Sph^{n-1}_+)}\geq\frac12 .
\end{equation}

At $s=0$, equation~\eqref{eq:s5-stationary-energy-taylor} gives
\[
 |\Sph^{n-2}|\bigl(H(w)(t_j)-H(C_n)\bigr)
       =D_j^2\mathcal Q_c[y_j](0)+O(D_j^3).
\]
Since $C_n$ has trace $1$, the trace of $y_j$ at $s=0$ is
\[
 \widetilde y_j(0,\vartheta)
 =\frac{\tilde w(t_j,\vartheta)-1}{D_j}
\]
for $\vartheta\in\partial\Sph^{n-1}_+$. The $C^2$ estimate
for $y_j$ above, taken with $\xi=0$, gives
$|\widetilde y_j(0,\vartheta)|\leq d(t_j)/D_j\leq1$.
Thus the remainder in~\eqref{eq:s5-stationary-energy-taylor}
is uniform in $j$.
By~\eqref{eq:proof-normalized-limit-energy} with $\eps_*=1$,
the left-hand side is $O(e^{-t_j})$ by
\eqref{eq:proof-energy-limit}. Dividing by $D_j^2$ and using
\eqref{eq:proof-cylinder-ratios}, we obtain
\[
 |\mathcal Q_c[y_j](0)|
     \leq C\left(\frac{e^{-t_j}}{D_j^2}+D_j\right)
     \leq C\left(\frac{e^{-j}}{j^2}+D_j\right)\longrightarrow0 .
\]
The convergence of $y_j$ and its first derivatives on
$\{0\}\times\Sph^{n-1}_+$ gives $\mathcal Q_c[\Phi](0)=0$.
Since $\Phi$ belongs to $\mathcal K_{\mathrm{pol}}(C_n)$,
Proposition~\ref{s5:stationary-energy} implies $\Phi=0$,
contradicting~\eqref{eq:proof-cylinder-nonzero}. This proves
\eqref{eq:proof-cylinder-window}.

By~\eqref{eq:proof-cylinder-window}, $d$ satisfies
\eqref{eq:proof-general-window} with $\sigma=1/2$ for all
sufficiently large $t$. The preceding estimate gives
$d(t)\leq Ce^{-\delta t}$ for each
$0<\delta<\min\{1/2,(\log 2)/R\}$. Since $a_0>1$ and $C_n$ is
independent of $t$, this proves~\eqref{eq:main-cylinder-estimate}
when $\eps_*=1$, for every choice of $\tau_*$.

\medskip
\noindent\emph{Estimates in Fermi coordinates and uniqueness.}
Let $r=e^{-t}$ and let $U_*=U_{\eps_*,\tau_*}$ be defined
by~\eqref{eq:main-comparison}. We work on the closed
half-annulus
\[
 \{z\in\R^n:z_n\geq0,\ 1/2\leq|z|\leq2\}.
\]
For all sufficiently large $t$, we have $2r<\rho$.
Equations~\eqref{eq:main-change} and~\eqref{eq:main-comparison}
give, on this half-annulus,
\begin{equation}\label{eq:proof-annulus-return}
 \begin{aligned}
 r^{\frac{n-2}{2}}(U-U_*)(rz)
  ={}&A_\kappa^{-1}|z|^{-\frac{n-2}{2}}\biggl[
       w\left(t-\log|z|,\frac{z}{|z|}\right)\\
     &\qquad\qquad
       -v_{\eps_*}\left(t-\log|z|+\tau_*,\frac{z}{|z|}\right)
       \biggr].
 \end{aligned}
\end{equation}
The factor $A_\kappa^{-1}|z|^{-(n-2)/2}$ comes from
undoing the normalization in~\eqref{eq:main-change} at
$x=rz$. The map
\[
 z\longmapsto\left(-\log|z|,\frac{z}{|z|}\right)
\]
gives the cylindrical coordinates relative to $t$, since
$-\log|rz|=t-\log|z|$. Its image is contained in
$[-1,1]\times\Sph^{n-1}_+$ because
$\bigl|\log |z|\bigr|\leq\log 2<1$. The scalar factor
and the coordinate map are smooth away from $z=0$,
with bounded derivatives through order two on the
half-annulus. For the spherical component, these bounds
can be taken in smooth charts of the hemisphere, including
charts at the equator.

The expression in brackets in~\eqref{eq:proof-annulus-return}
is the difference estimated in~\eqref{eq:main-cylinder-estimate},
evaluated at $(-\log|z|,z/|z|)$. The product and chain rules
therefore bound the $C^2$ norm of the left-hand side on the
half-annulus by $Ce^{-\delta t}=Cr^\delta$. For a multi-index
$\alpha$ with $|\alpha|\leq2$, differentiation at fixed $r$
gives
\[
 \partial_z^\alpha
 \bigl[r^{\frac{n-2}{2}}(U-U_*)(rz)\bigr]
 =
 r^{\frac{n-2}{2}+|\alpha|}
 (\partial_x^\alpha(U-U_*))(rz).
\]
Evaluating at $|z|=1$ gives~\eqref{eq:main-coordinate-estimate}
with $|x|=r$. The estimates include $z_n=0$ because the
$C^2$ bound holds up to the equator.

For $\eps_*<1$, positivity and periodicity give a positive
minimum of $v_{\eps_*}$ on
$[0,T_{\eps_*}]\times\Sph^{n-1}_+$. For $\eps_*=1$, the
positive minimum is that of $C_n$ on $\Sph^{n-1}_+$.
Thus~\eqref{eq:main-comparison} gives
$U_*(x)\geq C^{-1}|x|^{-(n-2)/2}$. Dividing the estimate
for $j=0$ in~\eqref{eq:main-coordinate-estimate} by $U_*(x)$
gives
\[
 \left|\frac{U(x)}{U_*(x)}-1\right|\leq C|x|^\delta .
\]

To prove uniqueness, suppose that $U_{\eps_*',\tau_*'}$
satisfies~\eqref{eq:main-coordinate-estimate} with an
exponent $\delta'>0$. For $s\in[-1,1]$ and
$\theta\in\Sph^{n-1}_+$, put $z=e^{-s}\theta$. The points
$z$ range over the closed half-annulus
\[
 \{z\in\R^n:z_n\geq0,\ e^{-1}\leq|z|\leq e\}.
\]
For sufficiently small $r$, the assumed weighted bounds
apply at $x=rz$. For $|\alpha|\leq2$, they give
\[
 \begin{aligned}
 &\left|\partial_z^\alpha
 \bigl[r^{\frac{n-2}{2}}(U-U_{\eps_*',\tau_*'})(rz)\bigr]\right|\\
 &\qquad\leq Cr^{\delta'}
                 |z|^{\delta'-\frac{n-2}{2}-|\alpha|}
 \leq C'r^{\delta'} .
 \end{aligned}
\]
The last bound is uniform on this half-annulus because
$e^{-1}\leq|z|\leq e$. By~\eqref{eq:main-change}
and~\eqref{eq:main-comparison},
\[
 \begin{aligned}
 &w(t+s,\theta)-v_{\eps_*'}(t+s+\tau_*',\theta)\\
 &\qquad=A_\kappa e^{-\frac{n-2}{2}s}
          r^{\frac{n-2}{2}}
          (U-U_{\eps_*',\tau_*'})(re^{-s}\theta).
 \end{aligned}
\]
The factor $A_\kappa e^{-(n-2)s/2}$ and the map
$(s,\theta)\mapsto e^{-s}\theta$ have bounded derivatives
through order two on $[-1,1]\times\Sph^{n-1}_+$. The chain
rule thus gives
\[
 \bigl\|w(t+\cdot,\cdot)
       -v_{\eps_*'}(t+\cdot+\tau_*',\cdot)\bigr\|
       _{C^2([-1,1]\times\Sph^{n-1}_+)}
 \leq Ce^{-\delta't}.
\]
In particular, at $s=0$, the differences of the functions,
their first derivatives and their boundary traces tend uniformly
to zero. Substituting in~\eqref{eq:proof-hamiltonian}
and using the constancy of the Hamiltonian of $v_{\eps_*'}$
gives
\[
 H(w)(t)-H(v_{\eps_*'})\longrightarrow0
 \qquad\text{as }t\to\infty.
\]
By~\eqref{eq:proof-energy-limit},
$H(v_{\eps_*'})=H_0/|\Sph^{n-2}|$. Together with
\eqref{eq:proof-normalized-limit-energy} and
Theorem~\ref{s4:global-energy-coordinate}, this gives
$\eps_*'=\eps_*$.

Suppose $\eps_*<1$. The two cylindrical estimates imply
that the norm of the difference of
$v_{\eps_*}(t+\cdot+\tau_*,\cdot)$ and
$v_{\eps_*}(t+\cdot+\tau_*',\cdot)$ in
$C^2([-1,1]\times\Sph^{n-1}_+)$ tends to zero as
$t\to\infty$. This norm is $T_{\eps_*}$-periodic in $t$.
For any fixed $t$, it therefore has the same value at
$t+kT_{\eps_*}$ for every positive integer $k$. Taking
$k\to\infty$ shows that this value is zero. The two
translates consequently agree on $\R\times\Sph^{n-1}_+$,
so $\tau_*'-\tau_*$ is a period of $v_{\eps_*}$. Since
$T_{\eps_*}$ is its least positive period, we have
$\tau_*'-\tau_*\in T_{\eps_*}\mathbb Z$. If $\eps_*=1$,
every translate is $C_n$.
\end{proof}

\subsection{Conformal deformations of the flat model}
\label{subsec:deformed-boundary-models}

We compose tangential translations with the Kelvin transform
in~\eqref{eq:flat-kelvin} and identify the resulting term of order
$e^{-t}$ in Emden--Fowler coordinates. After multiplication by
$\sin^{(n-2)/2}\beta$, the term has the coefficient $z_-$ from
\eqref{eq:s5-conformal-amplitudes} in degree $1$. We then construct
a metric for which the deformed flat solutions leave a nonzero
term of order $|x|$ in the relative error in boundary Fermi
coordinates. For the corresponding construction in the problem
without boundary, see
\cite[Sections~2.1 and~5]{KMPS} and
\cite[Theorem~2 and Section~5]{Marques}.

\begin{proposition}\label{prop:deformed-boundary-models}
Let $U_*=U_{\eps,\tau}$ be the flat solution on
$(\R^n_+\setminus\{0\},\delta)$ defined by~\eqref{eq:main-comparison},
where $0<\eps\leq1$ and $\tau\in\R$. Fix
$a=(\overline a,0)$ with $\overline a\in\R^{n-1}$ and
choose $\rho_a>0$ such that $|a|\rho_a<1/2$. For
$x\in\R^n_+$ with $0<|x|\leq\rho_a$, define
\[
 \begin{aligned}
 U_{*,a}(x)
 &=\bigl(1-2a\cdot x+|a|^2|x|^2\bigr)^{-\frac{n-2}{2}}\\
 &\quad\times U_*\left(
       \frac{x-|x|^2a}{1-2a\cdot x+|a|^2|x|^2}
       \right).
 \end{aligned}
\]
Let $w_a$ be the Emden--Fowler transform of $U_{*,a}$
as in~\eqref{eq:main-change}, and let $h_-$ and $z_-$ be
the functions in~\eqref{eq:s5-conformal-amplitudes}
associated with $v_\eps(t+\tau,\theta)$. Write
$\theta=(\bar\theta,\theta_n)$. For $t>-\log\rho_a$ and
$\theta\in\Sph^{n-1}_+$, we have
\[
 w_a(t,\theta)=v_\eps(t+\tau,\theta)
       +e^{-t}(\overline a\cdot\bar\theta)h_-(t,\beta)
       +R_a(t,\theta).
\]
For each $J>0$, there is a constant $C_{\eps,\tau,a,J}$
such that, for all sufficiently large $t$,
\[
 \|R_a(t+\cdot,\cdot)\|_{C^2([-J,J]\times\Sph^{n-1}_+)}
 \leq C_{\eps,\tau,a,J}e^{-2t}.
\]
Since $\bar\theta=\sin\beta\,\vartheta$, multiplying the
term of order $e^{-t}$ by $\sin^{(n-2)/2}\beta$ gives
$z_-(t,\beta)(\overline a\cdot\vartheta)$.
\end{proposition}

\begin{proof}
By~\eqref{eq:main-change} and~\eqref{eq:main-comparison},
we have $w_a(t,\theta)=F(t,e^{-t},\theta)$, where
\[
 F(t,r,\theta)
 =|\theta-ra|^{-\frac{n-2}{2}}\,
 v_\eps\left(t+\tau+\log|\theta-ra|,
             \frac{\theta-ra}{|\theta-ra|}\right).
\]
We expand $F$ in $r$ at $r=0$, with $t$ and $\theta$ fixed.
Its value there is $v_\eps(t+\tau,\theta)$, and
\[
 \left.\partial_r\log|\theta-ra|\right|_{r=0}
 =-a\cdot\theta.
\]
Since $a_n=0$, the last coordinate of
$(\theta-ra)/|\theta-ra|$ is $\cos\beta/|\theta-ra|$.
For $0<\beta\leq L$, this gives
\[
 \left.\partial_r
 \arccos\left(\frac{\cos\beta}{|\theta-ra|}\right)
 \right|_{r=0}
 =-(a\cdot\theta)\cot\beta.
\]
By cylindrical symmetry, $v_\eps$ depends on the spherical
variable only through $\beta$. The product and chain rules
therefore give the following derivative at $r=0$:
\[
 (\overline a\cdot\bar\theta)
 \left(
 \frac{n-2}{2}v_\eps(t+\tau,\beta)
 -\partial_t v_\eps(t+\tau,\beta)
 -\cot\beta\,\partial_\beta v_\eps(t+\tau,\beta)
 \right).
\]
This equals $(\overline a\cdot\bar\theta)h_-(t,\beta)$
by~\eqref{eq:s5-conformal-amplitudes}. The derivative is
smooth in spherical charts, and $h_-$ is smooth up to
$\beta=0$ by Proposition~\ref{s5:positive-pair}. Hence
the identity holds at $\beta=0$ as well.

Taylor's formula gives the remainder
\[
 R_a(t,\theta)
 =e^{-2t}\int_0^1(1-\lambda)
       \partial_r^2 F(t,\lambda e^{-t},\theta)\dd\lambda .
\]
For $0\leq r\leq\rho_a$, we have $|\theta-ra|\geq1-|a|\rho_a>1/2$,
and $(\theta-ra)/|\theta-ra|$ belongs to $\Sph^{n-1}_+$ because
$a_n=0$. These facts and the uniform smooth bounds for $v_\eps$
from Section~\ref{subsec:companion-family} bound the derivatives of
$F$ through order four uniformly in smooth spherical charts,
including charts at the equator. Each of the first two derivatives
in $t$ of $\lambda e^{-t}$ is bounded uniformly in $\lambda\in[0,1]$
for sufficiently large $t$. Differentiating under the integral
and using the chain rule bounds $R_a$ and its derivatives
through order two by $Ce^{-2t}$.
Taking the supremum on $[-J,J]$ after translating $t$ gives the
$C^2$ bound for $R_a(t+\cdot,\cdot)$.

Finally, using $\bar\theta=\sin\beta\,\vartheta$, we get
\[
 \begin{aligned}
 &\sin^{(n-2)/2}\beta\,
   e^{-t}(\overline a\cdot\bar\theta)h_-(t,\beta)\\
 &\qquad
 =e^{-t}\sin^{n/2}\beta\,h_-(t,\beta)
       (\overline a\cdot\vartheta)
 =z_-(t,\beta)(\overline a\cdot\vartheta),
 \end{aligned}
\]
where the last equality is the definition of $z_-$
in~\eqref{eq:s5-conformal-amplitudes}.
\end{proof}

Marques~\cite[Theorem~2]{Marques} obtains a relative error of
order $O(|x|^\gamma)$ for some $\gamma>1$ by using deformed
Fowler solutions in the problem without boundary. The following
proposition establishes a first-order obstruction to such an
estimate in boundary Fermi coordinates. The proof uses the
boundary trace to determine the possible approximating solution
and then evaluates the error along the inward normal.

\begin{proposition}\label{prop:first-order-geometric-obstruction}
For every $n\geq3$ and $\kappa>0$, there exist $\rho>0$, a
smooth conformally flat metric $g$
(Definition~\ref{def:conformally-flat}) written in boundary
Fermi coordinates on $B_\rho^+$, and a positive solution $U$
satisfying the hypotheses of Theorem~\ref{thm:fixed-profile}.
Its Emden--Fowler transform $w$, defined by
\eqref{eq:main-change}, satisfies, for every $J>0$,
\[
 \|w(t+\cdot,\cdot)-C_n\|_{C^2([-J,J]\times\Sph^{n-1}_+)}
 \longrightarrow 0
 \qquad\text{as }t\to\infty.
\]

There are no parameters $0<\eps\leq1$, $\tau\in\R$,
$a=(\overline a,0)$ with $\overline a\in\R^{n-1}$,
and $\sigma>1$ for which
\[
 U(y)=\bigl(1+O(|y|^\sigma)\bigr)U_{*,a}(y)
\]
as $y\to0$ in $\R^n_+$. Here $U_{*,a}$ is obtained from
$U_*=U_{\eps,\tau}$ by Proposition~\ref{prop:deformed-boundary-models}.
\end{proposition}

\begin{proof}
Fix $\lambda\ne0$ and choose $\rho>0$ such that
$1+\lambda y_n>0$ for $y\in\R^n_+$ with $|y|\leq\rho$.
Consider the metric
\[
 g=\dd y_n^2+(1+\lambda y_n)^2
                  \sum_{i=1}^{n-1}\dd y_i^2.
\]
The induced boundary metric is Euclidean, and the normal
coordinate lines are unit-speed geodesics. Hence $y$ is a
boundary Fermi chart. The outward normal is
$-\partial_{y_n}$ at $y_n=0$, so $h_g(0)=-\lambda$.
The change of variables
\[
 x_n=\lambda^{-1}\log(1+\lambda y_n),
 \qquad \overline x=\overline y
\]
writes $g=e^{2\lambda x_n}\delta$. Define
\[
 U(y)=(1+\lambda y_n)^{-\frac{n-2}{2}}
     U_{1,0}\bigl(\overline y,\lambda^{-1}\log(1+\lambda y_n)\bigr),
\]
where $U_{1,0}$ is the cylindrical solution in
\eqref{eq:main-comparison}. Conformal covariance of $L_g$
and $B_g$ gives~\eqref{eq:main-geometric}. The trace is
$u(\overline y)=A_\kappa^{-1}|\overline y|^{-(n-2)/2}$, so
\eqref{eq:main-trace-bounds} holds.

By homogeneity of $U_{1,0}$, for $z$ in any fixed closed
half-annulus and sufficiently small $r>0$,
\[
 r^{\frac{n-2}{2}}U(rz)
   =(1+\lambda r z_n)^{-\frac{n-2}{2}}
     U_{1,0}\left(\overline z,
                  \frac{\log(1+\lambda r z_n)}{\lambda r}\right).
\]
The right-hand side tends to $U_{1,0}(z)$ in $C^2$ on
that half-annulus. By~\eqref{eq:main-change}, $w(t+\cdot,\cdot)$
therefore tends to $C_n$ in $C^2([-J,J]\times\Sph^{n-1}_+)$
for every $J>0$ as $t\to\infty$.

Suppose that the estimate in the statement holds for some
$\eps$, $\tau$, $a$ and $\sigma>1$. The expansion in
Proposition~\ref{prop:deformed-boundary-models} gives, uniformly
in $\vartheta\in\Sph^{n-2}$ as $r\to0$,
\[
 A_\kappa r^{\frac{n-2}{2}}U_{*,a}(r\vartheta,0)
        =\tilde v_\eps(-\log r+\tau)+O(r).
\]
The exact trace of $U$ now implies that
$\tilde v_\eps(t+\tau)\to1$ as $t\to\infty$. By periodicity,
$\tilde v_\eps\equiv1$, which gives $\eps=1$. On the flat face,
\[
 \bigl|\overline y-|\overline y|^2\overline a\bigr|^2
 =|\overline y|^2\bigl(1-2\overline a\cdot\overline y
                       +|\overline a|^2|\overline y|^2\bigr).
\]
Using the definition of $U_{*,a}$ and the trace of $U_{1,0}$, we get
\[
 \begin{aligned}
 \frac{U_{*,a}(\overline y,0)}{U_{1,0}(\overline y,0)}
 &=\bigl(1-2\overline a\cdot\overline y
                  +|\overline a|^2|\overline y|^2\bigr)^{-\frac{n-2}{4}}\\
 &=1+\frac{n-2}{2}\overline a\cdot\overline y+O(|\overline y|^2).
 \end{aligned}
\]
If $\overline a\ne0$, taking $\overline y$ in the direction of
$\overline a$ gives a nonzero linear term, which contradicts
$\sigma>1$. Thus $a=0$ and $U_{*,a}=U_{1,0}$. Along the inward
normal, homogeneity of $U_{1,0}$ gives
\[
 \begin{aligned}
 \frac{U(0,s)}{U_{1,0}(0,s)}
 &=\left(\frac{s}
 {(1+\lambda s)\lambda^{-1}\log(1+\lambda s)}\right)^{\frac{n-2}{2}}\\
 &=1-\frac{(n-2)\lambda}{4}s+O(s^2)
 \end{aligned}
\]
as $s\to0$ through positive values. The coefficient of $s$
is nonzero, contradicting the assumed estimate.
\end{proof}

If $\mathrm{II}_g(p)=0$, the expansion in
\eqref{eq:prelim-fermi-expansion} gives
$g_{ij}(x)-\delta_{ij}=O(|x|^2)$ and
$h_g(\overline x,0)=O(|\overline x|)$. The rescaling argument in the
proof of Proposition~\ref{prop:geometric-inputs} then gives
\eqref{eq:proof-defect-decay} with $e^{-2t_0}$ in place of
$e^{-t_0}$, as detailed in Section~\ref{subsec:app-geometric-defects}.
To obtain a higher-order expansion by this approach, the improved
bound must be combined with estimates for the Jacobi
equation in exponentially weighted spaces, separating the
conformal fields of degree $1$ from the other decaying fields.

\appendix
\section{Computations in cylindrical coordinates}
\label{app:cylindrical-computations}

Section~\ref{subsec:app-geometric-defects} expands the rescaled
Laplacian in boundary Fermi coordinates.
Section~\ref{subsec:app-period-shift} records the period shift of
the auxiliary function used in the pairing argument.
Section~\ref{subsec:app-degree-one} verifies the equations for the
fields of degree $1$ and estimates the endpoint term as $\beta\to0$
in Picone's identity.
Section~\ref{subsec:app-hamiltonian-variations} computes the
first and second derivatives of the Hamiltonian.

\subsection{Expansion of the rescaled metric Laplacian}
\label{subsec:app-geometric-defects}

We compute the linear term of the rescaled Laplacian in the
boundary Fermi chart and use it to improve the error estimates
in~\eqref{eq:proof-defect-decay} when $\mathrm{II}_g(p)=0$.

In the boundary Fermi chart of~\eqref{eq:prelim-fermi-map},
let $(g_{r_0})_{ij}(z)=g_{ij}(r_0z)$. By
\eqref{eq:proof-fermi}, we have $(g_{r_0})_{nn}=1$ and
$(g_{r_0})_{an}=(g_{r_0})_{na}=0$ for $1\leq a\leq n-1$.
Write $g_{r_0}^{ij}$ for the entries of the inverse matrix.
The identity
\[
 \sum_{i=1}^n g_{r_0}^{ji}(g_{r_0})_{in}
   =g_{r_0}^{jn}=\delta_{jn}
\]
gives $g_{r_0}^{nn}=1$ and $g_{r_0}^{an}=g_{r_0}^{na}=0$ for
$1\leq a\leq n-1$.

For $1\leq a,b\leq n-1$, the expansion
\eqref{eq:prelim-fermi-expansion} gives
\[
 (g_{r_0})_{ab}(z)
   =\delta_{ab}-2r_0z_n\mathrm{II}_{ab}(p)+O(r_0^2).
\]
The coefficient of $r_0$ is the derivative with respect to
$r_0$ at $r_0=0$. To compute the corresponding coefficient
for the inverse, differentiate
$\sum_{d=1}^{n-1}g_{r_0}^{ad}(g_{r_0})_{db}=\delta_{ab}$
at $r_0=0$. Both matrices equal the identity there, so
\[
 0=\left.\partial_{r_0}g_{r_0}^{ab}(z)\right|_{r_0=0}
      -2z_n\mathrm{II}_{ab}(p).
\]
Taylor's formula therefore gives
\[
 g_{r_0}^{ab}(z)
   =\delta_{ab}+2r_0z_n\mathrm{II}_{ab}(p)+O(r_0^2).
\]

In the determinant expansion, each product other than the
product of the diagonal entries contains at least two
off-diagonal entries. Each such entry is $O(r_0)$. Thus
\[
 \begin{aligned}
 \det g_{r_0}(z)
 &=\prod_{a=1}^{n-1}(g_{r_0})_{aa}(z)+O(r_0^2)\\
 &=1-2r_0z_n\sum_{a=1}^{n-1}\mathrm{II}_{aa}(p)+O(r_0^2)\\
 &=1-2r_0z_n(n-1)h_g(p)+O(r_0^2).
 \end{aligned}
\]
Since $\sqrt{\det g_{r_0}}+1=2+O(r_0)$, the identity
\[
 \sqrt{\det g_{r_0}}-1
 =\frac{\det g_{r_0}-1}{\sqrt{\det g_{r_0}}+1}
\]
gives
$\sqrt{\det g_{r_0}(z)}=1-r_0z_n(n-1)h_g(p)+O(r_0^2)$.
The remainders in these expansions satisfy the same
$O(r_0^2)$ estimate after any fixed number of $z$-derivatives
on a fixed half-annulus.

For a smooth function $\psi$ in the rescaled Fermi coordinates,
the product rule gives
\[
 \begin{aligned}
 \Delta_{g_{r_0}}\psi
 &=\frac{1}{\sqrt{\det g_{r_0}}}
   \sum_{i,j=1}^n
     \partial_{z_i}\bigl(\sqrt{\det g_{r_0}}\,
                         g_{r_0}^{ij}\partial_{z_j}\psi\bigr)\\
 &=\sum_{i,j=1}^n g_{r_0}^{ij}\partial_{z_i z_j}\psi\\
 &\quad+\sum_{j=1}^n\left(
   \sum_{i=1}^n\partial_{z_i}g_{r_0}^{ij}
   +\frac{\sum_{i=1}^n g_{r_0}^{ij}
               \partial_{z_i}\sqrt{\det g_{r_0}}}
          {\sqrt{\det g_{r_0}}}\right)\partial_{z_j}\psi .
 \end{aligned}
\]
For $i<n$, the coefficient of $r_0$ in $g_{r_0}^{ij}$ is
independent of $z_i$. For $i=n$, we have $g_{r_0}^{nj}=\delta_{nj}$.
Hence $\sum_{i=1}^n\partial_{z_i}g_{r_0}^{ij}=O(r_0^2)$.
The expansion of $\sqrt{\det g_{r_0}}$ gives
\[
 \partial_{z_i}\sqrt{\det g_{r_0}}
   =-r_0(n-1)h_g(p)\delta_{in}+O(r_0^2).
\]
Multiplication by $g_{r_0}^{ij}$, summation over $i$, and
division by $\sqrt{\det g_{r_0}}=1+O(r_0)$ give
\[
 \frac{\sum_{i=1}^n g_{r_0}^{ij}
             \partial_{z_i}\sqrt{\det g_{r_0}}}
      {\sqrt{\det g_{r_0}}}
   =-r_0(n-1)h_g(p)\delta_{jn}+O(r_0^2).
\]
It follows that
\[
 \begin{aligned}
 \Delta_{g_{r_0}}-\Delta_\delta
 ={}&2r_0z_n\sum_{a,b=1}^{n-1}
         \mathrm{II}_{ab}(p)\partial_{z_a z_b}\\
    &-r_0(n-1)h_g(p)\partial_{z_n}+O(r_0^2).
 \end{aligned}
\]
Here $O(r_0^2)$ is an operator of order at most two whose
coefficients satisfy the preceding estimates on fixed half-annuli.

The change of variables $x=r_0z$ gives
\[
 r_0^2(\Delta_gU)(r_0z)
   =\bigl[\Delta_{g_{r_0}}(U(r_0\,\cdot))\bigr](z).
\]
The same identity holds for $\Delta_\delta$. These are the
rescaling identities used to estimate $f_g$ in
Proposition~\ref{prop:geometric-inputs}. If $\mathrm{II}_g(p)=0$,
the coefficient of $r_0$ in $\Delta_{g_{r_0}}-\Delta_\delta$
vanishes, and $h_g(p)=0$. Fix $\ell>0$ and write $t=t_0+s$ and
$r_0=e^{-t_0}$. Smoothness gives $h_g(r_0e^{-s}\vartheta,0)=O(r_0)$ for $|s|\leq\ell$,
with the same bound for every fixed number of derivatives in
$(s,\vartheta)$. Together with the uniform rescaled estimates in
\eqref{eq:proof-scaled-regularity}, these bounds improve both
errors in~\eqref{eq:proof-defect-decay} to $O(r_0^2)$.

\subsection{The shift of the auxiliary function}
\label{subsec:app-period-shift}

Use $d$, $E_\eps$, $F_\eps$, and $G_\eps$ from the proof of
Proposition~\ref{s4:axial-polynomial-kernel}. The $T_\eps$-periodicity
of $E_\eps$ and $F_\eps$ gives, for every $t\in\R$ and
$\theta\in\Sph^{n-1}_+$,
\[
 \begin{aligned}
 &(t+T_\eps)E_\eps(t+T_\eps,\theta)
       +\frac d2(t+T_\eps)^2F_\eps(t+T_\eps,\theta)\\
 &\qquad-tE_\eps(t,\theta)-\frac d2t^2F_\eps(t,\theta)\\
 &\quad=T_\eps\bigl(E_\eps(t,\theta)+dtF_\eps(t,\theta)\bigr)
                   +\frac{dT_\eps^2}{2}F_\eps(t,\theta)\\
 &\quad=T_\eps G_\eps(t,\theta)
                   +\frac{dT_\eps^2}{2}F_\eps(t,\theta).
 \end{aligned}
\]
At $t=0$ and $\beta=L$, this difference equals $T_\eps>0$
by~\eqref{eq:s4-axial-parity-evaluations}.

\subsection{The equations in degree \texorpdfstring{$1$}{1} and Picone's identity}
\label{subsec:app-degree-one}

We verify the interior equation and boundary condition for $z_\pm$
from~\eqref{eq:s5-conformal-amplitudes} and compute the Taylor
expansion of $h_\pm$ at $\beta=0$. We then estimate the term at
$\beta=\delta$ in the integration by parts used in Picone's identity.

Let $v$ be the cylindrically symmetric solution considered in
Proposition~\ref{s5:positive-pair}. In expressions involving
$\Delta_{\Sph^{n+1}_+}$, a function of $\beta$ is regarded as
a cylindrically symmetric function on $\Sph^{n+1}_+$. For
$0<\beta<L$, its Laplacian on this hemisphere is
$\partial_\beta^2+n\cot\beta\,\partial_\beta$. The corresponding
formula on $\Sph^{n-1}_+$ and~\eqref{eq:s3-normalized-model} give
\[
 \left(-\partial_t^2-\Delta_{\Sph^{n+1}_+}
                +\frac{n^2}{4}-1\right)v
 =(n-2)v-2\cot\beta\,\partial_\beta v.
\]
Differentiating the equation for $v$ with respect to $\beta$
gives
\[
 \begin{aligned}
 &\left(-\partial_t^2-\partial_\beta^2
    -(n-2)\cot\beta\,\partial_\beta
    +\frac{(n-2)^2}{4}\right)\partial_\beta v\\
 &\qquad=-(n-2)\csc^2\beta\,\partial_\beta v .
 \end{aligned}
\]
Using $(\cot\beta)'=-\csc^2\beta$ in the product rule, we obtain
\[
 \begin{aligned}
 &\left(-\partial_t^2-\Delta_{\Sph^{n+1}_+}
       +\frac{n^2}{4}-1\right)
           \bigl(\cot\beta\,\partial_\beta v\bigr)\\
 &\qquad=2\partial_\beta^2 v+(n-2)\cot\beta\,\partial_\beta v .
 \end{aligned}
\]
Indeed, the coefficient of $\partial_\beta^2v$ is
$2(\csc^2\beta-\cot^2\beta)=2$. The terms with
$\cot\beta\,\csc^2\beta\,\partial_\beta v$ have coefficients
$-(n-2)$, $-2$, and $n$, and their sum is zero. It follows
from the equation for $v$ that
\[
 \begin{aligned}
 &\left(-\partial_t^2-\Delta_{\Sph^{n+1}_+}
       +\frac{n^2}{4}-1\right)
       \left(\frac{n-2}{2}v-\cot\beta\,\partial_\beta v\right)\\
 &\qquad=\frac{(n-2)^2}{2}v
          -2(n-2)\cot\beta\,\partial_\beta v
          -2\partial_\beta^2v
 =2\partial_t^2 v .
 \end{aligned}
\]
Differentiating the identity for the operator applied to $v$
with respect to $t$ gives
\[
 \begin{aligned}
 &\left(-\partial_t^2-\Delta_{\Sph^{n+1}_+}
                +\frac{n^2}{4}-1\right)\partial_t v\\
 &\qquad=2\partial_t
     \left(\frac{n-2}{2}v-\cot\beta\,\partial_\beta v\right).
 \end{aligned}
\]
Adding or subtracting this identity from the preceding formula
and using the definition of $h_\pm$ in
\eqref{eq:s5-conformal-amplitudes} gives
\eqref{eq:s5-lifted-equation} for $0<\beta<L$.

The operator in~\eqref{eq:s5-mode-operator} with $\ell=1$
has potential $n(n-2)\csc^2\beta/4$. The identities
\[
 \begin{aligned}
 \frac{\partial_\beta(\sin^{n/2}\beta)}{\sin^{n/2}\beta}
   &=\frac n2\cot\beta,\\
 \frac{\partial_\beta^2(\sin^{n/2}\beta)}{\sin^{n/2}\beta}
   &=\frac{n(n-2)}4\csc^2\beta-\frac{n^2}{4}
 \end{aligned}
\]
show how this potential cancels when $\mathcal L_1$ acts on
$\sin^{n/2}\beta\,h_\pm$. More precisely,
\[
 \begin{aligned}
 &(\mathcal L_1-1)(\sin^{n/2}\beta\,h_\pm)\\
 &\qquad=\sin^{n/2}\beta
      \left(-\partial_t^2-\Delta_{\Sph^{n+1}_+}
                   +\frac{n^2}{4}-1\right)h_\pm
 =\pm2\sin^{n/2}\beta\,\partial_t h_\pm .
 \end{aligned}
\]
This is~\eqref{eq:s5-conformal-intertwining}. By the product
rule, applying $\mathcal L_1$ to $z_\pm$ gives
\[
 \mathcal L_1z_\pm
   =e^{\pm t}(\mathcal L_1-1\mp2\partial_t)
      (\sin^{n/2}\beta\,h_\pm)=0.
\]

At $\beta=L$, the boundary condition in~\eqref{eq:s3-normalized-model}
gives $\partial_\beta v(t,L)=c\tilde v(t)^{n/(n-2)}$.
Differentiating the formula for $h_\pm$ in
\eqref{eq:s5-conformal-amplitudes} with respect to $\beta$
and evaluating at $L$ gives
\[
 \begin{aligned}
 h_\pm(t,L)
 &=\frac{n-2}{2}\tilde v(t)\pm\tilde v'(t),\\
 \partial_\beta h_\pm(t,L)
 &=\frac n2\,\partial_\beta v(t,L)
      \pm\partial_t\partial_\beta v(t,L)\\
 &=\frac{nc}{n-2}\tilde v(t)^{2/(n-2)}
      \left(\frac{n-2}{2}\tilde v(t)\pm\tilde v'(t)\right).
 \end{aligned}
\]
The term $+\partial_\beta v(t,L)$ in the second line comes
from differentiating $-\cot\beta\,\partial_\beta v$.
Since $\sin^{n/2}L=1$ and
$\partial_\beta(\sin^{n/2}\beta)|_{\beta=L}=0$, this gives
the condition for $z_\pm$ in~\eqref{eq:s5-mode-operator}.

The even extensions of $v$ and $h_\pm$ near $\beta=0$ are
described in Section~\ref{subsec:s5-polynomial}.
Taylor's formula for the even extension of $v$ gives
\[
 \begin{aligned}
 v(t,\beta)
   &=v(t,0)+\frac{\beta^2}{2}\partial_\beta^2v(t,0)
             +O(\beta^4),\\
 \partial_t v(t,\beta)
   &=\partial_t v(t,0)
       +\frac{\beta^2}{2}\partial_t\partial_\beta^2v(t,0)
       +O(\beta^4),\\
 \partial_\beta v(t,\beta)
   &=\beta\,\partial_\beta^2v(t,0)
      +\frac{\beta^3}{6}\partial_\beta^4v(t,0)+O(\beta^5).
 \end{aligned}
\]
Using $\cot\beta=\beta^{-1}-\beta/3+O(\beta^3)$ in the last
line gives
\[
 \cot\beta\,\partial_\beta v(t,\beta)
 =\partial_\beta^2v(t,0)
  +\beta^2\left(\frac16\partial_\beta^4v(t,0)
              -\frac13\partial_\beta^2v(t,0)\right)
  +O(\beta^4).
\]
Substitution in~\eqref{eq:s5-conformal-amplitudes} yields
\[
 \begin{aligned}
 h_\pm(t,\beta)
 &=\left[\frac{n-2}{2}v-\partial_\beta^2v
                     \pm\partial_t v\right](t,0)\\
 &\quad+\beta^2\left[
       \left(\frac{n-2}{4}+\frac13\right)\partial_\beta^2v
       -\frac16\partial_\beta^4v
       \pm\frac12\partial_t\partial_\beta^2v
                    \right](t,0)+O(\beta^4).
 \end{aligned}
\]
More generally, for each integer $k\geq0$, evenness gives
\[
 h_\pm(t,\beta)
   =\sum_{j=0}^k
      \frac{\partial_\beta^{2j}h_\pm(t,0)}{(2j)!}\beta^{2j}
      +O(\beta^{2k+2}).
\]
In geodesic normal coordinates centered at $(0,\ldots,0,1)$
on $\Sph^{n+1}_+$, $\beta$ is the norm of the coordinate
vector. Each power $\beta^{2j}$ is therefore a polynomial
in these coordinates. By Taylor's formula and the chain
rule, the derivatives of the remainder of order
$m\leq 2k+1$ in these coordinates are
$O(\beta^{2k+2-m})$. Since $k$ is arbitrary, the radial
extension is smooth at the center of the chart, and
\eqref{eq:s5-lifted-equation} holds there by continuity.

In the integration by parts for~\eqref{eq:s5-picone} over
$I\times[\delta,L]$, the lower endpoint contributes
\[
 -\int_I
 \frac{\partial_\beta\phi_0(t,\delta)}{\phi_0(t,\delta)}
 Z(t,\delta)^2\,\dd t,
 \qquad \phi_0=\sqrt{z_+z_-}.
\]
For the smooth approximations in the proof of
Proposition~\ref{s5:nonaxial-gap}, this term vanishes when
$\delta$ is sufficiently small.
For the coefficient $Z_1$ of a Jacobi field, smoothness at
$\beta=0$ and the zero mean of a degree-$1$ spherical harmonic
give $y_1(t,\beta)=O(\beta)$, as
in Section~\ref{subsec:s5-polynomial}. Multiplication by
$\sin^{(n-2)/2}\beta$ gives $Z_1(t,\beta)=O(\beta^{n/2})$.
The estimate is uniform for $t\in[0,T]$ in the periodic case.
For $I=\R$, apply the same estimate to $Z_1$ after multiplication
by a fixed smooth cutoff in $t$.

By~\eqref{eq:s5-conformal-amplitudes},
$\phi_0=\sin^{n/2}\beta\sqrt{h_+h_-}$. Differentiation gives
\[
 \frac{\partial_\beta\phi_0}{\phi_0}
  =\frac n2\cot\beta
   +\frac{\partial_\beta h_+}{2h_+}
   +\frac{\partial_\beta h_-}{2h_-}.
\]
The Taylor expansions above give
$\partial_\beta h_\pm(t,\beta)=O(\beta)$, and
\eqref{eq:s5-lifted-positive-pair} bounds $h_\pm$ below by
a positive constant. It follows that
$\partial_\beta\phi_0/\phi_0=O(\beta^{-1})$.
Thus for $Z=Z_1$ in the periodic case, or for the cutoff
of $Z_1$ when $I=\R$, the endpoint term satisfies
\[
 \left|\int_I \frac{\partial_\beta\phi_0(t,\delta)}
                        {\phi_0(t,\delta)}
                    Z(t,\delta)^2\,\dd t\right|
 \leq C\delta^{n-1}\longrightarrow0 .
\]

\subsection{Differentiation of the Hamiltonian}
\label{subsec:app-hamiltonian-variations}

We record the first variation of the Hamiltonian
\eqref{eq:proof-hamiltonian} in an arbitrary direction and its
second derivative at the cylindrical solution $C_n$.

Fix $t$ and let $v$ and $y$ be smooth near
$\{t\}\times\Sph^{n-1}_+$, with $\tilde v(t,\vartheta)>0$ on
$\partial\Sph^{n-1}_+$. By compactness of the boundary section,
$\tilde v+s\tilde y$ remains positive there for sufficiently small
$|s|$. Differentiating $s\mapsto H(v+s y)(t)$ at $s=0$ gives
\[
 \begin{aligned}
 |\Sph^{n-2}|D H_v[y](t)
 &=\int_{\Sph^{n-1}_+}
    \left(v_t y_t-
      \langle\nabla_{\Sph^{n-1}_+}v,\nabla_{\Sph^{n-1}_+}y\rangle
      -\frac{(n-2)^2}{4}vy\right)\dd\theta\\
 &\quad+c\int_{\partial\Sph^{n-1}_+}
           \tilde v^{n/(n-2)}\tilde y\dd\sigma .
 \end{aligned}
\]
The derivative of the boundary integral has coefficient $c$
because the exponent $2(n-1)/(n-2)$ cancels the factor
$(n-2)/(2(n-1))$ in~\eqref{eq:proof-hamiltonian}.

Apply this formula to $w$ from~\eqref{eq:main-change}.
Integrating by parts on $\Sph^{n-1}_+$ gives the boundary
term $-\partial_{\nu_{\Sph^{n-1}_+}}w\,\tilde y$. Adding the
derivative of the boundary integral and using
\eqref{eq:proof-defect-system} gives $-f_{\partial,g}\tilde y$.
In the interior, the same equation gives
$\Delta_{\Sph^{n-1}_+}w-(n-2)^2w/4=f_g-w_{tt}$. Thus
\[
 \begin{aligned}
 |\Sph^{n-2}|D H_w[y](t)
 &=\int_{\Sph^{n-1}_+}(w_t y_t-w_{tt}y)\dd\theta\\
 &\quad+\int_{\Sph^{n-1}_+} f_g y\dd\theta
       -\int_{\partial\Sph^{n-1}_+} f_{\partial,g}\tilde y\dd\sigma .
 \end{aligned}
\]
Taking $y=w_t$ recovers the derivative of $H(w)(t)$ used to prove
\eqref{eq:proof-energy-limit}. When the two errors vanish, the
identity reduces to~\eqref{eq:s5-full-first-variation}.

At $C_n$, the derivative vanishes by
\eqref{eq:s5-full-first-variation} and $\partial_t C_n=0$.
Since $C_n$ has trace $1$, differentiating twice in the
directions $y_1$ and $y_2$ gives
\[
 \begin{aligned}
 &|\Sph^{n-2}|D^2 H_{C_n}[y_1,y_2](t)\\
 &\quad=\int_{\Sph^{n-1}_+}\left(
    \partial_t y_1\,\partial_t y_2
    -\langle\nabla_{\Sph^{n-1}_+}y_1,\nabla_{\Sph^{n-1}_+}y_2\rangle
    -\frac{(n-2)^2}{4}y_1y_2\right)\dd\theta\\
 &\qquad+\frac{nc}{n-2}
         \int_{\partial\Sph^{n-1}_+}\tilde y_1\tilde y_2\dd\sigma .
 \end{aligned}
\]
For $y_1=y_2=y$, division by $2$ gives the quadratic
term $\mathcal Q_c[y](t)$ in~\eqref{eq:s5-stationary-quadratic-form}.

\end{document}